\documentclass[10pt]{article}

\usepackage{amsmath}
\usepackage{amsthm}
\usepackage{amssymb}
\usepackage{amsfonts}
\usepackage{graphicx}
\usepackage{verbatim}
\usepackage{mathrsfs}
\usepackage{fancyhdr}
\usepackage{latexsym}
\usepackage{graphicx}
\usepackage{dsfont}
\usepackage{color}
\usepackage[hidelinks]{hyperref}
\usepackage{bm}
\usepackage{drawmatrix}
\usepackage{enumitem}
\usepackage[percent]{overpic}
\usepackage{subcaption}

\newcommand{\R}{\mathbb R}
\newcommand{\N}{\mathbb N}
\newcommand{\Z}{\mathbb Z}
\newcommand{\C}{\mathbb C}
\newcommand{\E}{\mathbb E}
\newcommand{\V}{\mathbb V}
\newcommand{\1}{\mathds 1}

\newcommand{\brho}{\bar\varrho}
\newcommand{\bb}{\bar{b}}
\DeclareMathOperator{\atanh}{atanh}

\newcommand{\cC}{{\mathcal C}}  

\def\txtd{{\textnormal{d}}}
\def\txte{{\textnormal{e}}}
\def\txti{{\textnormal{i}}}

\def\txts{{\textnormal{s}}}
\def\txtu{{\textnormal{u}}}
\def\txtD{{\textnormal{D}}}

\theoremstyle{plain}
\newtheorem{theorem}{Theorem}[section]

\newtheorem{lemma}[theorem]{Lemma}

\newtheorem{remark}[theorem]{Remark}
\newtheorem{definition}[theorem]{Definition}

\allowdisplaybreaks

\title{Computing invariant sets of random differential equations using polynomial chaos}

\author{Maxime Breden \thanks{Technical University of Munich, 
Faculty of Mathematics, Research Unit ``Multiscale and Stochastic 
Dynamics", 85748 Garching b. M\"unchen, Germany. \texttt{maxime.breden@tum.de}} 
\and Christian Kuehn \thanks{Technical University of Munich, Faculty 
of Mathematics, Research Unit ``Multiscale and Stochastic Dynamics", 
85748 Garching b. M\"unchen, Germany. \texttt{ckuehn@ma.tum.de}}}

\begin{document}

\maketitle

\begin{abstract}
Differential equations with random parameters have gained significant prominence in recent years due to their importance in mathematical modelling and data assimilation. In many cases, random ordinary differential equations (RODEs) are studied by using Monte-Carlo methods or by direct numerical simulation techniques using polynomial chaos (PC), i.e., by a series expansion of the random parameters in combination with forward integration. Here we take a dynamical systems viewpoint and focus on the \emph{invariant sets} of differential equations such as steady states, stable/unstable manifolds, periodic orbits, and heteroclinic orbits. We employ PC to compute representations of all these different types of invariant sets for RODEs. This allows us to obtain fast sampling, geometric visualization of distributional properties of invariants sets, and uncertainty quantification of dynamical output such as periods or locations of orbits. We apply our techniques to a predator-prey model, where we compute steady states and stable/unstable manifolds. We also include several benchmarks to illustrate the numerical efficiency of adaptively chosen PC depending upon the random input. Then we employ the methods for the Lorenz system, obtaining computational PC representations of periodic orbits, stable/unstable manifolds and heteroclinic orbits.
\end{abstract}

\begin{center}
\begin{small}
\end{small}
\end{center}

\tableofcontents

\section{Introduction}

In this work, we study random nonlinear dynamical systems. More precisely, we focus on nonlinear random ordinary differential equations (RODEs) of the form
\begin{equation}
\label{eq:ODE_rand_coef}
\frac{\txtd x}{\txtd t}=\dot{x}=f(x,p),\qquad x=x(t)\in\R^m,
\end{equation}
where $f:\R^m\times\R^d \to \R^m$ is a smooth vector field and $p=p(\omega)\in\R^d$ denotes random parameters with given distributions. From the viewpoint of applications, it is frequently natural to assume that the parameters $p$ are only known from measurements, which naturally carry an associated probability distribution. Then the challenge is to quantify the uncertainty in the ``output'' of the RODE~\eqref{eq:ODE_rand_coef} based upon the random input. Yet, we are still far away to fully understand the nonlinear dynamics of such RODEs. In classical uncertainty quantification problems~\cite{Grigoriu2,LeMKni10}, one is often interested~\cite{LucorKarniadakis,VenturiWanKarniadakis} in the moments as an output of the solution $\mathbb{E}[x(t)^k]$ for $k\geq 1$, where $\mathbb{E}$ denotes the expectation. Here we take a dynamical systems perspective focusing on the invariant sets (e.g.~steady states, periodic orbits, invariant manifolds, connecting orbits, etc.) of~\eqref{eq:ODE_rand_coef}, and especially in understanding their dependence on the noisy parameters $p$. In this paper, we study the invariant sets from a numerical viewpoint. However, the framework that we develop is well suited to the usage of rigorous numerics~\cite{Tuc11}, and in particular of a-posteriori validation techniques~\cite{BerLes15}, which could applied to obtain rigorous results about these stochastic invariant sets. This idea will be presented in a forthcoming work.

Let us start by mentioning that there exists many well developed techniques to numerically study invariant sets of \emph{deterministic} ODEs. Therefore, a natural way of studying the invariant sets of~\eqref{eq:ODE_rand_coef} would be to use a Monte-Carlo type approach: consider a large sample of values $p_i$ taken according to the distribution of $p$, and for each $i$ study the invariant sets of the deterministic ODE $\dot{x}=f(x,p_i)$. However, this approach is known to be very costly, because it requires a large sample to accurately represent the statistics of the invariant sets~\cite{Fis13}. In this work we make use of a different technique, namely polynomial chaos (PC) expansions~\cite{Wie38,GhaSpa91}, to accurately compute invariant sets of~\eqref{eq:ODE_rand_coef}. Roughly speaking, we view each invariant set of~\eqref{eq:ODE_rand_coef} as a curve parametrized by $p$ (or as a manifold if $p$ is more than one-dimensional), and compute such parameterization explicitly via a PC expansion. This can be thought of as a parameter continuation in $p$, but in an astute way that allows us to obtain not only the geometrical object (i.e. the curve/manifold of invariant sets), but also statistical properties of this object (e.g. mean position, variance, ...). Furthermore, our techniques naturally extend to numerical continuation algorithms~\cite{KrauskopfOsingaGalan-Vioque,KuehnSDEcont1} if the probability distributions of $p$ contain further parameters, which is a direction we will pursue in future work.

Let $X=X(p)$ be an invariant set of~\eqref{eq:ODE_rand_coef}, i.e., trajectories starting inside $X$ remain in $X$ for all $t\in\R$. In this work, we focus on the cases where $X(p)$ represents a steady state, a stable/unstable manifold, a periodic orbit or a heteroclinic orbit, which are among the most important objects in nonlinear dynamics~\cite{Kuz13,GH}. One key observation is that PC expansions can be used to unify the computational framework for invariant sets and that they provide a natural deterministic structure to view the dynamics of nonlinear random ODEs. Our goal is to find a series expansion of $X$ as a function of $p$:
\begin{equation}
\label{eq:expansion1}
X(p)=\sum_{n\in\N^d} X_n \phi_n(p).
\end{equation}
Note that the coefficients $X_n$ have to be found as solutions of suitable nonlinear problems to correctly represent the different invariant sets; we develop the details for each type of invariant set in this work. In practice, the choice of the expansion basis $\left(\phi_n\right)_{n\in\N^d}$ is of course also crucial, as it determines the quality of the approximation, or more precisely the number $N$ of coefficients required to obtain a good enough approximation. For the rest of this discussion, we assume for simplicity that $p$ is one-dimensional. 

If $f$ is analytic with respect to $p$, then we can expect $X$ to also be analytic as a function of $p$, at least around values of $p$ that are not bifurcation values. Therefore, in many cases one could think about writing $X$ as a Taylor series by taking $\phi_n(p)=p^n$. However, it is well known in approximation theory that faster convergence can be achieved by considering instead Chebyshev series or Legendre series (i.e.~taking $\phi_n(p)=T_n(p)$ or $\phi_n(p)=L_n(p)$, where $T_n$ and $L_n$ respectively denote Cheybshev and Legendre orthogonal polynomials), mainly because these expansions are less sensitive to potential poles of $X$ in the complex plane; see e.g.~\cite{Tre13}. Chebyshev or Legendre expansions also have the advantage of being convergent even when $X$ is only of class $\cC^k$, where the coefficients $X_n$ decay at an algebraic rate, rather than geometric in the analytic case. If $p$ is deterministic, then looking at the decay rates of the coefficients $X_n$ is a relevant benchmark, because it is related to the error in $\cC^0$-norm. For instance, if $X^N$ is the truncated series given by
\begin{equation*}
X^N(p)=\sum_{n=0}^N X_n \phi_n(p),
\end{equation*}
then for a Taylor expansion, a Chebyshev expansion or a Legendre expansion alike one has
\begin{equation*}
\sup_{p\in[-1,1]} \left\Vert X(p) - X^N(p) \right\Vert \leq \sum_{n=N+1}^\infty \left\Vert X_n\right\Vert,
\end{equation*}
because for each of these choices one has $\sup_{p\in[-1,1]} \left\vert \phi_n(p)\right\vert \leq 1$. However, when $p=p(\omega)$ is a random variable, one is more interested in controlling different quantities such as
\begin{equation*}
\left\Vert \E \left(X(p)\right) - \E \left(X^N(p)\right)\right\Vert \qquad \text{or} \qquad \left\Vert \V \left(X(p)\right) - \V \left(X^N(p)\right)\right\Vert,
\end{equation*}
where $\V$ denotes the variance. To minimize the error in the moments, the choice of the expansion is critical not only because it influences the decay of the coefficients, but also because the $\phi_n$ themselves appear in the error term. For instance, one has
\begin{equation*}
\left\Vert \E \left(X(p)\right) - \E \left(X^N(p)\right)\right\Vert \leq \sum_{n=N+1}^\infty \left\Vert X_n\right\Vert \left\vert \E \left(\phi_n(p)\right)\right\vert.
\end{equation*}
Therefore, two different expansions leading to the same decay of the coefficients may not lead to an error of the same order for the expectation or of the variance. Of course, one could change these weights $\left\vert \E \left(\phi_n(p)\right)\right\vert$ by rescaling the $\phi_n$, but this only shifts the problem because the rescaling would then affect the decay of the coefficients.

In this work, we use the PC paradigm to minimize such quantities, by choosing an expansion basis $\phi_n$ that is adapted to the distribution of the noisy parameter $p$. PC expansions have become a very important tool in uncertainty quantification in the last decades, and a review of its many applications is far beyond the scope of the present work. We instead refer the interested reader to the survey~\cite{Xiu09} and the book~\cite{LeMKni10}.

The paper is structured as follows. First, we review the PC methodology in Section~\ref{sec:PC}. In Section~\ref{sec:inv_sets} we introduce some classical techniques to numerically study periodic orbits, invariant manifold and connecting orbits of deterministic ODEs. Then we proceed to the main contributions of our work. We explain, how to combine the numerical study of invariant sets with PC expansions to study the dynamics of nonlinear RODEs. In Section~\ref{sec:LV} we focus on a relatively simple example, where many quantities can also be computed analytically, which allows us to benchmark our numerical computations in the context of steady states. Yet, we also go beyond explicit structures and compute the distribution of stable/unstable manifolds using the ``parameterization method'' in combination with PC. Furthermore, we explain the implications of the random structure of the invariant sets and how to gain information from the moments of the invariants sets very efficiently. In Section~\ref{sec:Lorenz} we then turn our attention to another example, where analytic computations are no longer available, and showcase the potential of our approach on the Lorenz system. For this system, we compute periodic orbits using a combined Fourier and PC ansatz, we compute again invariant stable/unstable manifolds, as well as heteroclinic connections. 

\section{A quick review on PC}
\label{sec:PC}

Let $\rho:\R\to\R$ be probability distribution function (PDF) having finite moments, i.e.
\begin{equation*}
 \int_\R s^n\rho(s) ~\txtd s <\infty \qquad \forall~n\in\N.
\end{equation*} 

Given normalization constants $\left(h_n\right)_{n\in\N}$, $h_n>0$ for all $n\in\N$, there exists a unique family of orthogonal polynomials  $\left(\phi_n\right)_{n\in\N}$ associated to the weight $\rho$, i.e. satisfying
\begin{equation*}
\left\langle\phi_{n_1},\phi_{n_2}\right\rangle_\rho := \int_\R \phi_{n_1}(s)\phi_{n_2}(s)\rho(s) ~\txtd s = h_{n_1}\delta_{{n_1},{n_2}},\qquad \forall~{n_1},{n_2}\in\N.
\end{equation*}

The most classical examples are:
\begin{itemize}
\item The Hermite polynomials $H_n$, which correspond to $\rho(s)=\frac{1}{\sqrt{2\pi}}\txte^{-\frac{s^2}{2}}$ and $h_n=n!$;
\item The Laguerre polynomials $L_n$, which correspond to $\rho(s)=\1_{s\in[0,+\infty)} \txte^{-s}$ and $h_n=1$;
\item The Jacobi polynomials $P_n^{\alpha,\beta}$, $\alpha,\beta>-1$, which correspond to $\rho(s)=\1_{s\in(-1,1)} \frac{(1-s)^\alpha(1+s)^\beta}{2^{\alpha+\beta+1}B(\alpha+1,\beta+1)}$ and $h_n=\frac{B(n+\alpha+1,n+\beta+1)}{(2n+\alpha+\beta+1)B(\alpha+1,\beta+1)B(n+1,n+\alpha+\beta+1)}$, where $B(x,y)=\frac{\Gamma(x)\Gamma(y)}{\Gamma(x+y)}$ is the Euler beta function.
\end{itemize}
Within the class of the Jacobi polynomials, we list a few remarkable cases (sometimes having different normalizations) that we make use of in this work:
\begin{itemize}
\item The Legendre polynomials $P_n$, which correspond to $\rho(s)=\frac{1}{2}\1_{s\in(-1,1)}$ and $h_n=\frac{1}{2n+1}$;
\item The Chebyshev polynomials of the first kind $T_n$, which correspond to $\rho(s)=\1_{s\in(-1,1)}\frac{1}{\pi\sqrt{1-s^2}}$ and $h_0=1$, $h_n=\frac{1}{2}$, for $n\geq 1$;
\item The Chebyshev polynomials of the second kind $U_n$, which correspond to $\rho(s)=\1_{s\in(-1,1)}\frac{2}{\pi}\sqrt{1-s^2}$ and $h_n=1$;
\item The Gegenbauer or ultraspherical polynomials $C_n^{\mu}$, $\mu>-\frac{1}{2}$, $\mu\neq 0$, which correspond to $\rho(s)=\1_{s\in(-1,1)}\frac{2^{2\mu-1}\mu B(\mu,\mu)}{\pi}(1-s^2)^{\mu-\frac{1}{2}}$ and $h_n=\frac{\mu}{n(n+\mu)B(n,2\mu)}$.
\end{itemize}
For a more complete description of PC choices and their relations to the Askey scheme, see~\cite{XiuKar02}.

\begin{figure}[htpb]
\centering
\includegraphics[width=\linewidth]{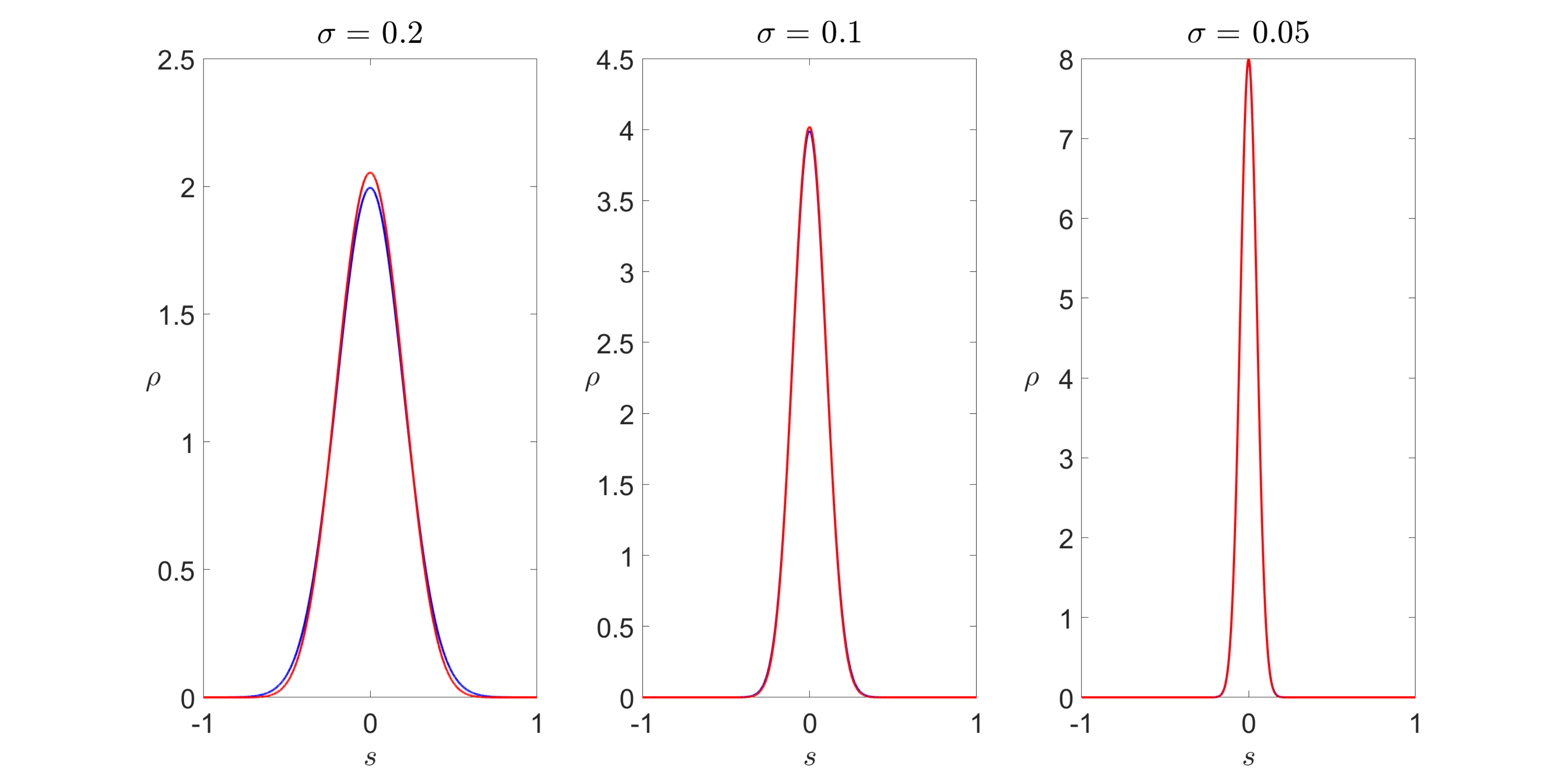}
\caption{Comparison between the Gaussian probability distribution $\rho(s)=\frac{1}{\sigma\sqrt{2\pi}}exp(-\frac{s^2}{2\sigma^2})$, in blue, and the symmetric Gamma probability distribution $\rho(s)=\1_{s\in(-1,1)}\frac{2^{2\mu-1}\mu B(\mu,\mu)}{\pi}(1-s^2)^{\mu-\frac{1}{2}}$, in red, for several values of $\sigma$, and with $\mu=\frac{1}{2}\left(1+\frac{1}{\sigma^2}\right)$.}
\label{fig:GaussVSGamma}
\end{figure}

\begin{remark}
In this work, we only consider parameters having a PDF with bounded support. Indeed, in most applications these parameters have a physical meaning, for instance they could represent a quantity which must always be nonnegative, and having PDF with unbounded supports like Gaussian would mean that said parameters would be negative with a positive probability, which is not realistic. However, many sources of uncertainties are still expected to have a Gaussian-like behavior, in the sense that their PDF should be concentrated around a point. In that case, a good compromise is to use Beta distributions, which are the weights associated to the Gegenbauer polynomials and provide good bounded approximations of Gaussian distributions, at least for small variances (see e.g.~\cite{Xiu04}, or the comparison on Figure~\ref{fig:GaussVSGamma}). Another widely used option is to use truncated Gaussian distributions.
\end{remark}

We recall that if $\rho$ has compact support, or decays at least exponentially fast at infinity (see e.g.~\cite{ErnMugStaUll12}), then $\left(\phi_n\right)_{n\in\N}$ is a Hilbert basis of $L^2(\rho\txtd s)$,
 i.e. any mesurable function $g$ such that
\begin{equation*}
\int_\R g^2(s) \rho(s)~\txtd s <\infty,
\end{equation*}
admits a unique series expansion of the form
\begin{equation*}
g=\sum_{n\in\N}g_n\phi_n,
\end{equation*}
where the series converges in $L^2(\rho\txtd s)$.

Now, assume that the noisy parameter $p=p(\omega)$ (still assumed to be one dimensional for the moment) has a PDF given by $\rho$. The PC paradigm then tells us that we should use the orthogonal polynomials associated to $\rho$ as a basis for the expansion~\eqref{eq:expansion1} of $X(p)$. Notice that with such a choice, the coefficients $X_n$ of the expansion directly provide us with the mean and variance of $X(p)$. Indeed, by orthogonality (and assuming $h_n=1$ for simplicity) we have
\begin{equation*}
\E \left(X(p)\right) = \sum_{n=0}^\infty X_n \E \left(\phi_n(p)\right) = \sum_{n=0}^\infty X_n \langle \phi_n,1\rangle_\rho = X_0,
\end{equation*}
and similarly
\begin{equation*}
\V \left(X(p)\right) = \E \left(X(p)^2\right) - \E \left(X(p)\right)^2 =  \sum_{m=0}^\infty\sum_{n=0}^\infty X_nX_m \langle \phi_m,\phi_n\rangle_\rho - X_0^2 = \sum_{n=1}^\infty X_n^2.
\end{equation*}
If $p=(p^{(1)},\ldots,p^{(d)})$ consists of several independent random variables, each with respective PDF $\rho^{(j)}$, $j=1,\ldots,d$, one can consider the PC basis $\left(\phi_n\right)_{n\in\N^d}$ constructed as a tensor product of the univariates bases. That is, for all $n=(n_1,\ldots,n_d)\in\N^d$ and all $s=(s_1,\ldots,s_d)\in\R^d$,
\begin{equation*}
 \phi_n(s):=\prod_{j=1}^d \phi_{n_j}^{(j)}(s_j),
\end{equation*}
where $(\phi_{n_j}^{(j)})_{n_j\in\N}$ is a basis of orthogonal polynomials associated to the weight $\rho^{(j)}$. In practice, there are several ways to compute the coefficients $X_n$ of a PC expansion~\eqref{eq:expansion1}. Notice that directly using the orthogonality relations to get
\begin{equation*}
X_n = \frac{1}{h_n}\langle X , \phi_n \rangle_\rho,
\end{equation*}
is usually not one of them, as $X$ is not known a-priori but is actually what we want to compute via the PC expansion. To formalize the discussion, assume that the quantity of interest $X(p)\in\R^m$ for which we want to find a PC expansion, solves a problem depending on a parameter $p\in\R^d$, of the form
\begin{equation*}
F(X(p),p)=0,
\end{equation*}
where $F:\R^m\times\R^d\to\R^m$. One common way to find the coefficients $X_n$ is to solve the system obtained by Galerkin projection:
\begin{equation}
\label{eq:Galerkin}
\left\langle F\left(\sum_{n_1=0}^\infty X_{n_1} \phi_{n_1}(p),p\right),\phi_{n_2}(p) \right\rangle_\rho = 0 \qquad \forall~{n_2}\in\N.
\end{equation}
Of course, in practice one only solves for a truncated expansion
\begin{equation*}
X^N(p)=\sum_{{n_1}=0}^N X_{n_1} \phi_{n_1}(p),
\end{equation*}
by considering a finite-dimensional projection of associated dimension, i.e. $\forall~n_2\leq N$ in~\eqref{eq:Galerkin}. This is the technique we make use of in this work. Another possible option is to use a collocation/interpolation approach~\cite{Xiu07}. This technique is based on first solving the deterministic problem several times for a well chosen sample of parameter values $\left(p_i\right)_{1\leq i \leq i}$, i.e. computing $X(p_i)$ that solves
\begin{equation*}
F(X(p_i),p_i)=0,\qquad \forall~1\leq i\leq I.
\end{equation*}
The polynomials chaos coefficients $X_n$ are then constructed by interpolation:
\begin{equation*}
X_n = \sum_{i=1}^I X(p_i)\phi_n(p_i)\alpha_i,
\end{equation*}
where $\alpha_i$ are weights associated to the interpolation points $p_i$ in such a way that, for any smooth function $g$
\begin{equation*}
\sum_{i=1}^I g(p_i)\alpha_i \approx \int g(s)\rho(s) ~\txtd s,
\end{equation*}
where $\rho$ is the PDF of $p$. For a detailed discussion about these two approaches and their respective merits and limitations, we refer to the survey~\cite{Xiu09}, the book~\cite{LeMKni10} and the references therein. In practice, after an accurate PC approximation
\begin{equation*}
X^N(p)=\sum_{{n}=0}^N X_{n} \phi_{n}(p)
\end{equation*}
of $X$ has been computed, we can then do Monte Carlo simulations for a very large sample of values of $p$. Indeed, instead of having to solve the deterministic problem $F(X(p_i),p_i)=0$ for each value $p_i$ of the sample, we only have to evaluate $\phi_n(p_i)$ to obtain $X^N(p_i)\approx X(p_i)$.

\section{Computation of periodic orbits, invariant manifolds and connecting orbits}
\label{sec:inv_sets}

In this section, we review some classical techniques to numerically study periodic orbits, invariant manifolds and connecting orbits of deterministic ODEs. An exhaustive review of these techniques is far beyond the scope of this work, and we only focus on one technique for each case, although alternative methods could certainly also be used; see e.g.~\cite{KrauskopfSurvey} for a comparison of methods for computing stable/unstable manifolds. All the techniques presented here are based on series expansions, and this choice is motivated by two main reasons. The first and most important one is that series expansions can easily and efficiently be combined with PC expansions, once we go to the stochastic setting. This is particularly relevant for limit cycles, and is to be contrasted with more classical long-term integration approaches, for which PC expansions are known to be ill-behaved. The second one is that series expansions are particularly well suited for a-posteriori validation techniques~\cite{BerLes15}, which we plan on developing in this context in a future work. In each case, we explain how an extra layer of PC expansion can be added, to keep track of the stochastic nature of the parameters. Illustrations for all these methods are presented in Section~\ref{sec:LV} and Section~\ref{sec:Lorenz}.

\subsection{Periodic orbits via Fourier series}
\label{sec:inv_periodic}

We start by considering the parameter dependent problem~\eqref{eq:ODE_rand_coef} from a deterministic point of view. Limit cycles of~\eqref{eq:ODE_rand_coef} can be efficiently studied and computed using Fourier series. That is, for a fixed $p$, we write a $T(p)$-periodic solution $t\mapsto X(t,p)$ of~\eqref{eq:ODE_rand_coef} as
\begin{equation}
\label{eq:Fourier_exp}
X(t,p)=\sum_{k\in\Z}X_k(p) \txte^{\txti k\Omega(p) t}, \quad X_k(p)\in\R^m,\ \Omega(p)\in\R,
\end{equation}
where $\Omega(p)=2\pi/T(p)$. To numerically find a periodic orbit, we thus solve for the coefficients $X_k(p)$ and $\Omega(p)$, which satisfy the system obtained by plugging~\eqref{eq:Fourier_exp} into~\eqref{eq:ODE_rand_coef}, namely
\begin{equation}
\label{eq:ODE_Fourier}
\txti k\Omega(p) X_k(p) = f_k(X(t,p),p) \qquad \forall k\in\Z,
\end{equation}
where $f_k(X(t,p),p)$ are the Fourier coefficients of $t\mapsto f(X(t,p),p)$.

\begin{remark}
In practice, it is helpful to complement the above system, with a \emph{phase condition}. For simplicity we choose to fix a section through which the solution has to pass at time $0$, but other options are available~\cite{Kuz13}. More precisely, we add a scalar equation of the form
\begin{equation*}
\left(X(0,p)-u(p)\right)\cdot v(p) =0,
\end{equation*}
where $u(p)$ is some (approximate) point on the orbit, $v(p)\approx f(u(p))$, and we use the dot product to denote the scalar product on $\R^n$ (and avoid potential confusion with the $L^2(\rho)$ scalar product $\langle,\rangle_\rho$). The phase condition allows to isolate the solution by eliminating time-shifts, which makes the system easier to solve in practice, especially using iterative methods. Indeed, if $t\mapsto X(t,p)$ is a periodic orbit of~\eqref{eq:ODE_rand_coef}, then so is any function of the form $t\mapsto X(t+\tau,p)$, but there is (locally) only one $\tau$ that also satisfies the phase condition.
\end{remark}

This approach of solving for the Fourier coefficients and the frequency has several advantages, compared to numerically integrating~\eqref{eq:ODE_rand_coef} and trying to find an (approximately) closed orbit. Indeed, the approach is not sensitive to the linear stability or instability of the limit cycle, and is therefore not susceptible of diverging if one tries to approximate an unstable periodic orbit, which is a definitive concern for time-integration based techniques. We illustrate this point in Section~\ref{sec:Lorenz} by computing unstable limit cycles belonging to the chaotic attractor of the Lorenz system. The Fourier series approach is also particularly well adapted to continuation algorithms, especially to compute limit cycles originating from a Hopf bifurcation, where the linearized analysis can predict the first Fourier coefficient and the frequency.

If we now consider that $p$ is random and has a given PDF $\rho_p$, it is natural to still consider a Fourier ansatz~\cite{MillmanKingBeran}. It is straightforward to extend the Fourier series approach by expanding each Fourier coefficient, together with the frequency, with PC. Namely, write
\begin{equation*}
X_k(p)=\sum_{n\in\Z} X_{k,n} \phi_n(p) \quad \text{and}\quad \Omega(p)=\sum_{n\in\N} \Omega_n \phi_n(p),
\end{equation*}
or equivalently
\begin{equation}
\label{eq:Fourier_PC_exp}
X(t,p)=\sum_{k\in\Z}\sum_{n\in\N} X_{k,n} \phi_n(p) \txte^{\txti kt \Omega (p)}.
\end{equation}
In practice, we of course consider truncated expansions
\begin{equation*}
X_k^N(p)=\sum_{n=0}^N X_{k,n} \phi_n(p) \quad \text{and}\quad \Omega(p)=\sum_{n=0}^N \Omega_n \phi_n(p).
\end{equation*}
for which we solve using~\eqref{eq:ODE_Fourier}. Each $X_k$ now belongs to $\R^{mN}$ instead of $\R^m$, and $\Omega$ belongs to $\R^N$ instead of $\R$. An explicit example of such system together with numerical solutions is given in Section~\ref{sec:Lorenz}. In this setting, the Fourier series approach also has the notable advantage of separating the random period (or equivalently the random frequency $\Omega(p)$), from the description of the random cycle given by the coefficients $X_k(p)$. In particular, we avoid the usual pitfalls related to \emph{phase-drift} and \emph{broadening of the spectrum}, which are the main reasons why limit cycles are hard to compute using PC combined with time integration~\cite{DesaiWitteveenSarkar}. A similar idea was introduced in~\cite{LeMMatKniHus10,SchHeuLeM14}, where a random time rescaling is used to compensate for the random period, which significantly improves the long time behavior of the PC expansions and allows to better capture stable limit cycles. Yet a Fourier series approach accomplishes that \emph{naturally}, and can also be used to study unstable limit cycles.

\subsection{Local invariant manifolds via Taylor series and the parameterization method}
\label{sec:inv_manifold}

\begin{figure}[h]
\centering
\begin{overpic}
[width=0.65\linewidth]{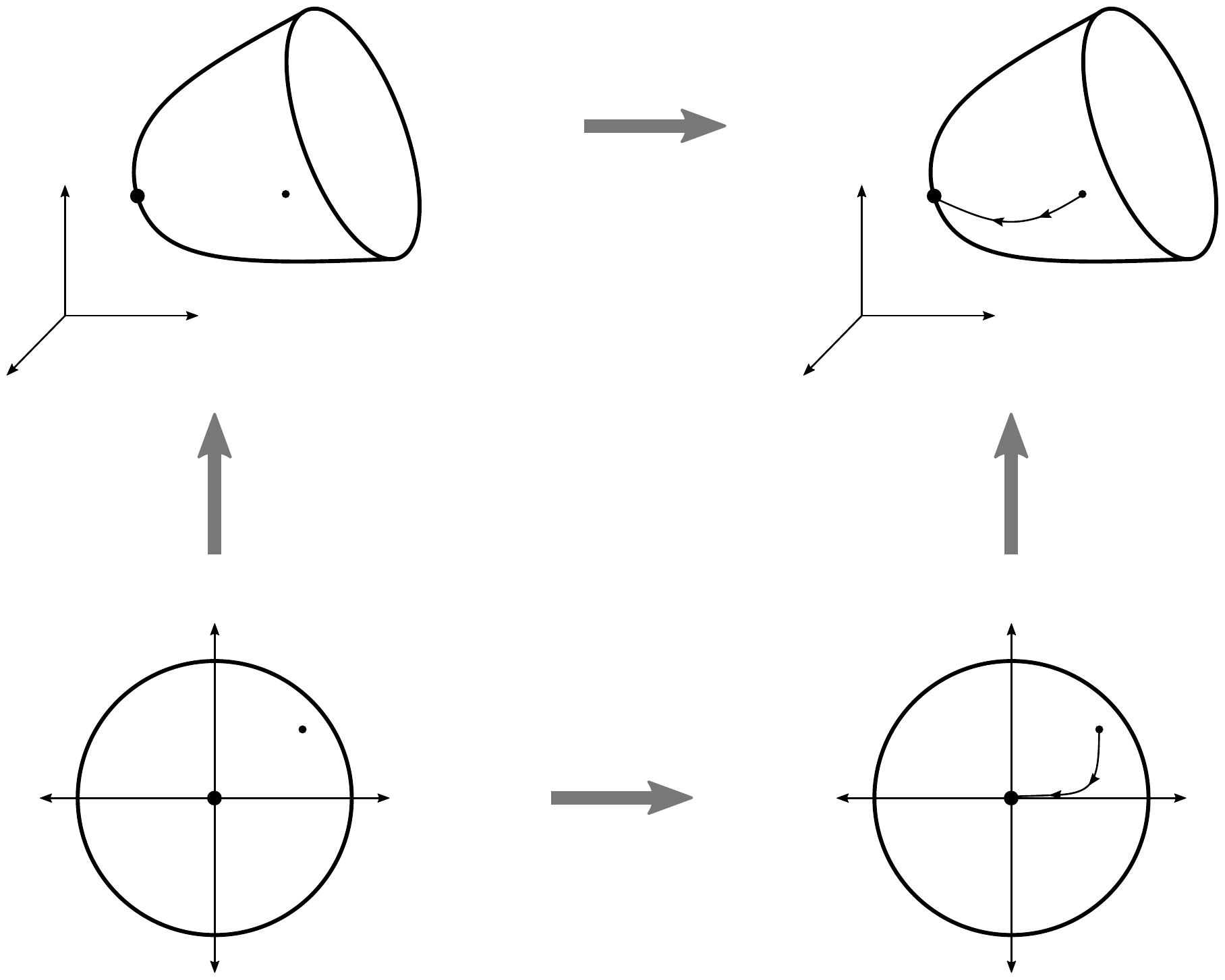}
\put (5,2) {\large$\bm{\R^{m_s}}$}
\put (70,2) {\large$\bm{\R^{m_s}}$}
\put (6,50) {\large$\bm{\R^m}$}
\put (71,50) {\large$\bm{\R^m}$}
\put (22,20) {\large$\bm{\theta}$}
\put (87,20) {\large$\bm{\theta}$}
\put (19,66) {\large$\bm{Q(\theta,p)}$}
\put (84,66) {\large$\bm{Q(\theta,p)}$}
\put (3,39) {\large$\bm{Q(\cdot,p)}$}
\put (84.5,39) {\large$\bm{Q(\cdot,p)}$}
\put (47,73) {\large$\bm{\varphi_p(t,\cdot)}$}
\put (46,9) {\large$\bm{\txte^{t\Lambda(p)}}$}
\end{overpic}
\caption{Schematic illustration of the parameterization method. We want the parameterization $Q$ to conjugate the nonlinear flow $\varphi_p$ to the linearized flow generated by $\Lambda(p)$ on the stable subspace.}
\label{fig:PM}
\end{figure}
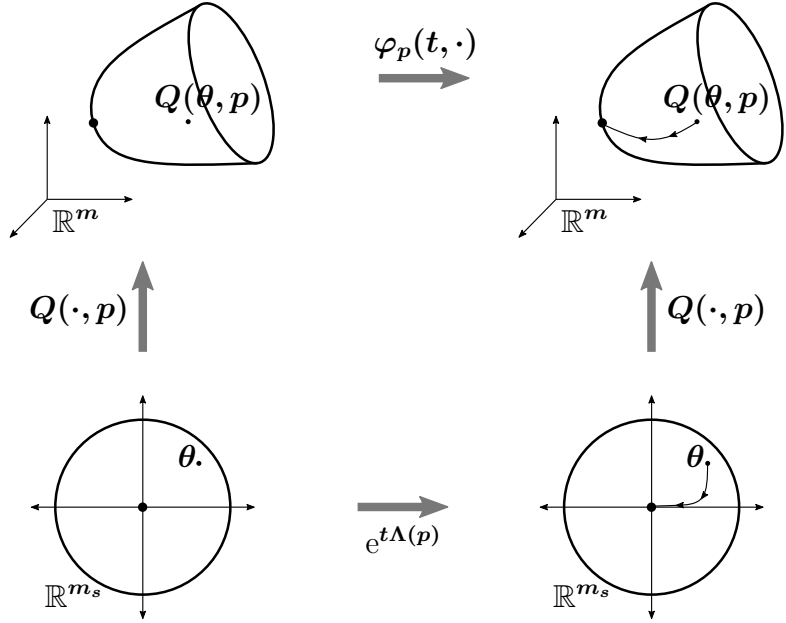

We again start by considering the parameter dependent problem~\eqref{eq:ODE_rand_coef} from a deterministic point of view, but now focus on studying local stable and unstable manifold attached to equilibrium points. In this section we only discuss the case of stable manifolds, but unstable manifolds can of course be studied with the same techniques. Let $\hat X(p)$ be an equilibrium point of~\eqref{eq:ODE_rand_coef}, i.e. $f(\hat X(p),p)=0$, and assume that the derivative $\txtD_xf(\hat X(p),p)$ has exactly $m_s\leq m$ eigenvalues $\lambda^{(1)}(p),\ldots,\lambda^{(m_s)}(p)$ with negative real part. For simplicity, we also assume that each of these eigenvalues is simple, and denote by $v^{(1)}(p),\ldots,v^{(m_s)}(p)$ associated eigenvectors. Our goal is to find a parameterization $Q(p):\R^{m_s}\to\R^m$ of the local stable manifold of $\hat X(p)$. We look for a power series representation of $Q$:
\begin{equation}
\label{eq:Q}
Q(\theta,p) = \sum_{\vert k\vert \geq 0} Q_k(p) \theta^k, \quad
\theta = \begin{pmatrix} \theta_1 \\ \vdots \\ \theta_{m_s}\end{pmatrix} \in \R^{m_s},\quad 
Q_k(p) = \begin{pmatrix} Q_k^{(1)}(p) \\ \vdots \\ Q_k^{(m)}(p) \end{pmatrix} \in \R^m, 
\end{equation}
with the classical multi-indexes notations $\vert k\vert = k_1+\cdots+k_{m_s}$ and $\theta^k=\theta_1^{k_1}\cdots\theta_{m_s}^{k_{m_s}}$. Since we want $\theta\mapsto Q(\theta,p)$ to be a parameterization of the local stable manifold of $\hat X(p)$, we must have
\begin{equation}
\label{eq:CI_para}
Q_0(p)=\hat X(p) \qquad \text{and} \qquad Q_{e_i}(p) = \gamma_i v^{(i)}(p),\ \forall~1\leq i\leq {m_s},
\end{equation}
where $\gamma_i$ are scaling that can be adjusted. To obtain the higher order terms, we follow the idea of the \emph{parameterization method}, introduced in~\cite{CabFonLla03,CabFonLla03bis,CabFonLla05} (see also the recent book~\cite{HarCanFigLuqMon16}). We want to obtain a parameterization that conjugates the dynamics on the stable manifold with the stable dynamics of the linearized system. More precisely, introducing the diagonal matrix $\Lambda(p)$ with diagonal entries $\lambda^{(1)}(p),\ldots,\lambda^{({m_s})}(p)$, we want $p$ to satisfy (see Figure~\ref{fig:PM})
\begin{equation}
\label{eq:conjugacy}
\varphi_p(t,Q(\theta,p))=Q(\txte^{\Lambda(p)t}\theta,p), \qquad \forall~ \left\Vert \theta\right\Vert_\infty\leq 1,\ \forall~t\geq 0,
\end{equation}
where $\varphi_p$ is the flow generated by the vector field $f(\cdot,p)$ and $\left\Vert \theta\right\Vert_\infty\leq 1 = \max_{1\leq i \leq {m_s}} \vert \theta_i\vert$.

\begin{remark}
If some of the stable eigenvalues are complex conjugate, say $\lambda^{(1)}=\bar\lambda^{(2)},\ldots,\lambda^{(2r-1)}=\bar\lambda^{(2r)}$, it is more convenient to first look for a complex valued parameterization $Q$ with $Q_k\in\C^m$, and then recover a real valued parameterization via
\begin{equation*}
Q_{\textnormal{real}}(\theta,p)=Q(\theta_1+\txti\theta_2,\theta_1-\txti\theta_2,\ldots,\theta_{2r-1}+\txti\theta_{2r},\theta_{2r-1}-\txti\theta_{2r},\theta_{2r+1},\ldots,\theta_{m_s},p),
\end{equation*}
see e.g.~\cite{BerMirRei16}.
\end{remark}

Finding a parameterization $Q$ satisfying~\eqref{eq:conjugacy} is interesting because it provides us with not only a local stable manifold but also an explicit description of the dynamics on this manifold. However, the formulation~\eqref{eq:conjugacy} is not the most convenient one to work with in order to determine the higher order coefficients of $Q$, because it involves the flow. To get rid of it, one can take a time derivative of~\eqref{eq:conjugacy} and evaluate at $t=0$, to obtain the following \emph{invariance equation}
\begin{equation}
\label{eq:invariance}
f(Q(\theta,p),p)=\txtD_\theta Q(\theta,p)\Lambda(p)\theta ,\quad \forall~ \left\Vert \theta\right\Vert_\infty\leq 1.
\end{equation}
One can check that, if $Q$ solves~\eqref{eq:invariance} and is such that $Q(0,p)=\hat X(p)$, then $Q$ satisfies~\eqref{eq:conjugacy}, therefore $Q(\cdot,p)$ is indeed a parameterization of the local unstable manifold of $\hat X(p)$. The invariance equation~\eqref{eq:invariance} is the one we are going to use to numerically find the coefficients $Q_k$. Plugging the expansion~\eqref{eq:Q} in the invariance equation~\eqref{eq:invariance}, we obtain
\begin{equation}
\sum_{\vert k\vert \geq 0} f_k(Q(\theta,p),p)\theta^k = \sum_{\vert k\vert \geq 0} (\lambda(p)\cdot k) Q_k(p)\theta^k,
\end{equation}
where $f_k(Q(\theta,p),p)$ are the Taylor coefficients of $\theta\mapsto f(Q(\theta,p),p)$ and $\lambda(p)\cdot k = \lambda^{(1)}(p)k_1+\ldots + \lambda^{({m_s})}(p) k_{m_s}$, therefore the coefficients $Q_k$ must satisfy
\begin{equation}
\label{eq:invariance_coeffs}
(\lambda(p)\cdot k) Q_k(p) = f_k(Q(\theta,p),p),\qquad \forall k\in\N^{m_s}.
\end{equation}
Notice that~\eqref{eq:CI_para} already ensures that~\eqref{eq:invariance_coeffs} is satisfied for $\vert k\vert \leq 1$. In order to solve~\eqref{eq:invariance_coeffs} for $\vert k\vert \geq 2$, let us first describe, how $f_k(Q(\theta,p),p)$ depends on $Q_k(p)$. Given a Taylor series of the form~\eqref{eq:Q}, we define for all $K\in\N$ the truncated series
\begin{equation*}
\pi_K Q (\theta,p) = \sum_{0 \leq \vert k\vert \leq K} Q_k(p) \theta^k.
\end{equation*}
Notice also that, for any $k\in\N^{m_s}$, $f_k(Q(\theta,p),p) = f_k(\pi_{\vert k\vert} Q(\theta,p),p)$. Besides, using a Taylor expansion of $f$ in the $x$ variable, we have for all $k\in\N^{m_s}\setminus\{0\}$
\begin{align*}
f\left(\sum_{0 \leq \vert l\vert \leq \vert k\vert } Q_l(p)\theta^l,p\right) &= f\left(\sum_{0 \leq \vert l\vert \leq \vert k\vert -1} Q_l(p)\theta^l + \sum_{\vert l\vert =\vert k\vert}Q_l(p)\theta^l,p\right)\\
&= f\left(\sum_{0 \leq \vert l\vert \leq \vert k\vert -1} Q_l(p)\theta^l,p\right) + \sum_{\vert l\vert =\vert k\vert} \txtD_xf\left(\sum_{0 \leq \vert j\vert \leq \vert k\vert-1} Q_j(p)\theta^j,p\right)Q_l(p)\theta^l \\
&\quad + \text{higher order terms},
\end{align*}
and looking at the coefficient of degree $k$ on each side we get
\begin{equation*}
f_k(Q(\theta,p),p) = f_k(\pi_{\vert k\vert} Q(\theta,p),p) = f_k(\pi_{\vert k\vert-1} Q(\theta,p),p) + \txtD_xf(Q_0(p),p)Q_k(p).
\end{equation*}
Therefore, assuming~\eqref{eq:CI_para}, having~\eqref{eq:invariance_coeffs} for all $\vert k\vert \geq 2$ is equivalent to having
\begin{equation}
\label{eq:homological}
\left((\lambda(p)\cdot k)-\txtD_xf(\hat X(p),p)\right) Q_k(p) = f_k(\pi_{\vert k\vert-1} Q(\theta,p),p), \quad \forall~\vert k\vert \geq 2.
\end{equation}
Assuming the following \emph{non-resonance} condition is satisfied:
\begin{equation*}
\lambda(p)\cdot k \neq \lambda^{(i)}(p) \qquad \forall~\vert k\vert \geq 2,\ \forall~1\leq i\leq r,
\end{equation*}
we see that~\eqref{eq:homological} has a unique solution that can be computed recursively via
\begin{equation}
\label{eq:Q_k_rec}
Q_k(p) = \left((\lambda(p)\cdot k)-\txtD_xf(\hat X(p),p)\right)^{-1}f_k(\pi_{\vert k\vert-1} Q(\theta,p),p) \qquad \forall~\vert k\vert \geq 2,
\end{equation}
since the right-hand side only depends on $Q_l(p)$ for $\vert l\vert <\vert k\vert$. We can therefore compute a truncated parameterization $\pi_KQ$ of arbitrary order, starting from~\eqref{eq:CI_para} and then computing~\eqref{eq:Q_k_rec} recursively for $k$ as large as desired. In practice, the weights $\gamma_i$ in~\eqref{eq:CI_para} are chosen in order to obtain a reasonnable decay of the coefficients $Q_k(p)$ (we refer to~\cite{BreLesMir16} for a detailed explanation of how this choice can be optimized). Explicit examples are presented in Sections~\ref{sec:LV} and~\ref{sec:Lorenz}.

\begin{remark}
In cases where resonant eigenvalues are present, a similar approach can still be used, but the conjugacy condition~\eqref{eq:conjugacy} defining $Q$ has to be adapted~\cite{BerMirRei16}.
\end{remark}

If we now consider that $p$ is random and has a given PDF $\rho_p$, this approach based on the parameterization method can also be easily generalized by adding a layer of PC expansion. Namely, we write
\begin{equation*}
Q(\theta,p) = \sum_{\vert k\vert \geq 0}\sum_{n\in\N} Q_{k,n} \phi_n(p) \theta^k,
\end{equation*}
or equivalently 
\begin{equation*}
Q_k(p)=\sum_{n\in\N} Q_{k,n} \phi_n(p).
\end{equation*}
In practice we consider a truncation
\begin{equation*}
Q_k^N(p)=\sum_{n=0}^N Q_{k,n} \phi_n(p), 
\end{equation*}
which we compute via
\begin{equation}
\label{eq:Q_k_rec_PC}
Q_k^N(p) = \left((\lambda^N(p)\cdot k)-\txtD_xf(\hat X^N(p),p)\right)^{-1}f_k(\pi_{\vert k\vert-1} Q^N(\theta,p),p) \qquad \forall~\vert k\vert \geq 2,
\end{equation}
where, compared to~\eqref{eq:Q_k_rec}, $Q_k^N(p)$ is now a vector of size $mN$ rather than $m$, and 
\begin{equation*}
\left((\lambda^N(p)\cdot k)-D_xf(\hat X^N(p),p)\right) 
\end{equation*}
can be interpreted as a block matrix, with $m\times m$ blocks having each size $N\times N$. For explicit computations we again refer to Sections~\ref{sec:LV} and~\ref{sec:Lorenz}.

\subsection{Heteroclinic orbits via Chebyshev series and projected boundaries}
\label{sec:inv_heteroclinic}

We go back to considering the parameter dependent problem~\eqref{eq:ODE_rand_coef} from a deterministic point of view, and extend the discussion of the previous subsection by looking at more global solutions, namely connecting orbits between equilibrium points.

Let $\hat X(p)$ and $\check X(p)$ be two equilibrium points of~\eqref{eq:ODE_rand_coef}. Assume that $\hat X(p)$ has an unstable manifold of dimension $\hat r$ and that $\check X(p)$ has a stable manifold of dimension $\check r$, such that $\hat r+\check r = m+1$. Then, if the two manifolds intersect, we can generically expect this intersection to be transverse in the phase space $\R^m$, in which case there exists a transverse heteroclinic orbit between $\hat X(p)$ and $\check X(p)$. We recall that a heteroclinic orbit between $\hat X(p)$ and $\check X(p)$ (or homoclinc orbit if $\hat X(p)=\check X(p)$) is a solution $t\mapsto x(t,p)$ of~\eqref{eq:ODE_rand_coef} such that
\begin{equation*}
\lim_{t\to -\infty} x(t,p) = \hat X(p) \quad \text{and}\quad \lim_{t\to +\infty} x(t,p) = \check X(p).
\end{equation*}
In this work, we compute such solution by solving a boundary value problem~\cite{AscherMattheijRussell} between the unstable manifold of $\hat X(p)$ and the stable manifold of $\check X(p)$, for which we first compute local parameterization as in Section~\ref{sec:inv_manifold}. This allows us to only solve~\eqref{eq:ODE_rand_coef} on a finite time interval, and recover the remaining parts of the orbit via the conjugacy satisfied by the parameterizations. More precisely, we want to find an orbit $t\mapsto X(t,p)$ such that
\begin{equation}
\left\{\begin{aligned}
\label{eq:BVP}
&\dot{X}(t,p) = f(X(t,p),p) \quad \forall~t\in[0,\tau(p)] \\
&X(0,p) \in W^\txtu(\hat X(p)) \\
&X(\tau(p),p) \in W^\txts(\check X(p)),
\end{aligned}\right.
\end{equation}
where $W^\txtu(\hat X(p))$ and $W^\txts(\check X(p))$ respectively denote the unstable manifold of $\hat X(p)$ and the stable manifold of $\check X(p)$. Note that, since we only consider autonomous vector fields in this work, only the length $\tau(p)$ of the time interval $[0,\tau(p)]$ is relevant, and the whole interval itself could of course be shifted. To numerically compute such an orbit, we use piece-wise Chebyshev series. In order to do so, we first introduce some notations.

For $k\in\N$, we denote by $T_k$ the Chebyshev polynomial of order $k$, defined for instance by $T_k(\cos(\theta))=\cos(k\theta)$. We introduce a partition $0=\tilde t^{(0)}<\tilde t^{(1)}<\ldots<\tilde t^{(J)}=1$ of $[0,1]$. We denote by $\tau(p)$ the (unknown) time the orbit spends between the two local manifolds and consider the partition of $[0,\tau(p)]$ given by $0=t^{(0)}(p)<t^{(1)}(p)<\ldots<t^{(J)}(p)=\tau(p)$, where $t^{(j)}(p)=\tau(p)\tilde t^{(j)}$ for all $1\leq j\leq J$. For all $k\in\N$ and $j=1,\ldots,J$, we also introduce the rescaled Chebyshev polynomial $T^{(j)}_k(p)$, defined as
\begin{equation*}
T^{(j)}_k(t,p)=T_k\left(\frac{2t-t^{(j)}(p)-t^{(j-1)}(p)}{t^{(j)}(p)-t^{(j-1)}(p)}\right).
\end{equation*}
We can then write the orbit between the two manifolds as
\begin{equation}
\label{eq:piecewise_cheb}
X(t,p) = X^{(j)}_0(p)+2\sum_{k=1}^\infty X^{(j)}_k(p) T^{(j)}_k(t,p), \quad \forall~t\in(t^{(j-1)}(p),t^{(j)}(p)),\ \forall~j\in\{1,\ldots,J\}.
\end{equation}
Finally, we assume that, following the methodology presented in Section~\ref{sec:inv_manifold}, a truncated parameterization $\hat Q(p)$ of the local unstable manifold of $\hat X(p)$ as well as a truncated parameterization $\check Q(p)$ of the local stable manifold of $\check X(p)$ have been computed. Rewriting the differential equation in~\eqref{eq:BVP} as an integral one, plugging in the expansion~\eqref{eq:piecewise_cheb} and using well known properties of the Chebyshev polynomials, namely
\begin{equation*}
T_k(1)=1,\quad T_k(-1)=(-1)^k\quad \text{and} \quad \int T_k =\frac{1}{2}\left(\frac{T_{k+1}}{k+1}-\frac{T_{k-1}}{k-1}\right),
\end{equation*}
we obtain
\begin{equation}
\label{eq:BVP_coeffs}
\left\{\begin{aligned}
&kX^{(j)}_k(p)=\tau(p)\frac{\tilde t^{(j)}-\tilde t^{(j-1)}}{4}\left(f^{(j)}_{k-1}(X(t,p),p)-f^{(j)}_{k+1}(X(t,p),p)\right),\quad \forall~k\geq 1,\ \forall~1\leq j\leq J \\
&X^{(j)}_0(p)+2\sum_{k=1}^\infty X^{(j)}_k(p) = X^{(j+1)}_{0}(p)+2\sum_{k=1}^\infty (-1)^kX^{(j+1)}_{k}(p),\quad \forall~1\leq j\leq M-1 \\
&X^{(1)}_0(p)+2\sum_{k=1}^\infty (-1)^kX^{(1)}_k(p) = \hat Q(\hat \theta(p),p) \\
&X^{(M)}_0(p)+2\sum_{k=1}^\infty X^{(J)}_k(p) = \check Q(\check \theta(p),p),
\end{aligned}\right.
\end{equation}
where $f^{(j)}_{k}(X(t,p),p)$ are the Chebyshev coefficients of $t\mapsto f(X(t,p),p)$ on $(t^{(j-1)}(p),t^{(j)}(p))$. The first line in~\eqref{eq:BVP_coeffs} corresponds to the differential equation $\dot{X}(t,p) = f(X(t,p),p)$ on each subinterval $(t^{(j-1)}(p),t^{(j)}(p))$, the second line insures that the solution connects continuously between two consecutive subintervals, and the last two lines correspond to the boundary conditions on each manifold.

\begin{remark}
\label{rem:para_BVP}
In~\eqref{eq:BVP_coeffs}, besides the unknown Chebyshev coefficients $X_k^{(j)}(p)$, we have $\hat r + \check r +1$ additional scalar unknowns: $\tau(p)\in\R$, $\hat\theta(p)\in\R^{\hat r}$ and $\check\theta(p)\in\R^{\check r}$. Still assuming $\hat r + \check r=m+1$, the system is then underdetermined and we can fix two of these $m+2$ parameters to recover a unique solution.
\end{remark}

Similarly to the two previous cases of Fourier and Taylor series, when we consider $p$ as a random parameter having a given distribution we only need to expand every unknown in~\eqref{eq:BVP_coeffs} using PC. We again refer to Section~\ref{sec:Lorenz} for an explicit example.

\section{First example: a Lotka-Volterra system}
\label{sec:LV}

In this section, we compute steady states and invariant manifolds of a Lotka-Volterra system of competition type
\begin{equation}
\label{eq:LV}
\left\{\begin{aligned}
\dot{x}&=(1-x-ay)x,\\
\dot{y}&=(1-y-bx)y,
\end{aligned}\right.
\end{equation}
where $b=b(\omega)$ is a bounded random variable having a given distribution. This basic example allows us to easily study the quality of our numerical computations by comparing them to analytic results.

\subsection{Analytical results (for benchmarking)}

\subsubsection{Deterministic framework}

Assuming $ab\neq 1$, system~\eqref{eq:LV} has a non trivial equilibrium given by
\begin{equation*}
(x_{\textnormal{eq}},y_{\textnormal{eq}})=\left(\frac{a-1}{ab-1},\frac{b-1}{ab-1}\right).
\end{equation*}
The eigenvalues of the Jacobian $\txtD f(x_{\textnormal{eq}},y_{\textnormal{eq}})$ are 
\begin{equation}
\label{eq:eigenval_LV}
\lambda^{(1)} = -1 \qquad \text{and} \qquad \lambda^{(2)} = \frac{(a-1)(b-1)}{ab-1}.
\end{equation}
and associated eigenvectors are given by
\begin{equation}
\label{eq:eigenvect_LV}
V_{\lambda^{(1)}}=\begin{pmatrix} a-1 \\ b-1 \end{pmatrix} \qquad \text{and} \qquad V_{\lambda^{(2)}}=\begin{pmatrix} a \\ -b \end{pmatrix}.
\end{equation}
In the case $a,b>1$, $\lambda^{(2)}$ is positive and the equilibrium is of saddle type. Its stable manifold is the line $y=\frac{b-1}{a-1} x$, $x>0$ and its unstable manifold is a (nonlinear) curve connecting to $(0,1)$ and $(1,0)$.

\subsubsection{Stochastic framework}
\label{sec:LV_sto}

We now assume that the parameter $b=b(\omega)$ is a bounded random variable. We write
\begin{equation*}
b(\omega) = \bb + \sigma p(\omega),
\end{equation*}
where $\sigma\geq 0$ and $p$ is a random variable taking values in $[-1,1]$. We assume $a\bb\neq 1$ and $\sigma \leq \left\vert\frac{a\bb-1}{a}\right\vert$, so that $ab(\omega)-1\neq 0$ a.e.. For specific distributions of $p$, one can compute analytically the first two moments of $x_{\textnormal{eq}}=\frac{a-1}{ab-1}$:
\begin{itemize}
\item If $p$ is uniformly distributed on $[-1,1]$, that is if its PDF is given by $\rho_p(s)=\frac{1}{2}\1_{s\in(-1,1)}$, one has
\begin{equation*}
\E(x_{\textnormal{eq}})=\frac{a-1}{a\sigma}\atanh\left(\frac{a\sigma}{a\bb-1}\right) \qquad \text{and}\qquad \E(x_{\textnormal{eq}}^2)=\frac{(a-1)^2}{(a\bb-1)^2-(a\sigma)^2}.
\end{equation*}
\item If $p$ is beta-distributed with parameters $(-\frac{1}{2},-\frac{1}{2})$, that is if its PDF is given by $\rho_p(s)=\1_{s\in(-1,1)}\frac{1}{\pi\sqrt{1-s^2}}$, one has
\begin{equation*}
\E(x_{\textnormal{eq}})=\frac{a-1}{a\bb-1}\frac{1}{\sqrt{1-\left(\frac{a\sigma}{a\bb-1}\right)^2}} \qquad \text{and}\qquad \E(x_{\textnormal{eq}}^2)=\left(\frac{a-1}{a\bb-1}\right)^2\frac{1}{\left(1-\left(\frac{a\sigma}{a\bb-1}\right)^2\right)^{\frac{3}{2}}}.
\end{equation*}
\item If $p$ is beta-distributed with parameters $(\frac{1}{2},\frac{1}{2})$, that is if its PDF is given by $\rho_p(s)=\1_{s\in(-1,1)}\frac{2}{\pi}\sqrt{1-s^2}$, one has
\begin{equation*}
\E(x_{\textnormal{eq}})=\frac{a-1}{a\bb-1}\frac{2}{1+\sqrt{1-\left(\frac{a\sigma}{a\bb-1}\right)^2}} \qquad \text{and}\qquad \E(x_{\textnormal{eq}}^2)=\frac{2\left(1-\sqrt{1-\left(\frac{a\sigma}{a\bb-1}\right)^2}\right)}{\left(\frac{a\sigma}{a\bb-1}\right)^2\sqrt{1-\left(\frac{a\sigma}{a\bb-1}\right)^2}}\left(\frac{a-1}{a\bb-1}\right)^2.
\end{equation*}
\end{itemize}

In order to focus the amount of comparisons and illustrations in the next section, we only consider the first component $x_{\textnormal{eq}}$ of the equilibrium, but similar analytical and numerical computations could of course also be carried out for $y_{\textnormal{eq}}$.

\subsection{Numerical results}

In this section, we compute using PC the equilibrium $(x_{\textnormal{eq}},y_{\textnormal{eq}})$ and its stable and unstable manifold. For given $a$, $\bb$ and $\sigma$, we use the techniques presented in Section~\ref{sec:inv_sets} on
\begin{equation}
\label{eq:LV_vector_field}
f((x,y),p)=\begin{pmatrix}
(1-x-ay)x,  &
(1-y-(\bb+\sigma p)x)y
\end{pmatrix},
\end{equation}
where for convenience we use $(x,y)$ instead of $(x_1,x_2)$.

\subsubsection{Analysis of convergence on the steady state problem}
\label{sec:LV_eq}

We first solve, for various choices of expansions, the steady state problem $f(x(p),y(p),p)=0$, and analyze how the choice of expansion together with the distribution of $p$ affect the convergence rates. More precisely, for a given expansion basis $\phi_n$ and truncation level $N$, we look for coefficients $\left(x_n\right)_{0\leq n <N}$ and $\left(y_n\right)_{0\leq n <N}$ such that
\begin{equation*}
f(x^N(p),y^N(p),p)\approx 0,
\end{equation*}
where
\begin{equation*}
x^N(p)=\sum_{0\leq n< N}x_n \phi_n(p) \qquad \text{and} \qquad y^N(p)=\sum_{0\leq n< N}y_n \phi_n(p).
\end{equation*}
By taking the scalar product with respect to each $\phi_n$, $0\leq n<N$, we obtain the following system
\begin{equation}
\label{eq:eq_PC}
\left\{\begin{aligned}
&x_n-(x\ast x)_n-a(x\ast y)_n =0\\
&y_n-(y\ast y)_n-\bb y_n-\sigma (p\ast x\ast y)_n =0
\end{aligned}\right.
\qquad \forall~0\leq n<N.
\end{equation}
Here $\ast$ denotes the product associated to the basis $\phi_n$, i.e. $(x\ast y)_{0\leq n< 2N-1}$ is the unique sequence such that 
\begin{equation*}
x^N(p)y^N(p) = \sum_{0\leq n< 2N-1}(x\ast y)_n \phi_n(p).
\end{equation*}
If $\phi_n$ is chosen so that $\phi_1(p)=p$, which will often be the case in practice, then the expansion of $p$ only has a single non-zero coefficient: $p_n=\delta_{n,1}$.

\begin{remark}
We emphasize that $\ast$ of course depends on the choice of the basis $\phi_n$. We refer to the Appendix for more details about this product structure, and discussions on how it can be computed in practice.
\end{remark}

\begin{figure}[h]
\centering
\begin{subfigure}{0.32\linewidth}
\includegraphics[width=\linewidth]{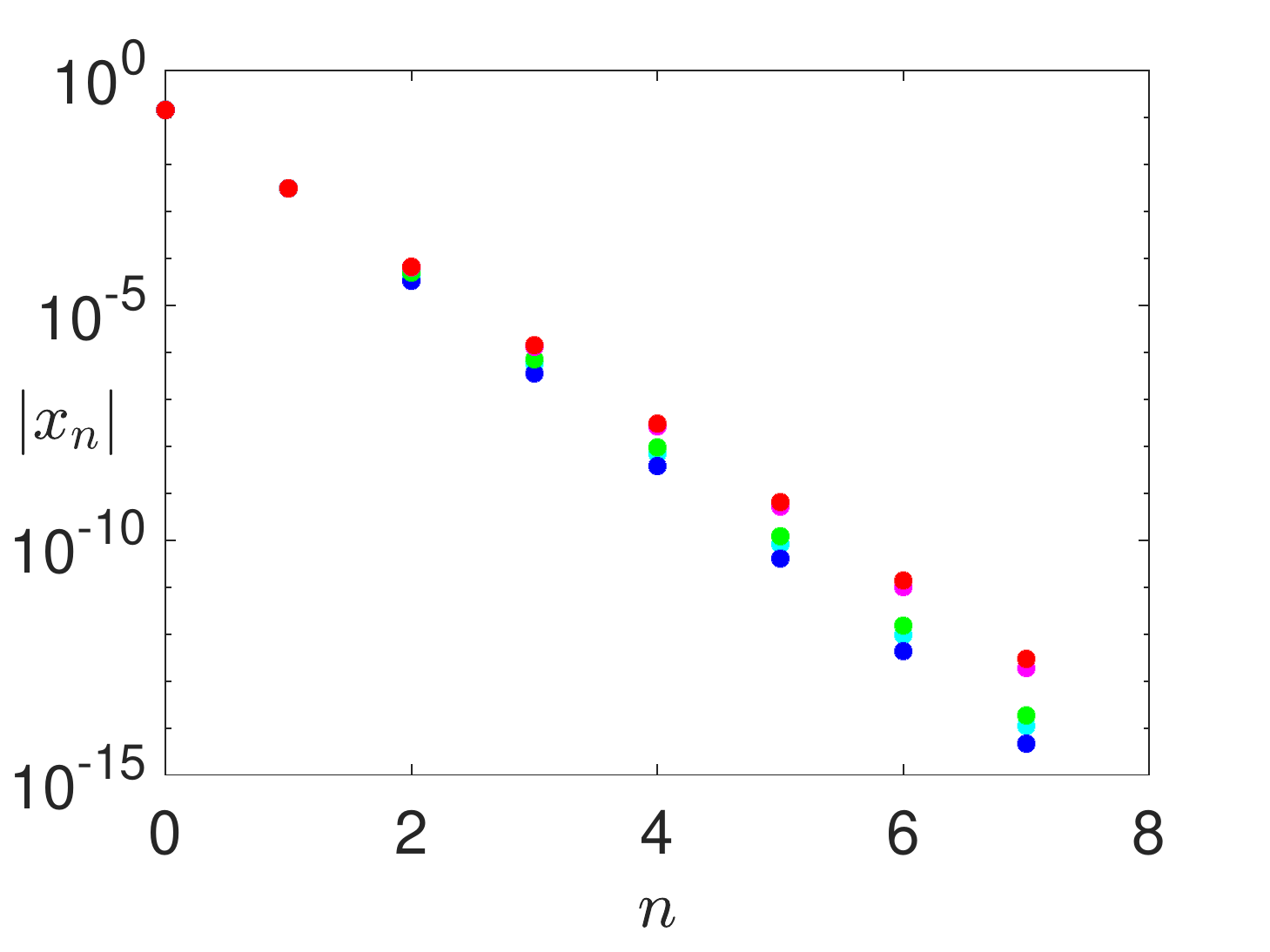}
\caption{$\sigma=0.1$}
\label{fig:decay_coeffs_0_1}
\end{subfigure}
\begin{subfigure}{0.32\linewidth}
\includegraphics[width=\linewidth]{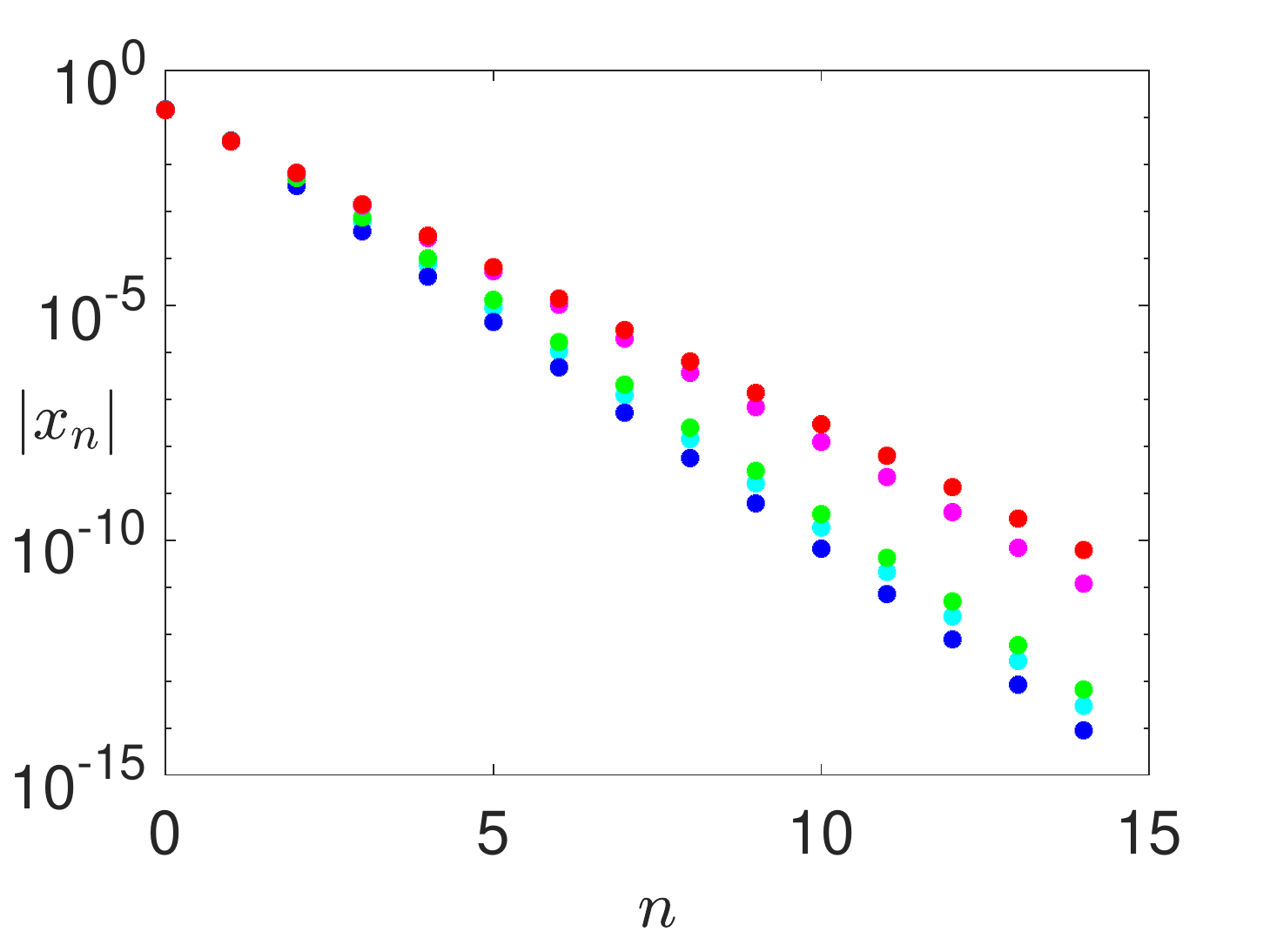}
\caption{$\sigma=1$}
\label{fig:decay_coeffs_1}
\end{subfigure}
\begin{subfigure}{0.32\linewidth}
\includegraphics[width=\linewidth]{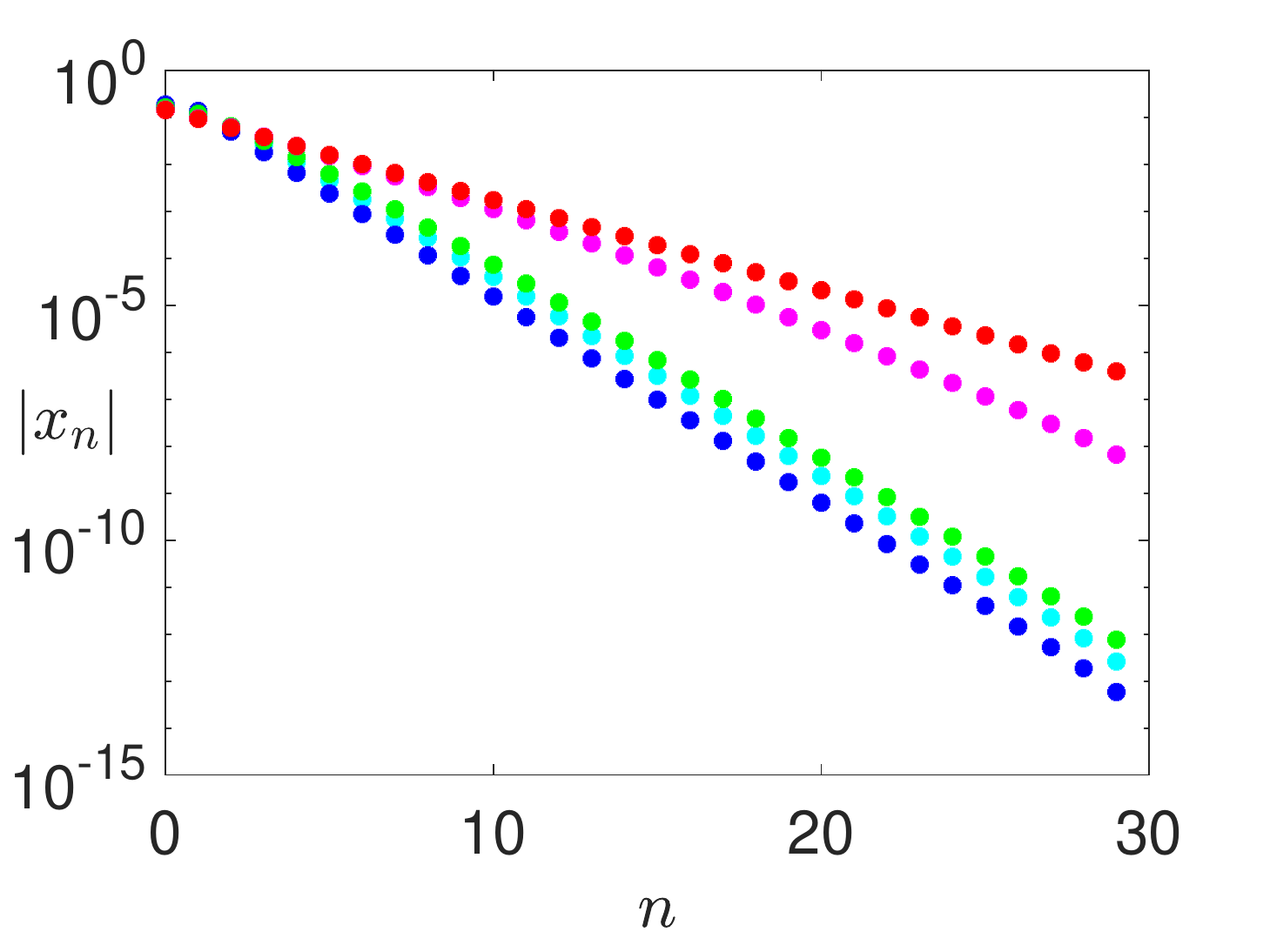}
\caption{$\sigma=3$}
\label{fig:decay_coeffs_3}
\end{subfigure}
\caption{We represent the absolute value of the coefficients $x_n$ solving~\eqref{eq:eq_PC} with respect to $n$. The computations were done for $a=3$, $\bb=5$ and $\sigma=0.1$ for Figure~\ref{fig:decay_coeffs_0_1}, $\sigma=1$ for Figure~\ref{fig:decay_coeffs_1} and $\sigma=3$ for Figure~\ref{fig:decay_coeffs_3}. In each case, the coefficients are computed for different expansions: Chebyshev of the first kind in blue, Legendre in cyan, Chebyshev of the second kind (rescaled) in green, Gegenbauer with $\mu=20$ (rescaled) in magenta and Taylor in red.}
\label{fig:decay_coeffs}
\end{figure}

We solve~\eqref{eq:eq_PC} using Newton's method, and the initial conditions given by the \emph{deterministic} equilibrium, i.e., for $\sigma=0$ so that
\begin{equation*}
x_0=\frac{a-1}{a\bb-1},\ y_0=\frac{\bb-1}{a\bb-1} \qquad \text{and}\qquad x_n=0,\ y_n=0 \quad \forall~1\leq n<N.
\end{equation*}
The results for different choices of basis $\phi_n$ are presented in Figure~\ref{fig:decay_coeffs}. To make the comparison fair, we consider rescaled version of the Chebyshev polynomials of the second kind and of the Gegenbauer polynomials (see the Appendix), so that all $\phi_n$ used satisfy $\sup_{s\in[-1,1]} \vert \phi_n(s)\vert =1$. The decay of the coefficients is therefore indicative of the truncation error in $\cC^0$-norm. As expected, the Chebyshev and Legendre expansions converge faster than the Taylor expansions, with the Gegenbauer one lying somewhere in between. As mentioned previously, we focus here only on the first component $x$, but the same behavior is observed for the second component $y$.

\begin{figure}[h]
\centering
\begin{subfigure}{0.49\linewidth}
\includegraphics[width=\linewidth]{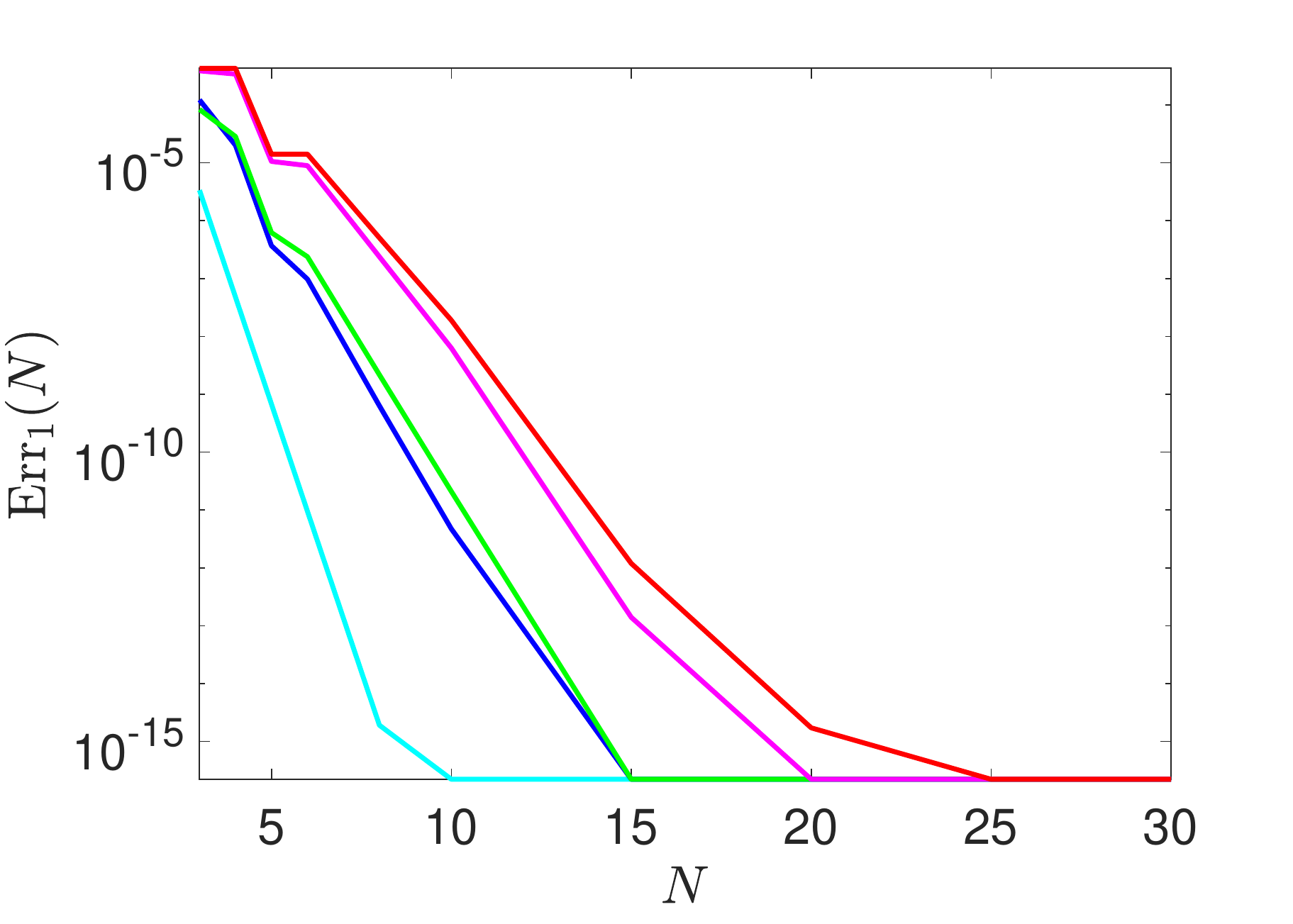}
\caption{$\rho_p(s)=\frac{1}{2}\1_{s\in(-1,1)}$}
\label{fig:err_xm1_uniform}
\end{subfigure}
\begin{subfigure}{0.49\linewidth}
\includegraphics[width=\linewidth]{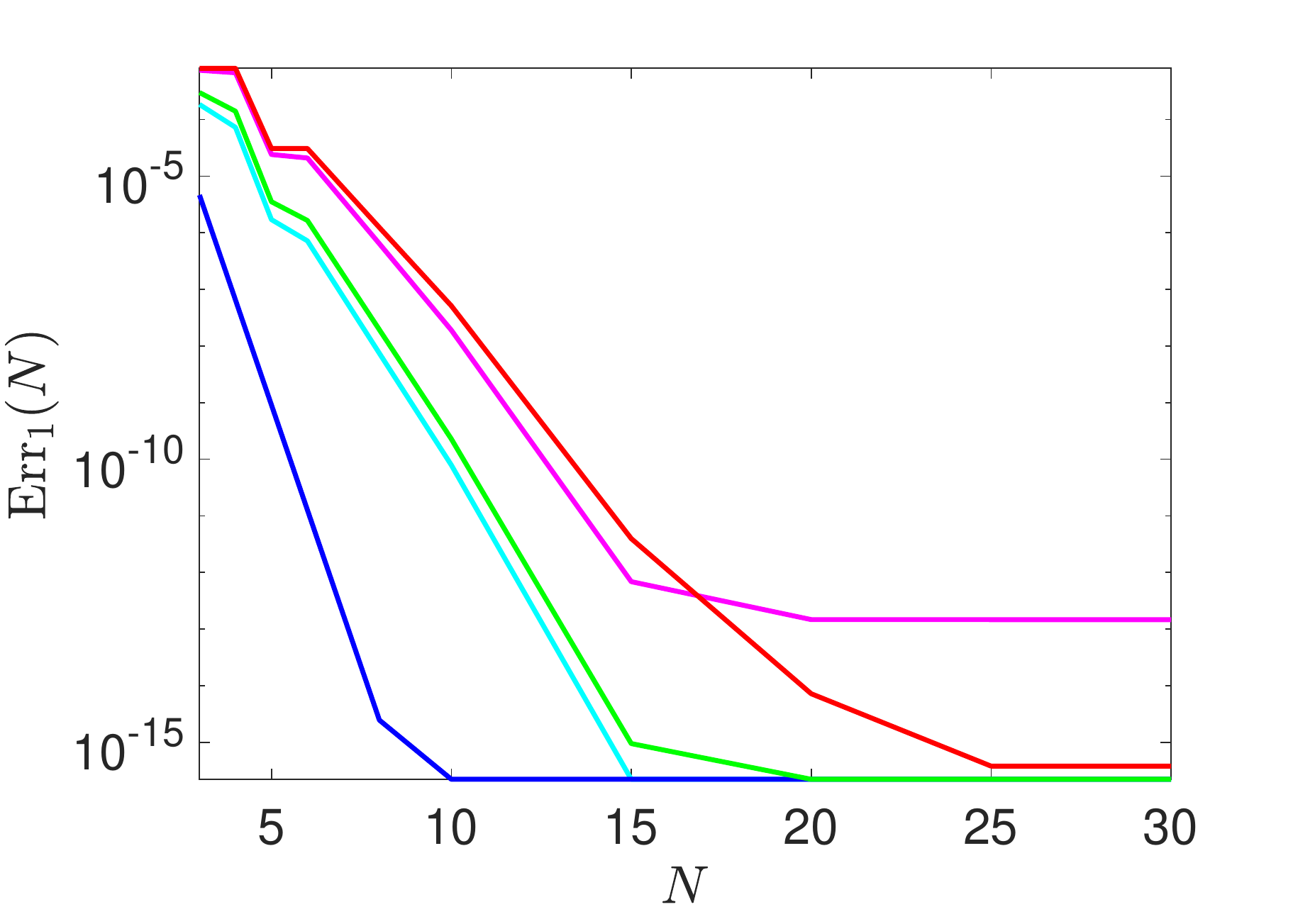}
\caption{$\rho_p(s)=\1_{s\in(-1,1)}\frac{1}{\pi\sqrt{1-s^2}}$}
\label{fig:err_xm1_betacheb1}
\end{subfigure}
\begin{subfigure}{0.49\linewidth}
\includegraphics[width=\linewidth]{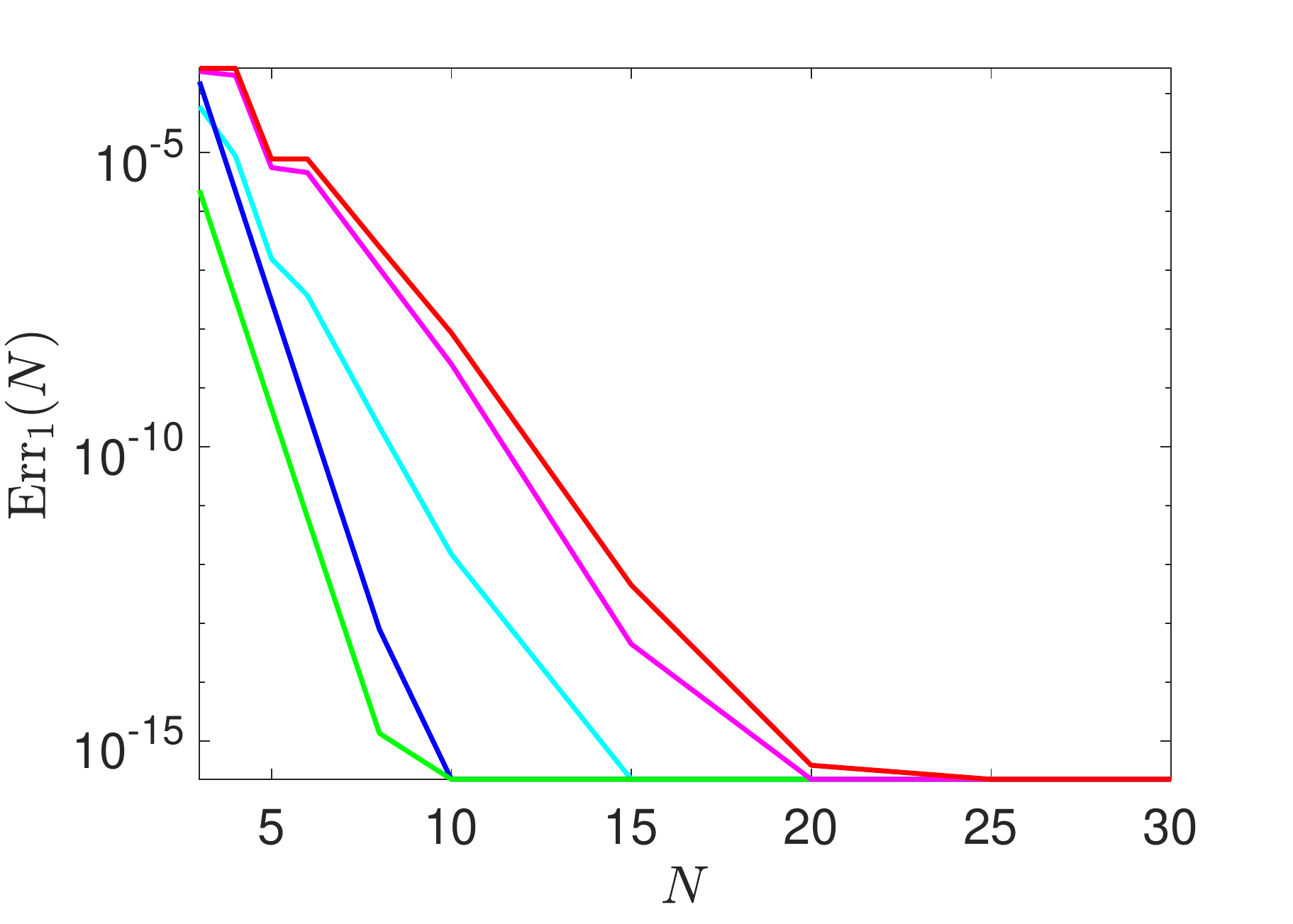}
\caption{$\rho_p(s)=\1_{s\in(-1,1)}\frac{2}{\pi}\sqrt{1-s^2}$}
\label{fig:err_xm1_betacheb2}
\end{subfigure}
\begin{subfigure}{0.49\linewidth}
\includegraphics[width=\linewidth]{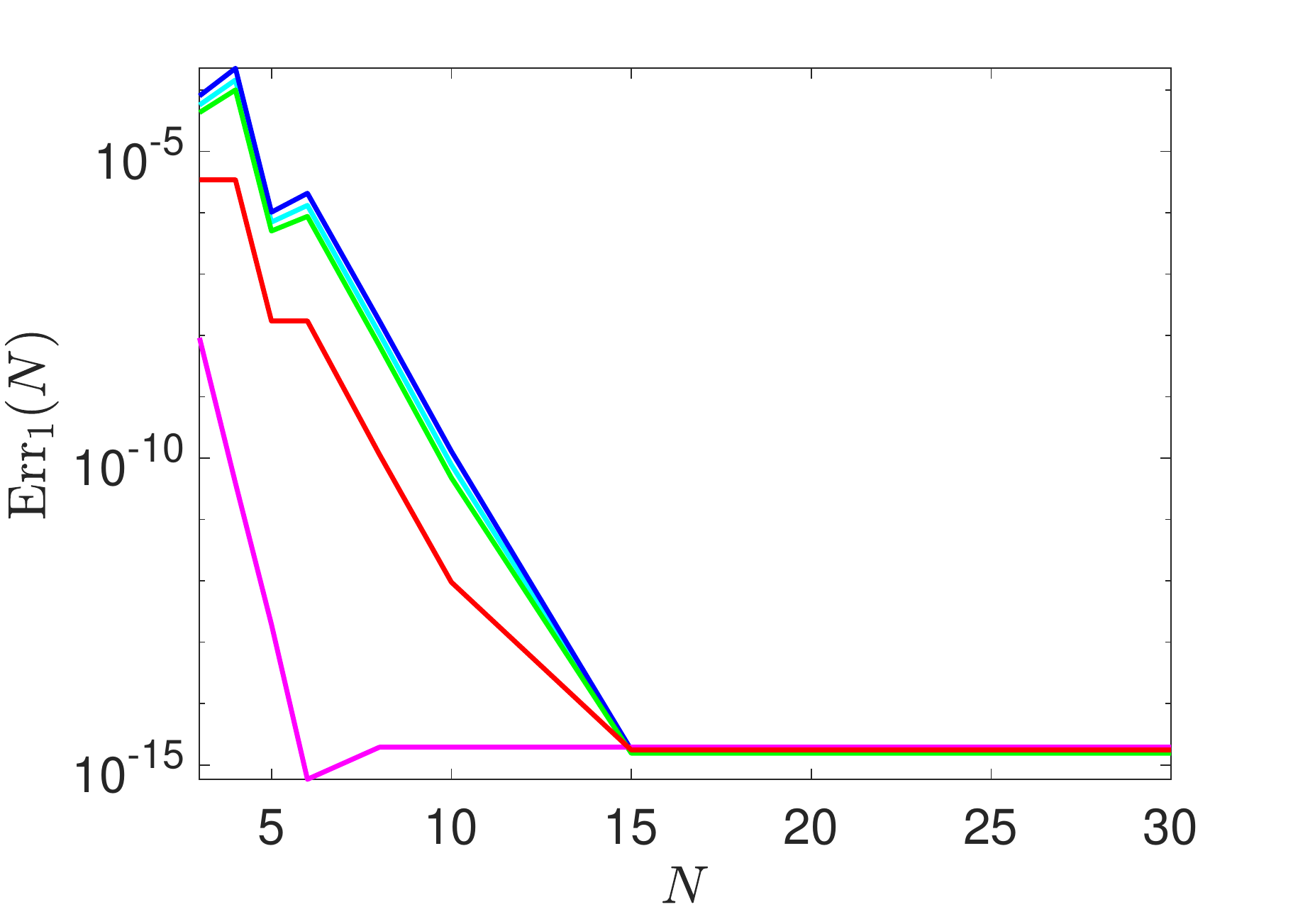}
\caption{$\rho_p(s)=\1_{s\in(-1,1)}\frac{2^{2\mu-1}\mu B(\mu,\mu)}{\pi}(1-s^2)^{\mu-\frac{1}{2}}$, with $\mu=20$}
\label{fig:err_xm1_betasym20}
\end{subfigure}
\caption{Relative error for the first moment in function of the truncation level $N$, for several distributions $\rho_p$ of the parameter $p$. In each case, we take $a=3$, $\bb=5$, $\sigma=1$, and use several expansions which are represented in different colors: Chebyshev of the first kind in blue, Legendre in cyan, Chebyshev of the second kind in green, Gegenbauer with $\mu=20$ in magenta and Taylor in red.}
\label{fig:err_xm1}
\end{figure}

\begin{figure}[h]
\centering
\begin{subfigure}{0.49\linewidth}
\includegraphics[width=\linewidth]{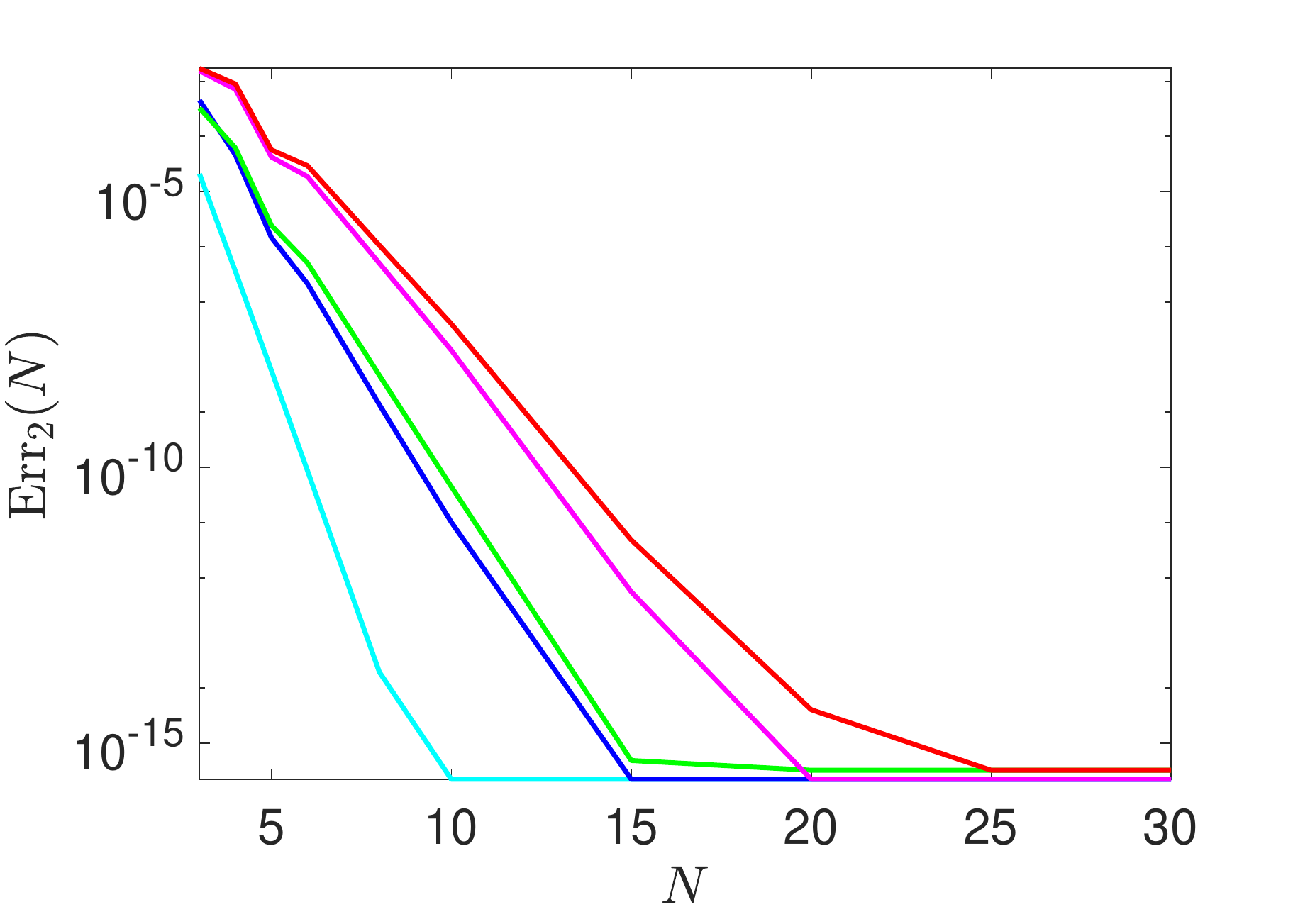}
\caption{$\rho_p(s)=\frac{1}{2}\1_{s\in(-1,1)}$}
\label{fig:err_xm2_uniform}
\end{subfigure}
\begin{subfigure}{0.49\linewidth}
\includegraphics[width=\linewidth]{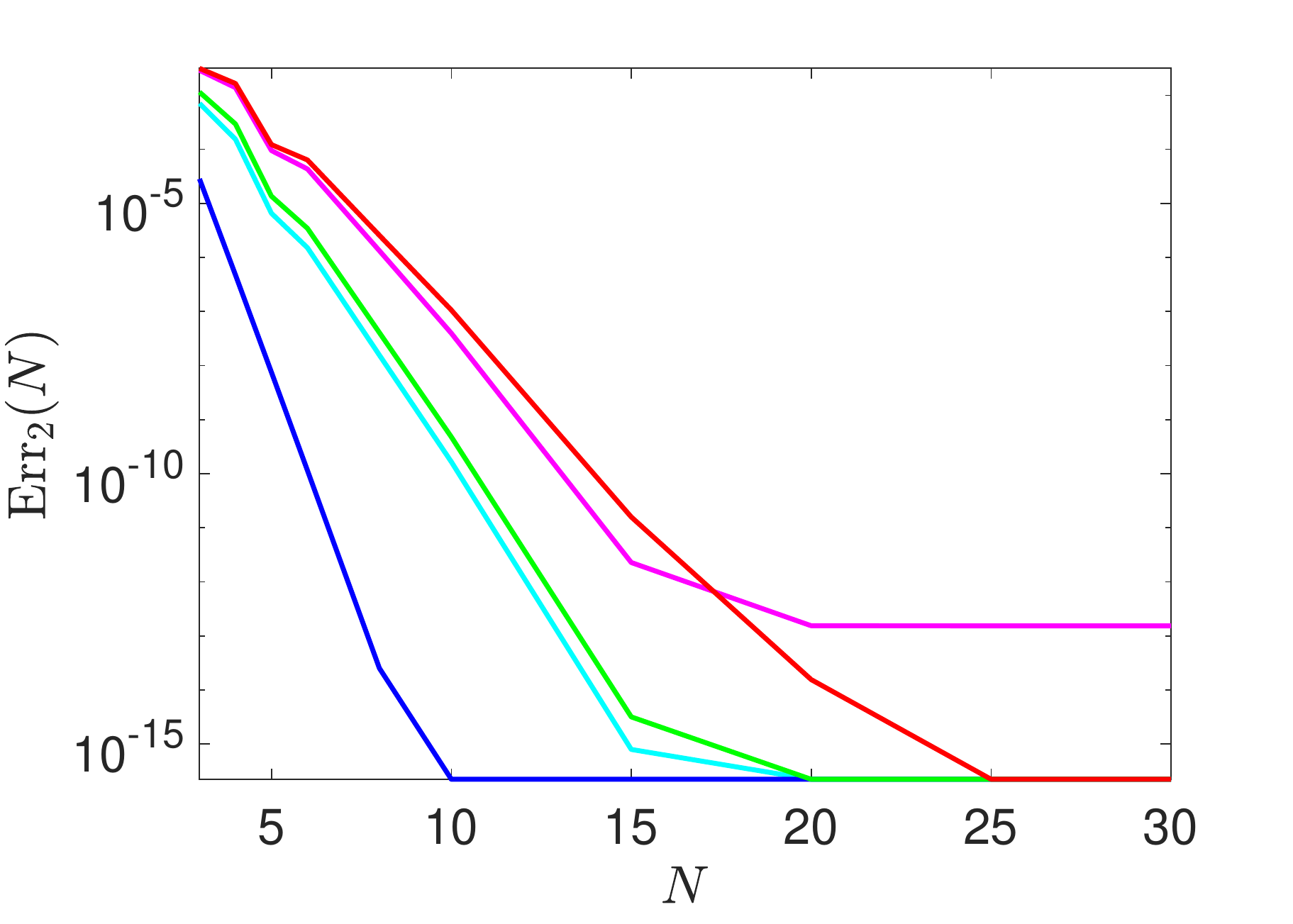}
\caption{$\rho_p(s)=\1_{s\in(-1,1)}\frac{1}{\pi\sqrt{1-s^2}}$}
\label{fig:err_xm2_betacheb1}
\end{subfigure}
\begin{subfigure}{0.49\linewidth}
\includegraphics[width=\linewidth]{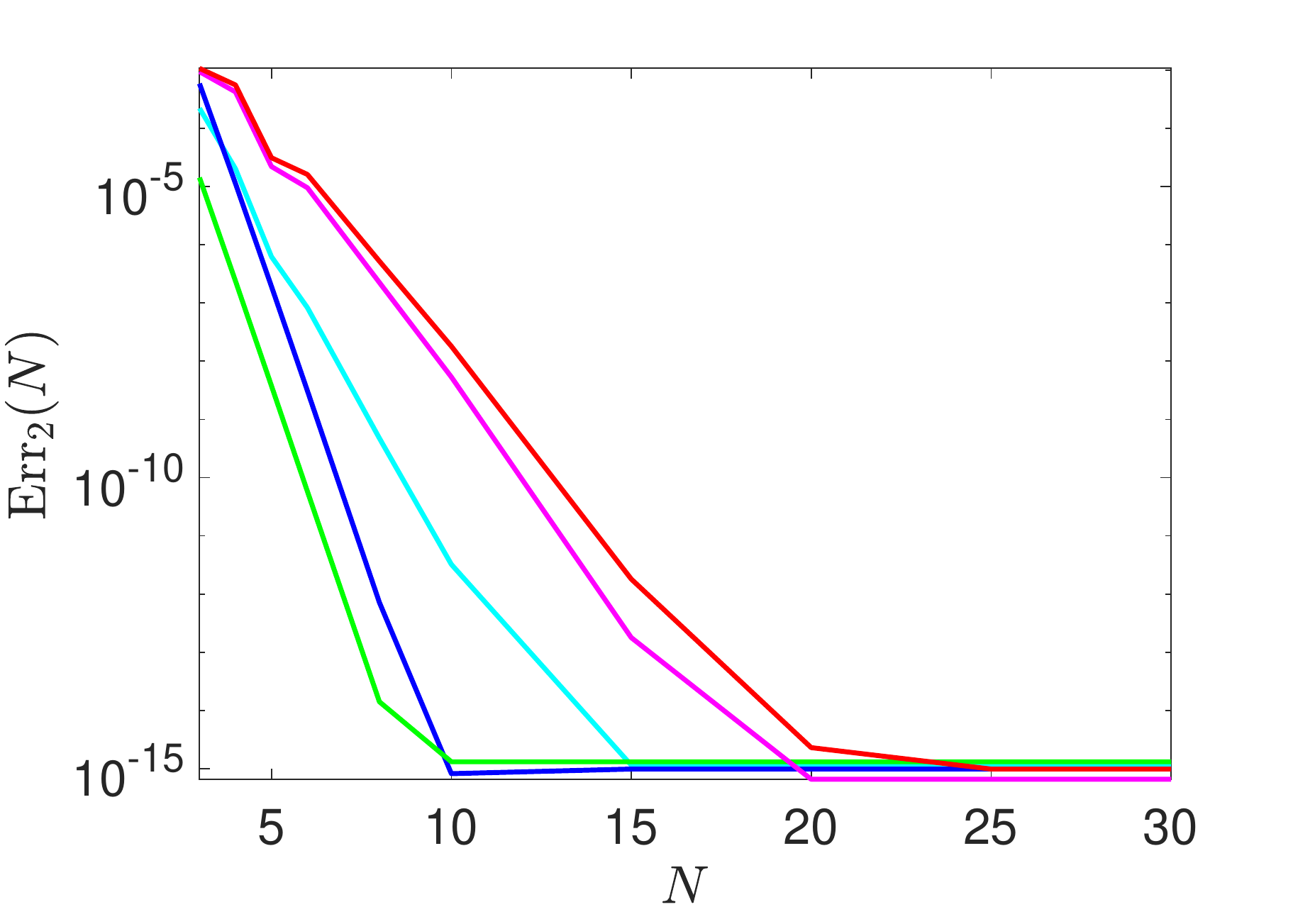}
\caption{$\rho_p(s)=\1_{s\in(-1,1)}\frac{2}{\pi}\sqrt{1-s^2}$}
\label{fig:err_xm2_betacheb2}
\end{subfigure}
\begin{subfigure}{0.49\linewidth}
\includegraphics[width=\linewidth]{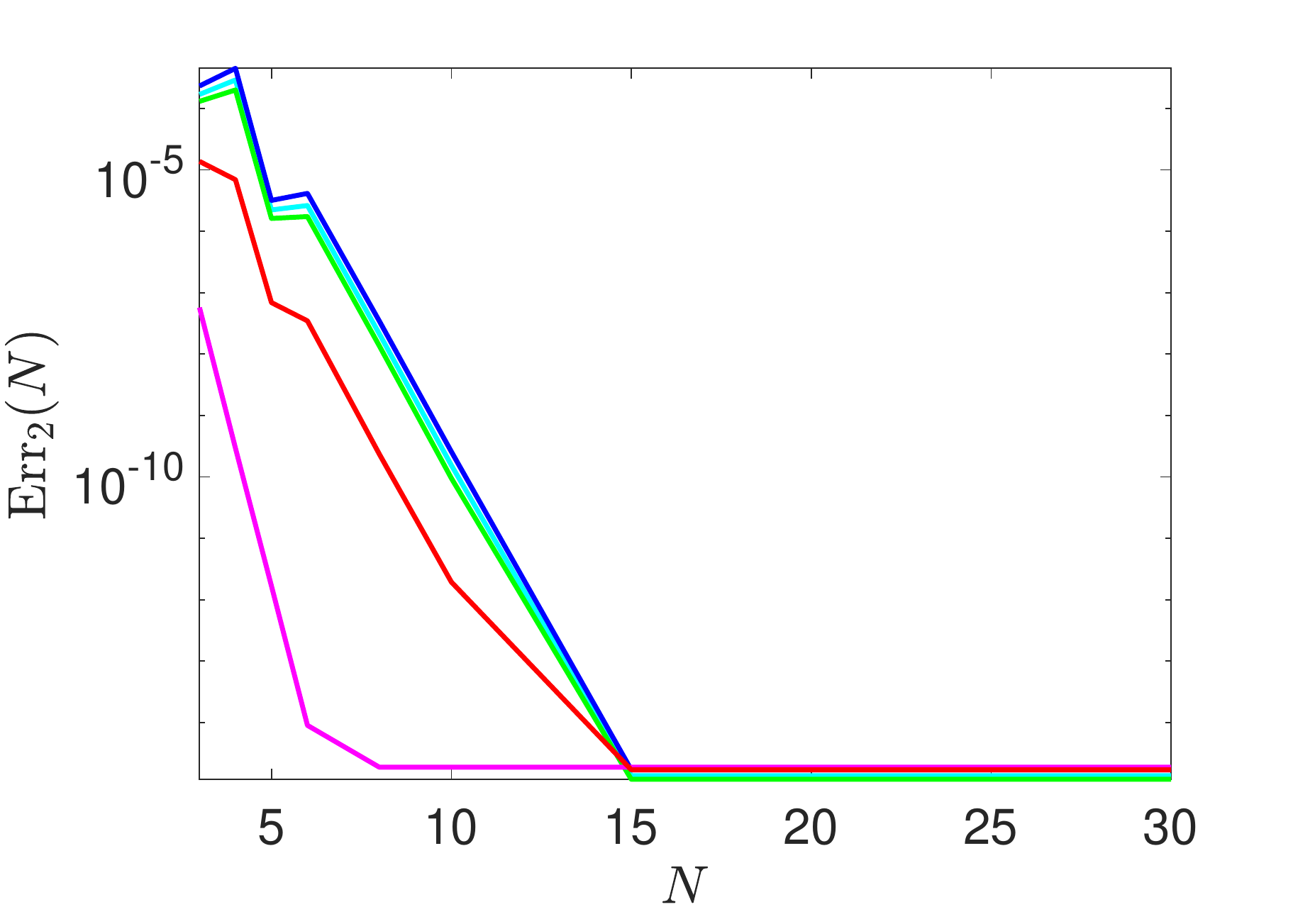}
\caption{$\rho_p(s)=\1_{s\in(-1,1)}\frac{2^{2\mu-1}\mu B(\mu,\mu)}{\pi}(1-s^2)^{\mu-\frac{1}{2}}$, with $\mu=20$}
\label{fig:err_xm2_betasym20}
\end{subfigure}
\caption{Relative error for the second moment in function of the truncation level $N$, for several distributions $\rho_p$ of the parameter $p$. In each case, we take $a=3$, $\bb=5$, $\sigma=1$, and use several expansions which are represented in different colors: Chebyshev of the first kind in blue, Legendre in cyan, Chebyshev of the second kind in green, Gegenbauer with $\mu=20$ in magenta and Taylor in red.}
\label{fig:err_xm2}
\end{figure}

However, as mentioned in the introduction, we are more interested in understanding the distribution of the solution, assuming $p$ has a prescribed distribution in $[-1,1]$. In Figure~\ref{fig:err_xm1} and Figure~\ref{fig:err_xm2} we display the relative error for the first two moments
\begin{equation*}
\textnormal{Err}_1(N)=\frac{\left\vert \E(x(p)) - \E(x^N(p)) \right\vert}{\left\vert \E(x(p)) \right\vert} \qquad \text{and} \qquad \textnormal{Err}_2(N)=\frac{\left\vert \E(x(p)^2) - \E(x^N(p)^2) \right\vert}{\left\vert \E(x(p)^2) \right\vert},
\end{equation*}
depending on the truncation level $N$, for several PDF $\rho_p$ of $p$ and again for several expansions. 

The limiting moments $\E(x(p))$ and $\E(x(p)^2)$ are computed analytically when possible (see Section~\ref{sec:LV_sto}) and using numerical integration otherwise. The truncated moments $\E(x^N(p))$ and $\E(x^N(p)^2)$ are computed from the coefficients of the expansion
\begin{equation*}
x^N = \sum_{n=0}^{N-1} x_n \phi_n.
\end{equation*}
If $\phi_n$ is the orthogonal basis associated to the distribution $\rho_p$, then as explained in Section~\ref{sec:PC}
\begin{equation*}
\E(x^N(p)) = x_0 \qquad \text{and} \qquad \E(x^N(p)^2) = \sum_{n=0}^{N-1} x_n^2 h_n.
\end{equation*}
Otherwise, we have
\begin{equation*}
\E(x^N(p)) = \sum_{n=0}^{N-1} x_n \int_{-1}^1 \phi_n(s) \rho_p(s) \txtd s \quad \text{and} \quad \E(x^N(p)) = \sum_{n=0}^{2N-2} (x\ast x)_n \int_{-1}^1 \phi_n(s) \rho_p(s) \txtd s,
\end{equation*}
and we first compute the above integrals to obtain the truncated moments. As expected from the paradigm of PC, we obtain faster convergence for these statistics when using the PC expansions associated to the distribution $\rho_p$ of the parameter $p$, i.e., a family of polynomial $\phi_n$ orthogonal with respect to $\left\langle \cdot,\cdot \right\rangle_{\rho_p}$. In particular, if $p$ has a uniform distribution then the Legendre expansion converges the fastest (Figure~\ref{fig:err_xm1_uniform} and Figure~\ref{fig:err_xm2_uniform}), if $p$ has an arcsine distribution then the Chebyshev expansion of the first kind converges the fastest (Figure~\ref{fig:err_xm1_betacheb1} and Figure~\ref{fig:err_xm2_betacheb1}), if $p$ has a Wigner semicircle distribution then the Chebyshev expansion of the second kind converges the fastest (Figure~\ref{fig:err_xm1_betacheb2} and Figure~\ref{fig:err_xm2_betacheb2}) and finally if $p$ has a beta distribution of parameter $(20,20)$ then the Gegenbauer expansion of parameter $\mu=20$ converges the fastest  (Figure~\ref{fig:err_xm1_betasym20} and Figure~\ref{fig:err_xm2_betasym20}).

\begin{remark}
These simple experiments confirm that it is in general a good option to choose the expansion basis associated to the distribution of the parameter, and we will systematically do so in the sequel. Nonetheless, we believe that in some cases, especially highly nonlinear ones, the cost of actually computing the nonlinear terms should also be taken into account (see the discussion in Section~\ref{sec:Appendix_fast}).
\end{remark}

\subsubsection{Computation of eigenvalues and eigenvectors}
\label{sec:LV_eig}

Now that we found an accurate representation of the steady state, we turn our attention to the dynamics around it. Before computing parameterizations of the local stable and unstable manifolds, which will be done in the next subsection, we first focus on the linearized system around the equilibrium. More precisely, we are interested in the usual eigenvalue problem
\begin{equation}
\label{eq:eigenproblem}
\left\{\begin{aligned}
&\txtD_{(x,y)}f(x(p),y(p),p) \begin{pmatrix} u(p) \\ v(p) \end{pmatrix} -\lambda(p) \begin{pmatrix} u(p) \\ v(p) \end{pmatrix} = 0 \\
& u(p)^2+v(p)^2=1, 
\end{aligned}\right.
\end{equation}
where $(x(p),y(p))$ is the previously computed equilibrium, and we solve for $u(p)$, $v(p)$ and $\lambda(p)$, the last equation being a normalization of the eigenvector allowing us to have a (locally) unique solution. In the previous subsection we obtained truncated expansions $x^N$ and $y^N$ of the equilibrium, and we now aim at doing the same for the eigendata. Namely, we write
\begin{equation*}
u^N(p)=\sum_{0\leq n<N}u_n \phi_n(p), \qquad v^N(p)=\sum_{0\leq n<N}v_n \phi_n(p) \qquad \text{and}\qquad \lambda^N(p)=\sum_{0\leq n<N}\lambda_n \phi_n(p),
\end{equation*}
and plug these expansions back in~\eqref{eq:eigenproblem}. For our explicit example~\eqref{eq:LV_vector_field}, this yields the following system for the coefficients $u_n$, $v_n$ and $\lambda_n$:
\begin{equation*}
\left\{\begin{aligned}
& u_n -2(x\ast u)_n-a(y\ast u)_n -a(x\ast v)_n - (\lambda\ast u)_n = 0 \\
& -\bb(y\ast u)_n -\sigma (p\ast y\ast u)_n +v_n -2(y\ast v)_n -\bb(x\ast v)_n -\sigma(p\ast x\ast v)_n -(\lambda\ast v)_n =0 \\
& (u\ast u)_n +(v\ast v)_n - \delta_{n,0} = 0
\end{aligned}\right. \qquad \forall~0\leq n <N.
\end{equation*}
Again we solve this system using Newton's method and the initial conditions given by numerically obtained eigenvalue and eigenvector for the deterministic problem~\eqref{eq:eigenval_LV}-\eqref{eq:eigenvect_LV}. We obtain PC expansions for the eigenvalues $\lambda^{(1)}$ and $\lambda^{(2)}$ and associated eigenvectors $V_{\lambda^{(1)}}=\left(u^{(1)},v^{(1)}\right)$ and $V_{\lambda^{(2)}}=\left(u^{(2)},v^{(2)}\right)$.

\begin{figure}[htbp]
\centering
\begin{subfigure}{0.49\linewidth}
\includegraphics[width=\linewidth]{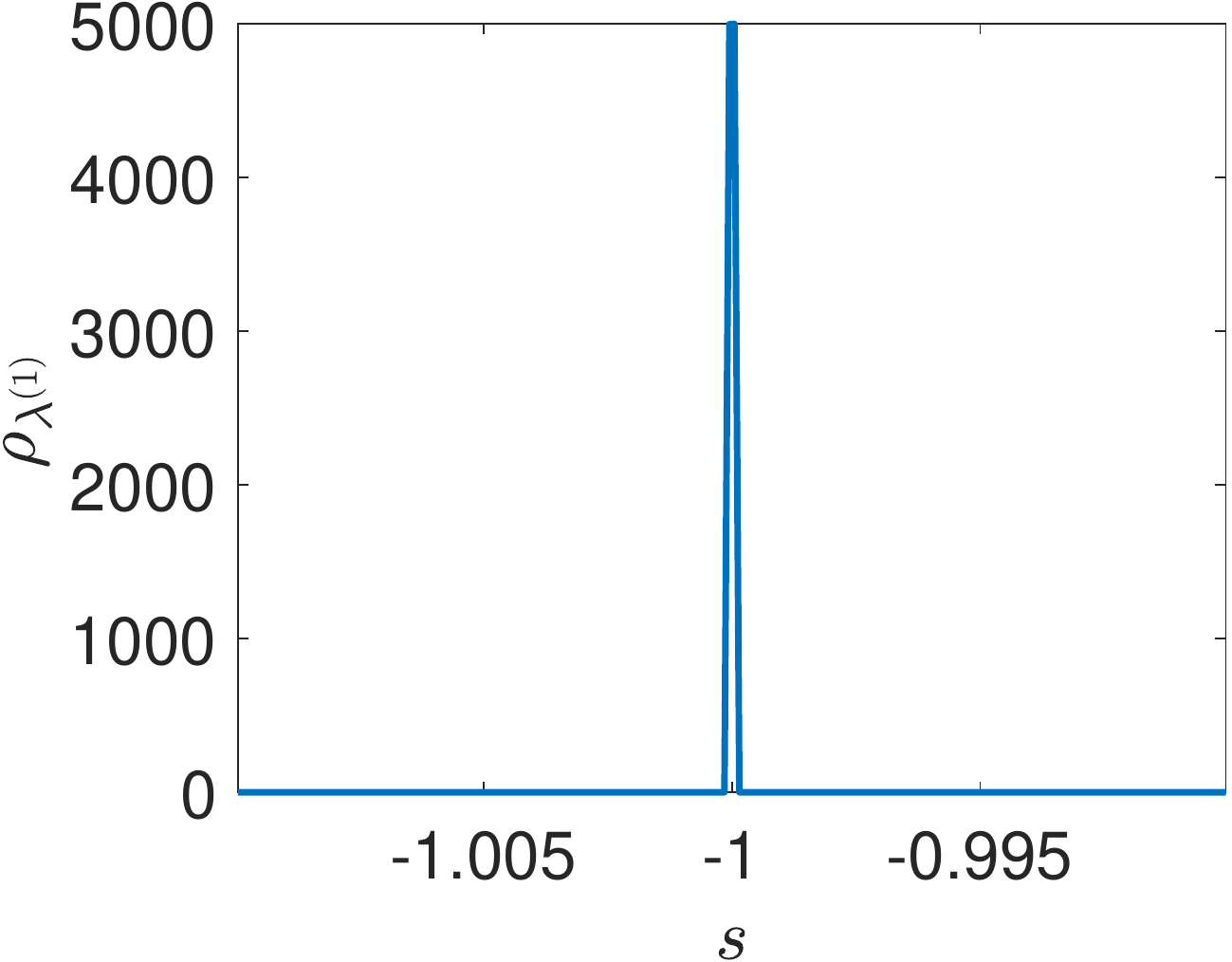}
\caption{PDF of $\lambda^{(1)}$ for $\rho_p(s)=\1_{s\in(-1,1)}\frac{1}{\pi\sqrt{1-s^2}}$.}
\label{fig:density_lambda_s_BetaCheb1}
\end{subfigure}
\hfill
\begin{subfigure}{0.49\linewidth}
\includegraphics[width=\linewidth]{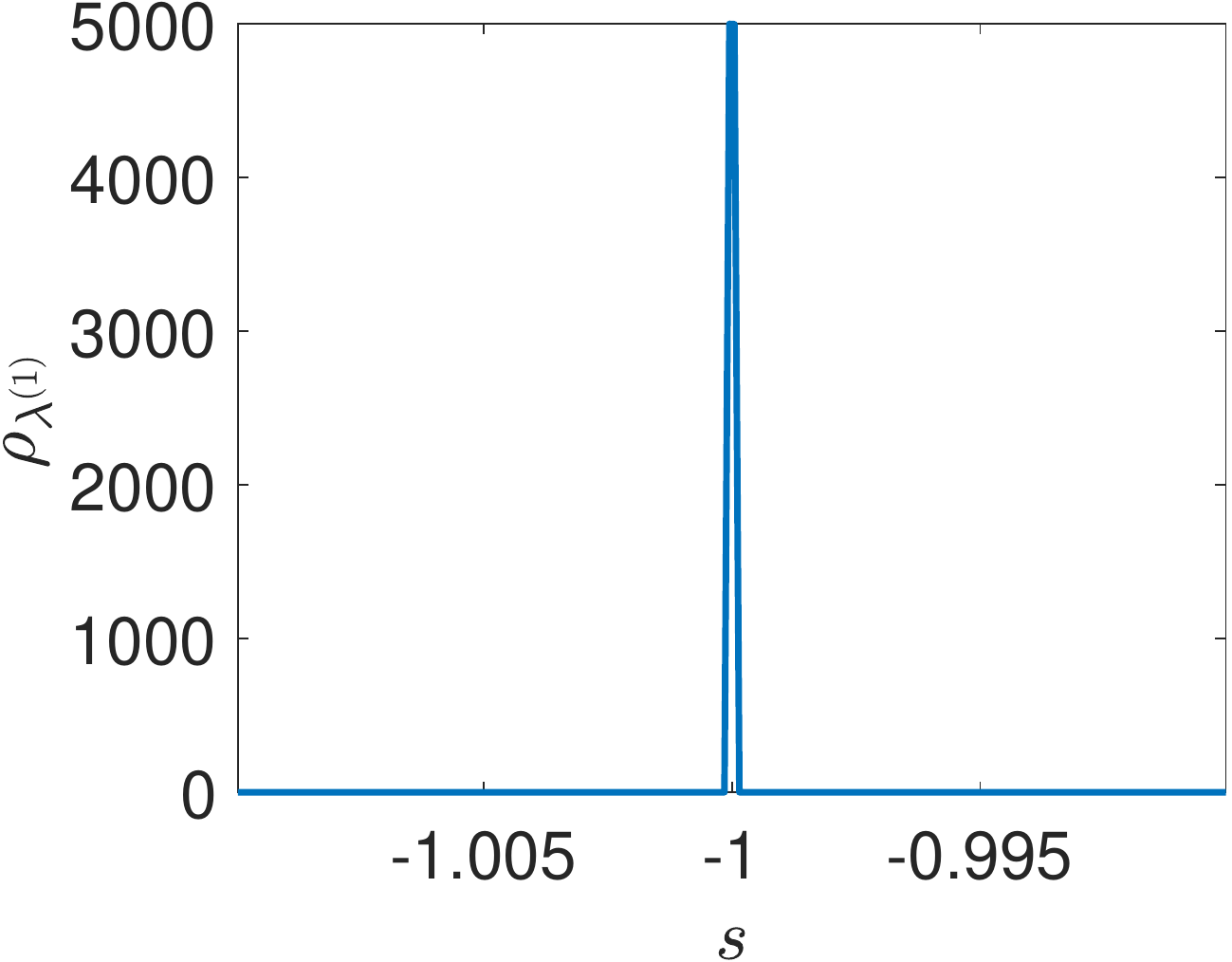}
\caption{PDF of $\lambda^{(1)}$ for $\rho_p(s)=\1_{s\in(-1,1)}\frac{2^{2\mu-1}\mu B(\mu,\mu)}{\pi}(1-s^2)^{\mu-\frac{1}{2}}$, with $\mu=20$.}
\label{fig:density_lambda_s_BetaSymmu}
\end{subfigure}\\
\vspace*{0.5cm}
\begin{subfigure}{0.49\linewidth}
\includegraphics[width=\linewidth]{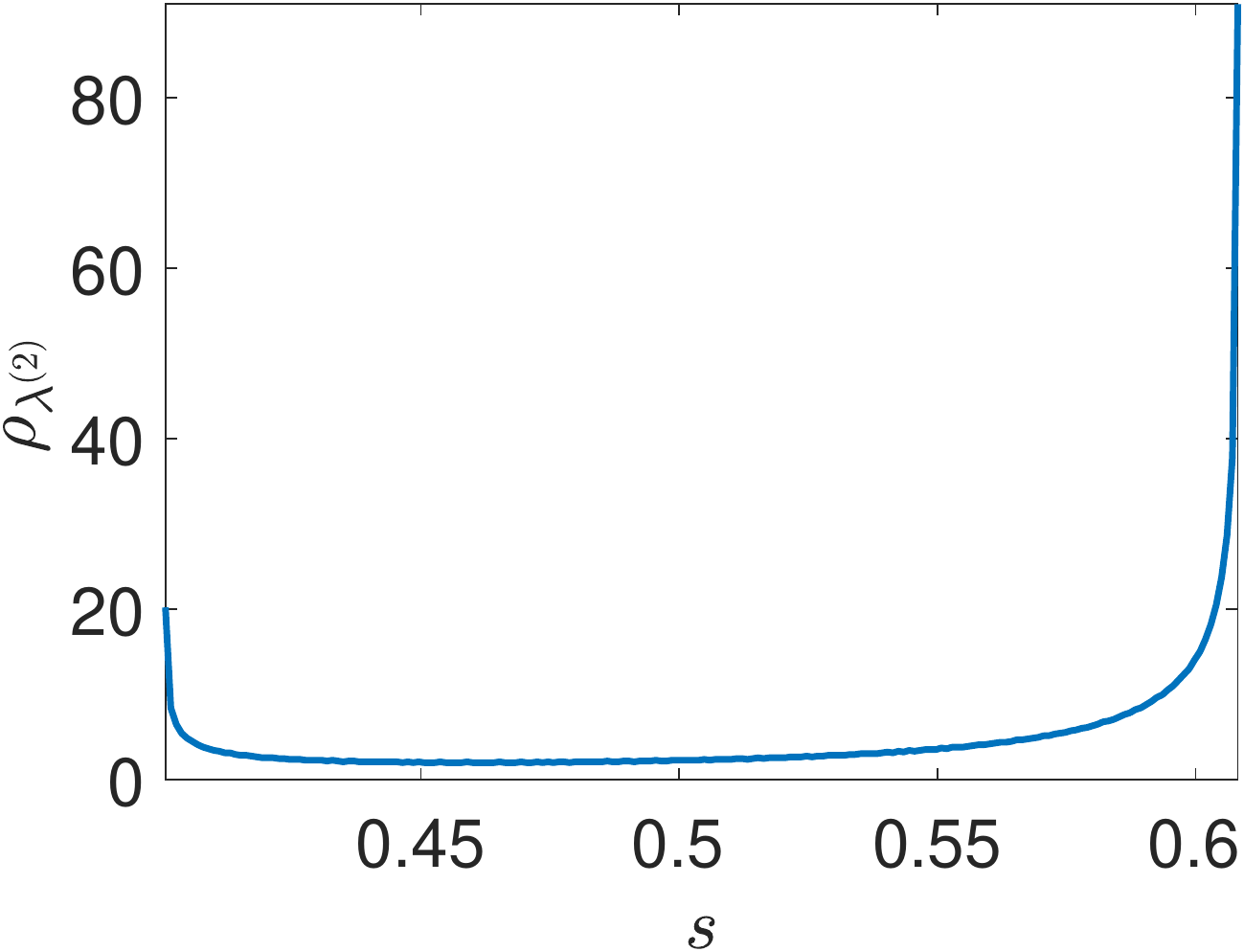}
\caption{PDF of $\lambda^{(2)}$ for $\rho_p(s)=\1_{s\in(-1,1)}\frac{1}{\pi\sqrt{1-s^2}}$.}
\label{fig:density_lambda_u_BetaCheb1}
\end{subfigure}
\hfill
\begin{subfigure}{0.49\linewidth}
\includegraphics[width=\linewidth]{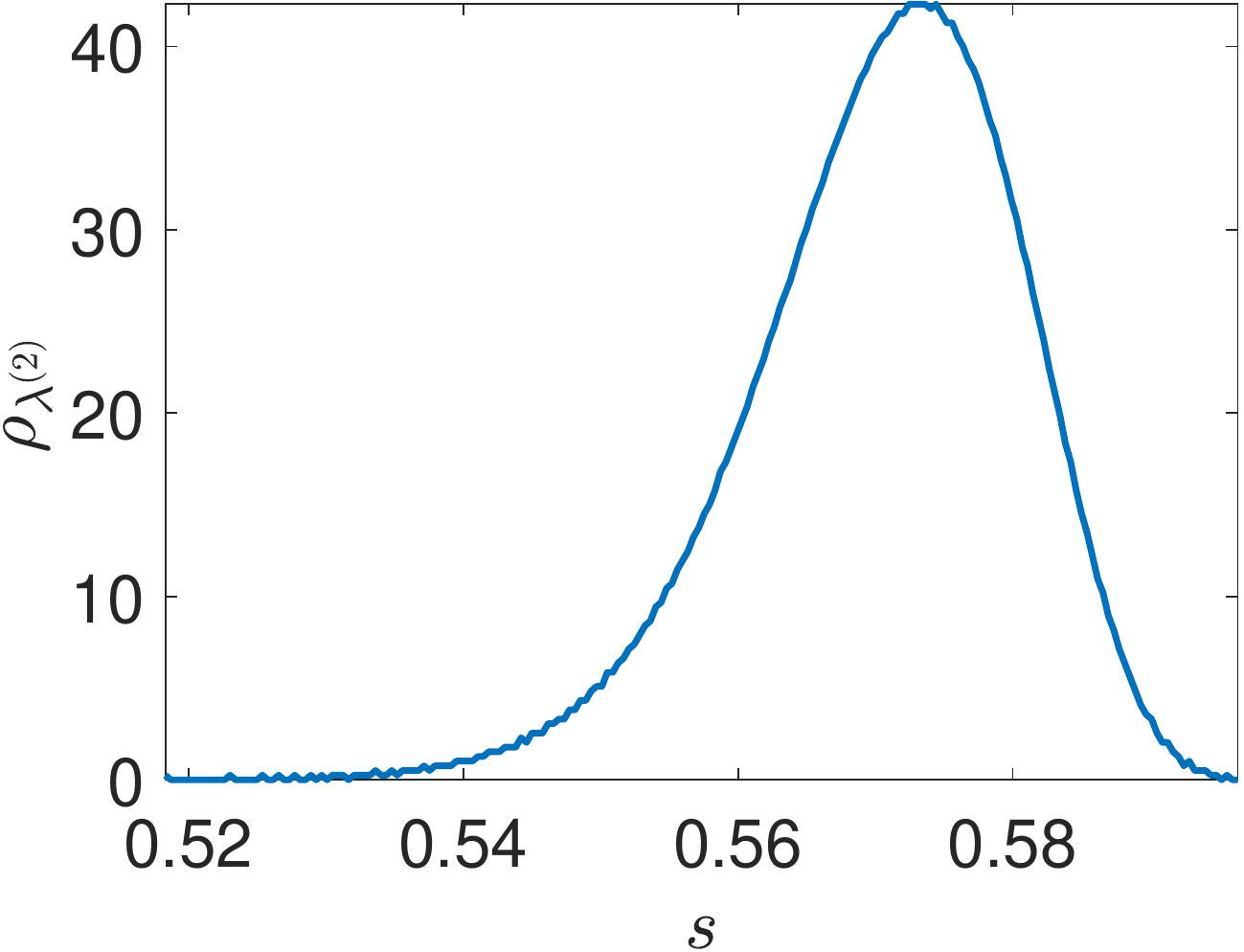}
\caption{PDF of $\lambda^{(2)}$ for $\rho_p(s)=\1_{s\in(-1,1)}\frac{2^{2\mu-1}\mu B(\mu,\mu)}{\pi}(1-s^2)^{\mu-\frac{1}{2}}$, with $\mu=20$.}
\label{fig:density_lambda_u_BetaSymmu}
\end{subfigure}
\caption{PDF of both eigenvalues, for different distributions $\rho_p$ of the parameter $p$. The PC coefficients of $\lambda^{(1)}$ and $\lambda^{(2)}$ were obtained for $a=3$, $\bb=5$, $\sigma=3$, and then used to compute $10,000$ samples used to recover the PDF. Notice that the results are consistent with the analytic computations~\eqref{eq:eigenval_LV}, as the PDF of $\lambda^{(1)}$ is supposed to be a Dirac at $-1$.}
\label{fig:density_lambda}
\end{figure}

\begin{figure}[h]
\centering
\begin{subfigure}{0.49\linewidth}
\includegraphics[width=\linewidth]{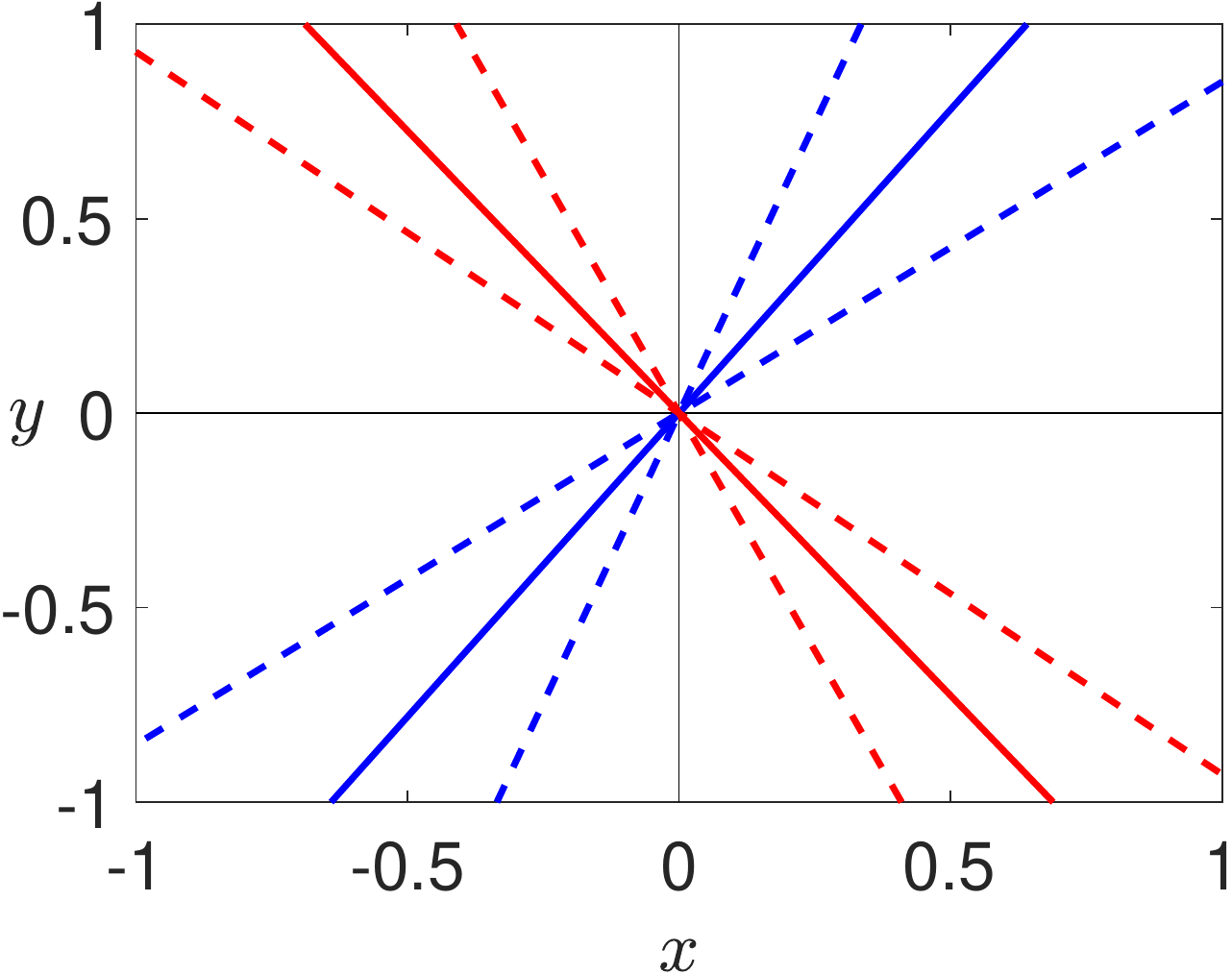}
\caption{$\rho_p(s)=\1_{s\in(-1,1)}\frac{1}{\pi\sqrt{1-s^2}}$.}
\label{fig:eigenspaces_tubes_BetaCheb1}
\end{subfigure}
\hfill
\begin{subfigure}{0.49\linewidth}
\includegraphics[width=\linewidth]{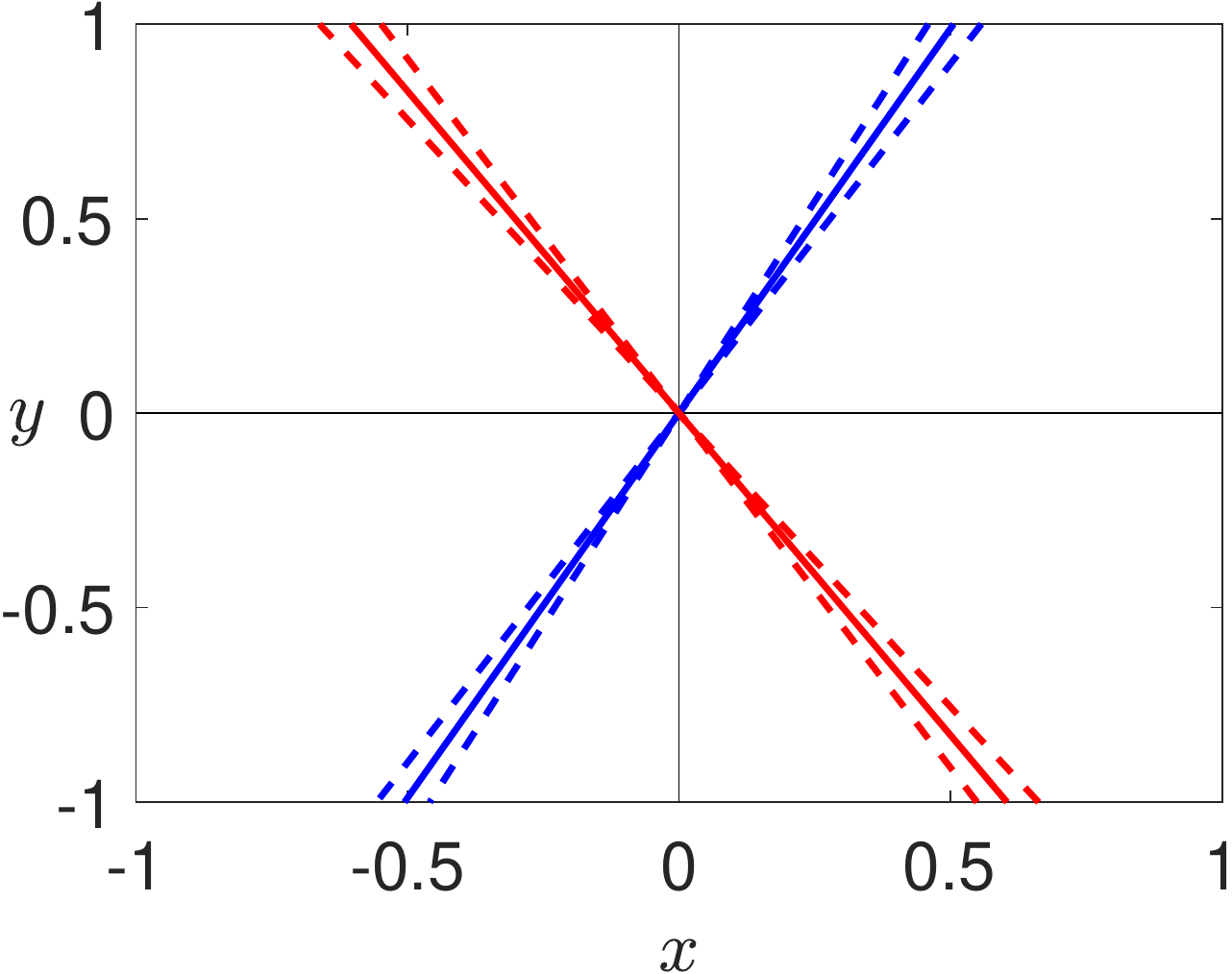}
\caption{$\rho_p(s)=\1_{s\in(-1,1)}\frac{2^{2\mu-1}\mu B(\mu,\mu)}{\pi}(1-s^2)^{\mu-\frac{1}{2}}$, with $\mu=20$.}
\label{fig:eigenspaces_tubes_BetaSymmu}
\end{subfigure}
\caption{We display here in solid lines the mean direction of the two eigenspaces (the stable one in blue and the unstable in red), given by $\left(\E(u^{(i)}),\E(v^{(i)})\right)$ for $i=1,2$, still for $a=3$, $\bb=5$ and $\sigma=3$. The dotted lines have directions given by $\left(\E(u^{(i)})\pm \sqrt{\V(u^{(i)})},\E(v^{(i)})\pm \sqrt{\V(v^{(i)})}\right)$.}
\label{fig:tubes_vects}
\end{figure}

From these PC expansions, we can then easily get statistical information about the eigenvalues and eigenvectors. For instance, as mentioned in the introduction, once the PC coefficients have been computed one can carry out Monte Carlo simulations basically for free (the only cost being the evaluation of the basis polynomials $\phi_n$, for $0\leq n<N$, at the sampled values of $p$). For a given distribution of $p$, this allows us to numerically compute the PDF of $\lambda$ (see Figure~\ref{fig:density_lambda}). Another approach, maybe more adapted to quantities of interest that are more than one-dimensional, is to use the mean and variance that can be obtained from the PC expansion to get a sense of where objects may lie in phase space (see Figure~\ref{fig:tubes_vects}). Since the PC expansion also allows us to compute higher order moments (see Section~\ref{sec:higher_order_terms}), this geometrical description can be made more quantitative by using these moments to obtain concentration inequalities.

\subsubsection{Computation of local stable/unstable manifolds}
\label{sec:LV_manifolds}

\begin{figure}[h]
\centering
\begin{subfigure}{0.49\linewidth}
\includegraphics[width=\linewidth]{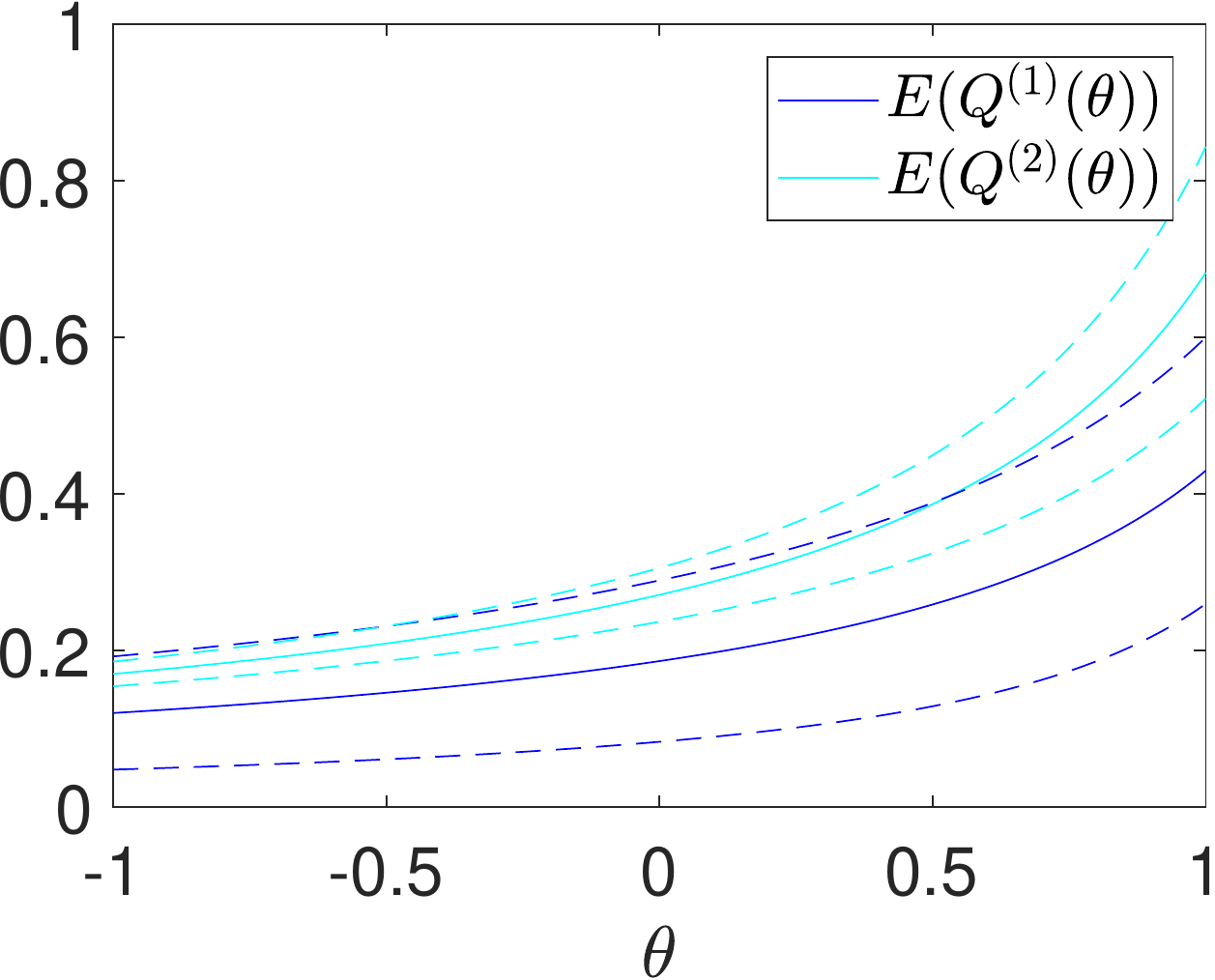}
\caption{Average position and standard deviation of the local stable manifold with respect to the parameter $\theta$.}
\label{fig:tube_u_BetaCheb1}
\end{subfigure}
\hfill
\begin{subfigure}{0.49\linewidth}
\includegraphics[width=\linewidth]{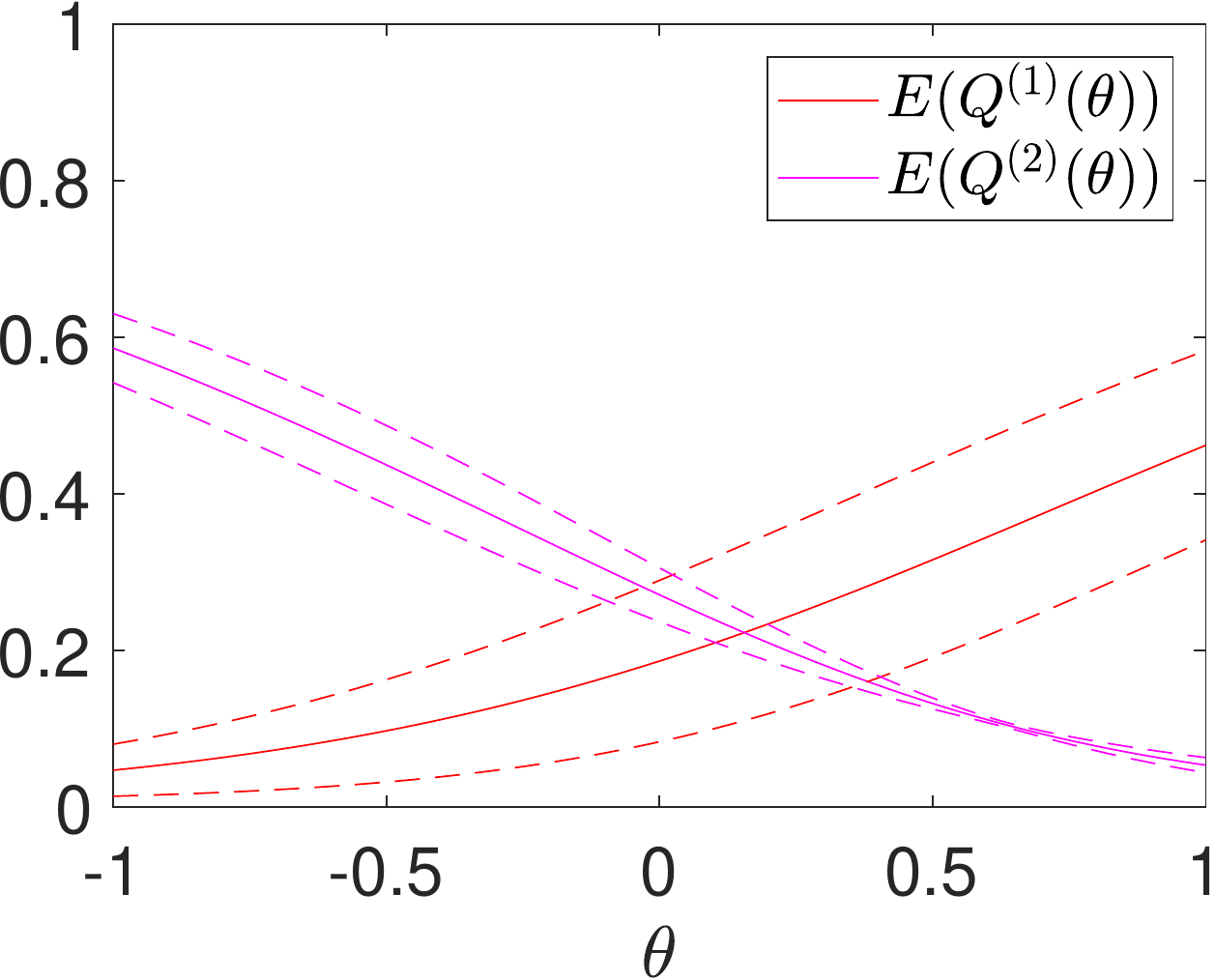}
\caption{Average position and standard deviation of the local unstable manifold with respect to the parameter $\theta$.}
\label{fig:tube_s_BetaCheb1}
\end{subfigure}\\
\vspace*{0.5cm}
\begin{subfigure}{0.49\linewidth}
\includegraphics[width=\linewidth]{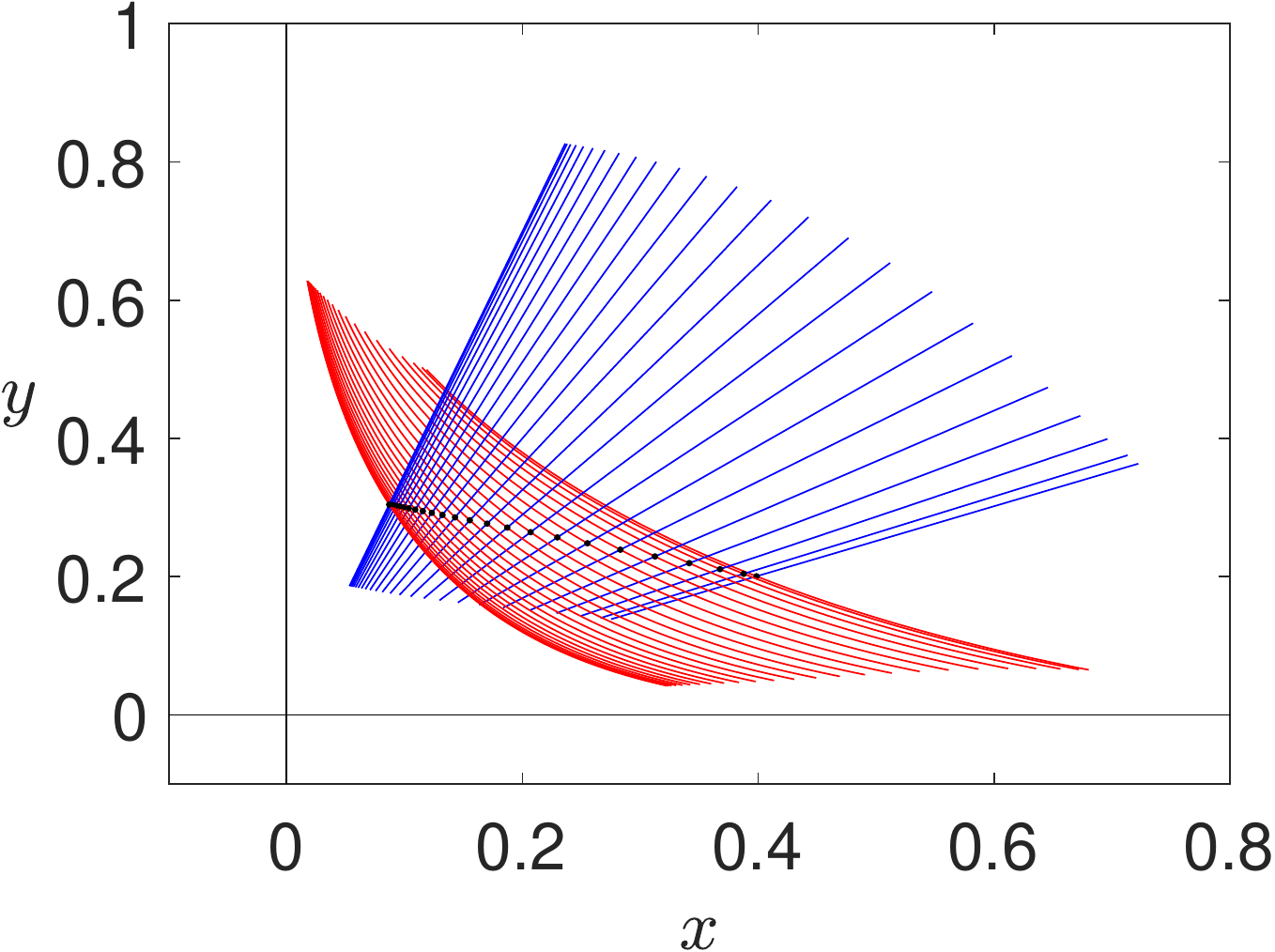}
\caption{Sampling of the manifolds in phase space.}
\label{fig:2D_BetaCheb1}
\end{subfigure}
\hfill
\begin{subfigure}{0.49\linewidth}
\includegraphics[width=\linewidth]{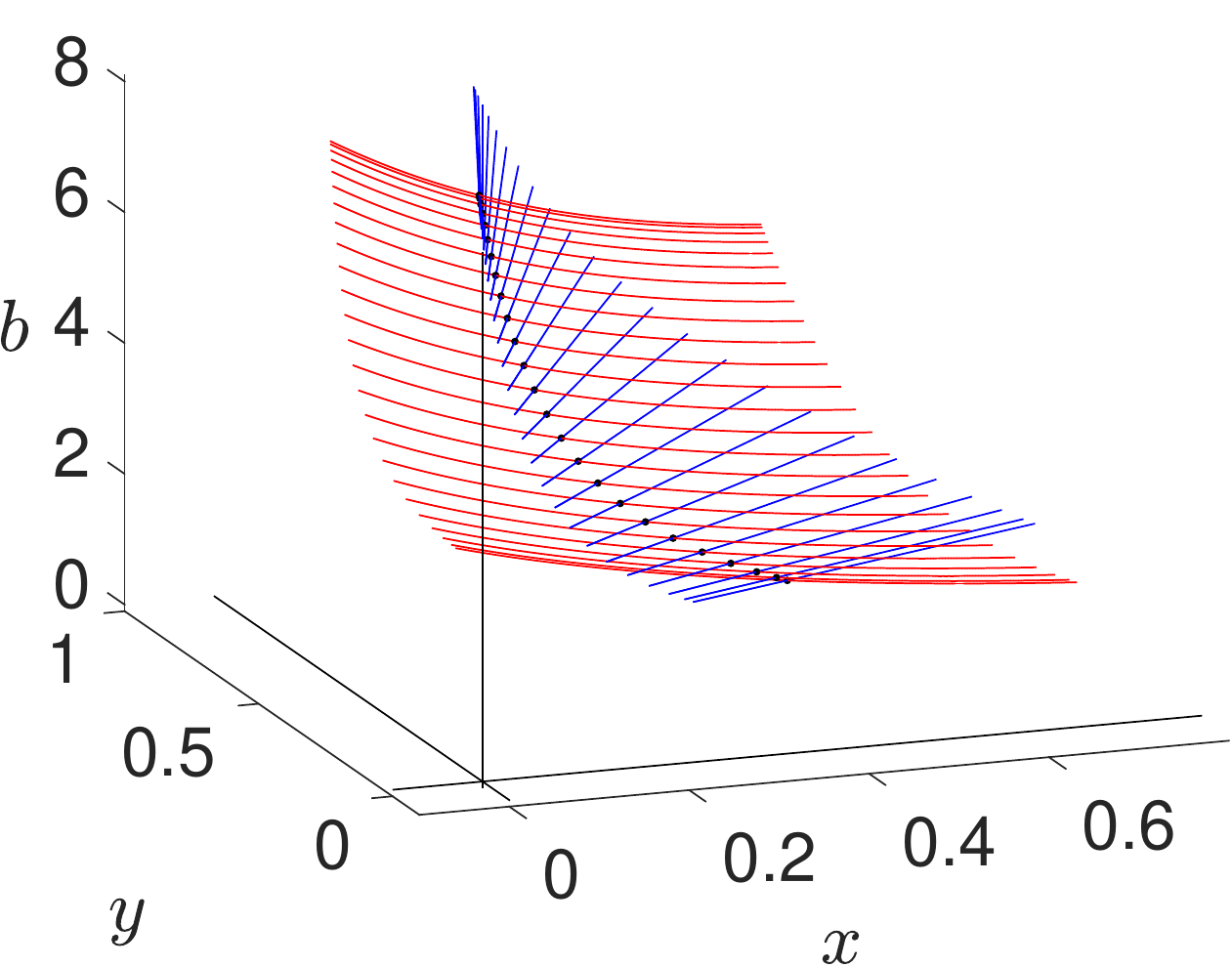}
\caption{Sampling of the manifolds in phase space $\times$ parameter space.}
\label{fig:3D_BetaCheb1}
\end{subfigure}
\caption{Representations of the local stable (in blue/cyan) and unstable (in red/magenta) manifolds of the equilibrium $(x_{\textnormal{eq}},y_{\textnormal{eq}})$. The computations were done for $a=3$, $\bb=5$, $\sigma=3$, $N=10$, $K=20$, $\gamma_1=0.2$, $\gamma_2=0.4$. The random parameter $p$ is assumed to have an arcsine distribution, i.e. $\rho_p(s)=\1_{s\in(-1,1)}\frac{1}{\pi\sqrt{1-s^2}}$ and therefore we used the Chebyshev polynomials of the first kind for the PC expansion. In Figure~\ref{fig:tube_u_BetaCheb1} (resp. Figure~\ref{fig:tube_s_BetaCheb1}) we display in full line the mean position of each component of the unstable (resp. stable) local manifold with respect to $\theta$: $\E\left(Q^{(i)}\right)(\theta)$, and in dotted line the associated standard deviation: $\E\left(Q^{(i)}\right)(\theta)\pm \sqrt{\V\left(Q^{(i)}\right)(\theta)}$. In Figure~\ref{fig:2D_BetaCheb1} we display manifolds for several values of $p$, sampled according to the arcsine distribution. The same manifolds are represented in Figure~\ref{fig:3D_BetaCheb1} with an additional dimension describing the value of $b=\bb+\sigma p$ corresponding to each sample.}
\label{fig:Manifold_LV_BetaCheb1}
\end{figure}

\begin{figure}[h]
\centering
\begin{subfigure}{0.49\linewidth}
\includegraphics[width=\linewidth]{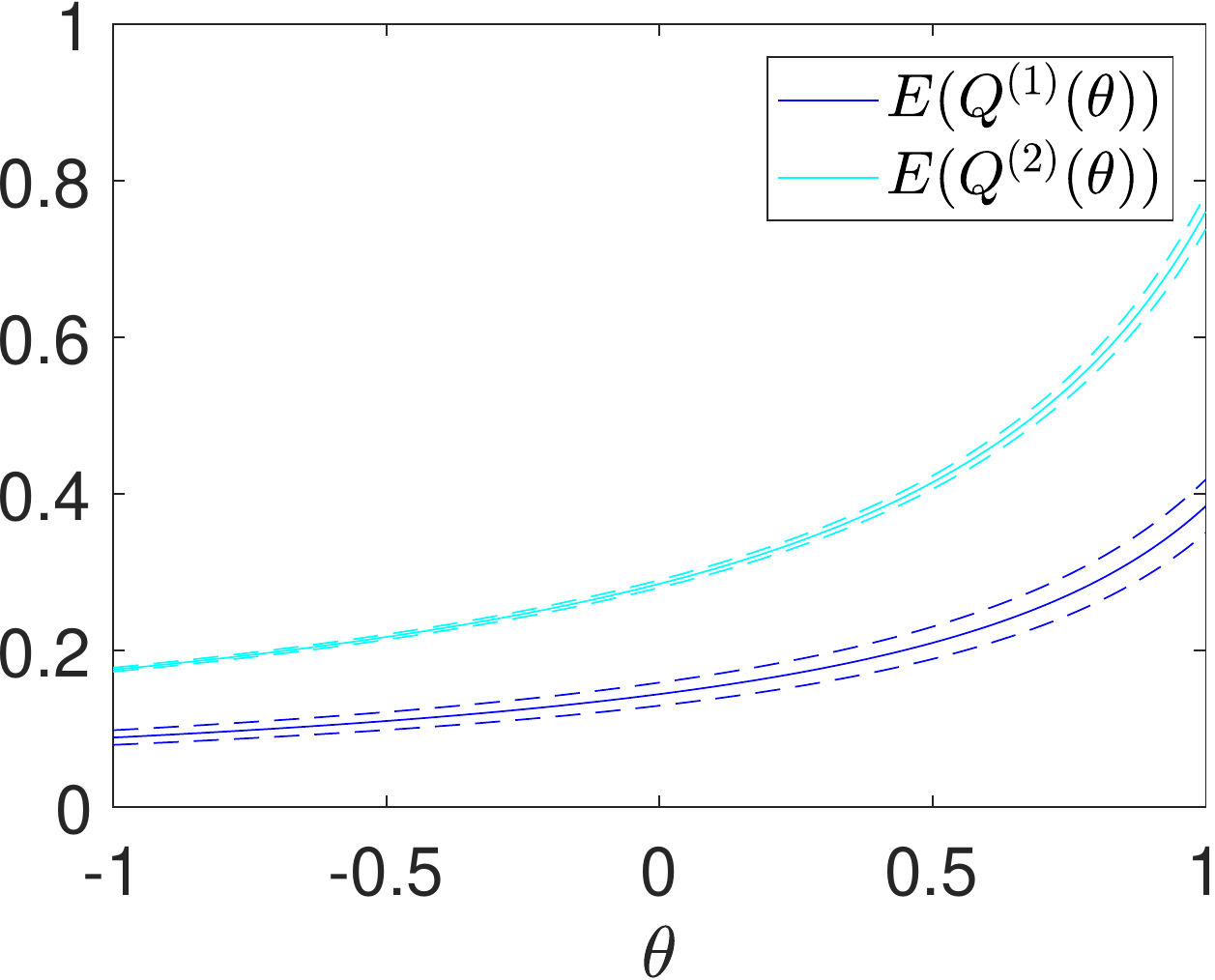}
\caption{Average position and standard deviation of the local stable manifold with respect to the parameter $\theta$.}
\label{fig:tube_u_BetaSymmu}
\end{subfigure}
\hfill
\begin{subfigure}{0.49\linewidth}
\includegraphics[width=\linewidth]{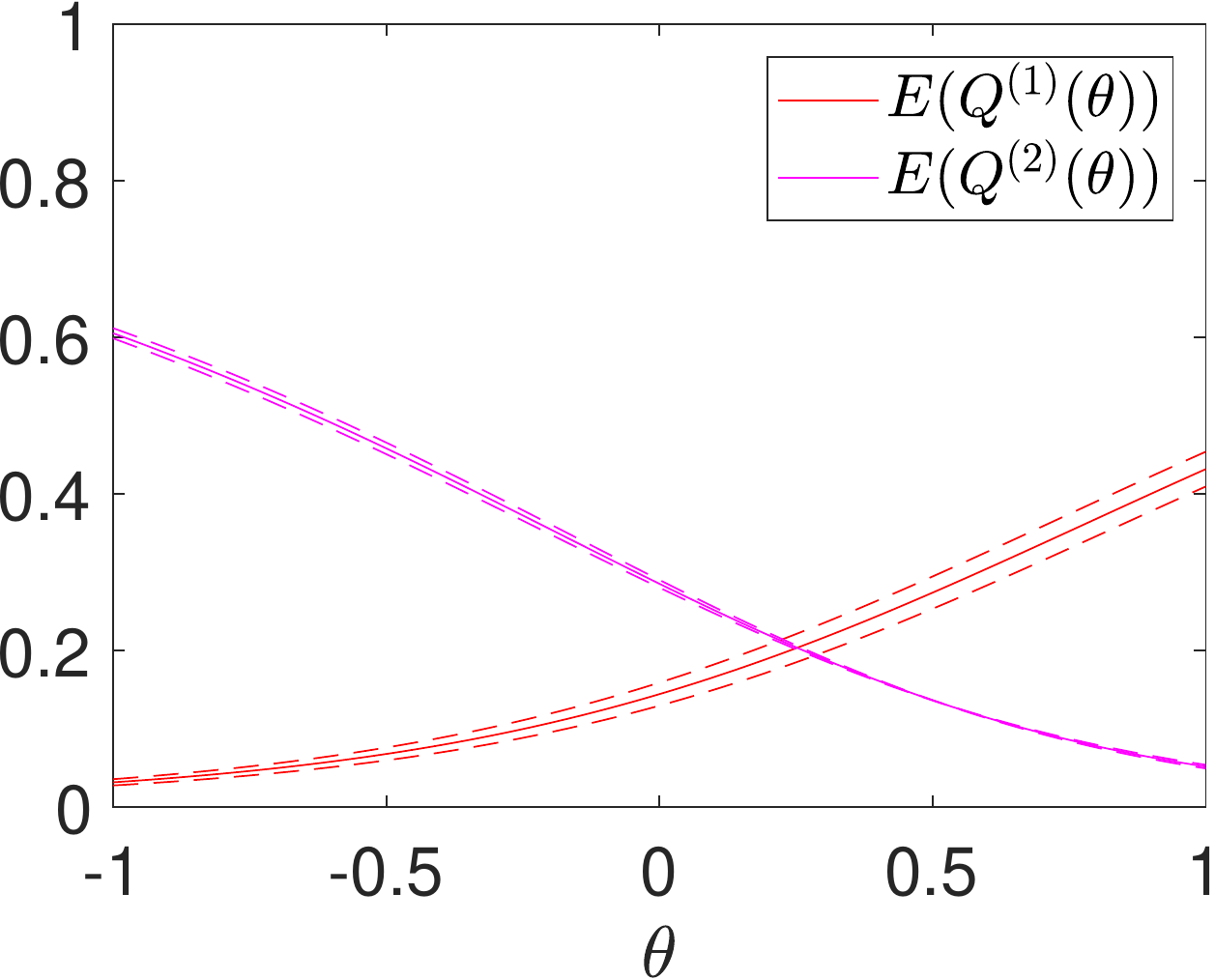}
\caption{Average position and standard deviation of the local unstable manifold with respect to the parameter $\theta$.}
\label{fig:tube_s_BetaSymmu}
\end{subfigure}
\begin{subfigure}{0.49\linewidth}
\includegraphics[width=\linewidth]{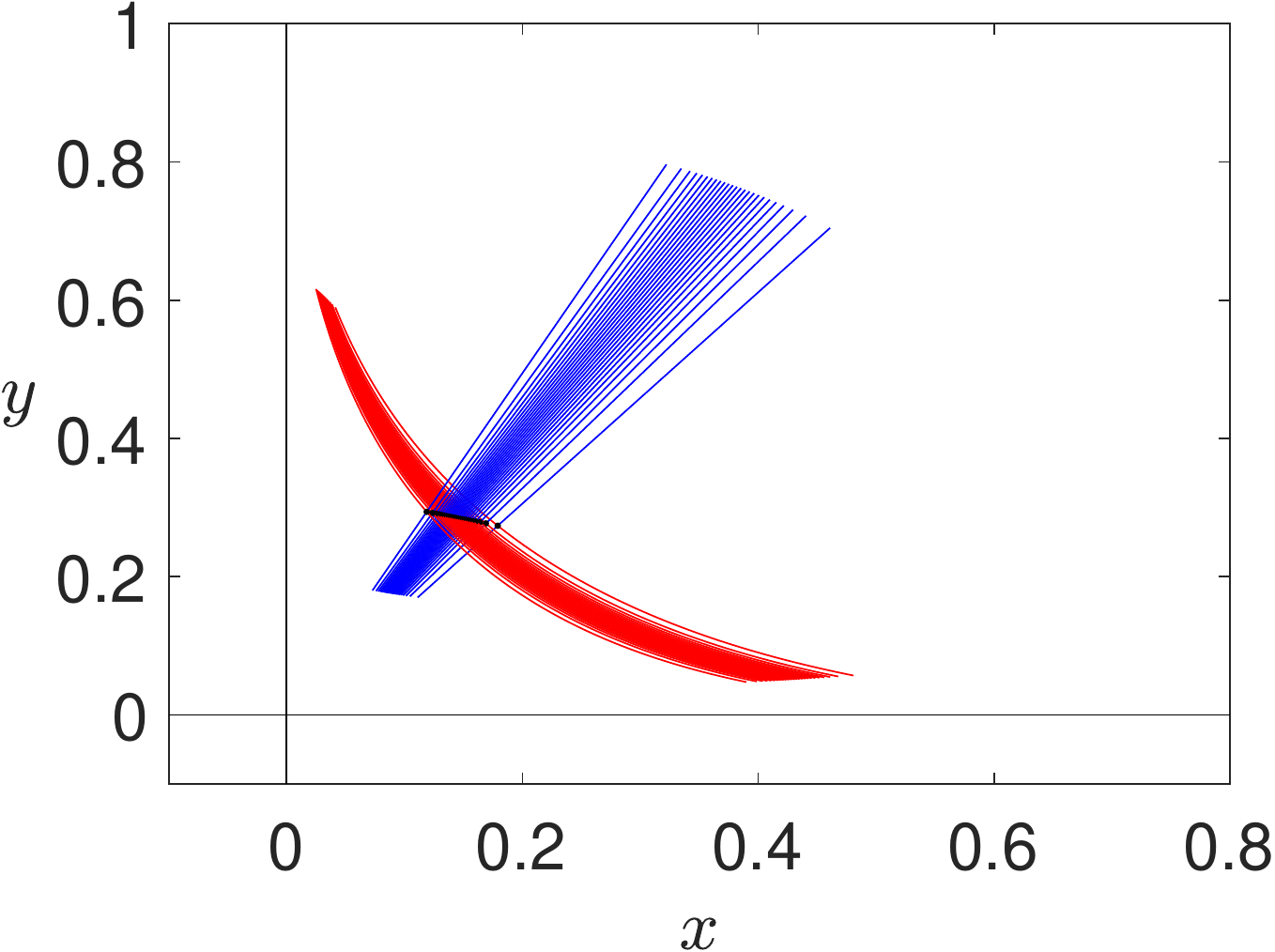}
\caption{Sampling of the manifolds in phase space.}
\label{fig:2D_BetaSymmu}
\end{subfigure}
\hfill
\begin{subfigure}{0.49\linewidth}
\includegraphics[width=\linewidth]{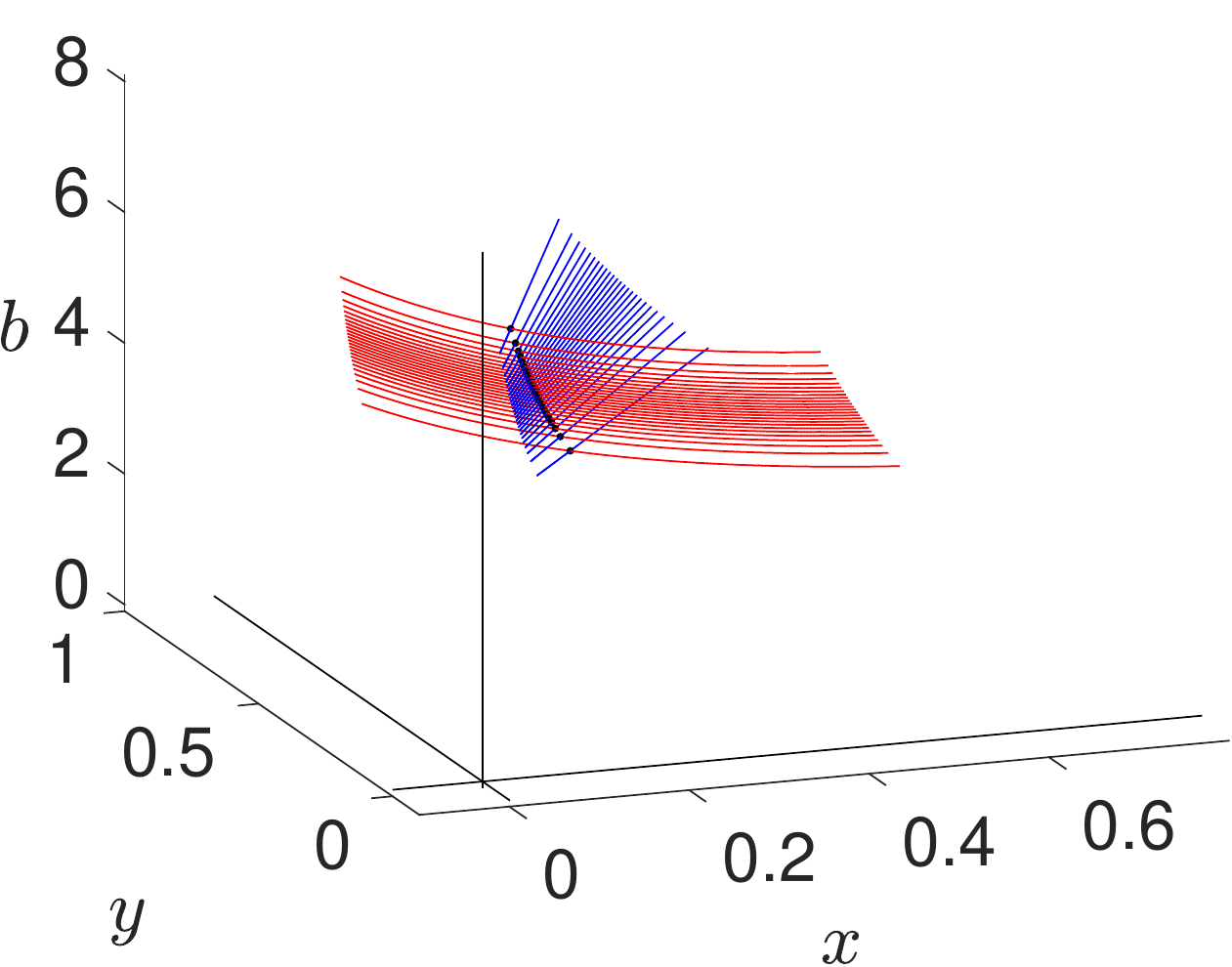}
\caption{Sampling of the manifolds in phase space $\times$ parameter space.}
\label{fig:3D_BetaSymmu}
\end{subfigure}
\caption{Representations of the local stable (in blue/cyan) and unstable (in red/magenta) manifolds of the equilibrium $(x_{\textnormal{eq}},y_{\textnormal{eq}})$. The computations were done for $a=3$, $\bb=5$, $\sigma=3$, $N=10$, $K=20$, $\gamma_1=0.2$, $\gamma_2=0.4$. The random parameter $p$ is assumed to have a beta distribution of parameter $(20,20)$, i.e. $\rho_p(s)=\1_{s\in(-1,1)}\frac{2^{2\mu-1}\mu B(\mu,\mu)}{\pi}(1-s^2)^{\mu-\frac{1}{2}}$ with $\mu=20$, and therefore we used the Gegenbauer polynomials of parameter $\mu=20$ for the PC expansion. In Figure~\ref{fig:tube_u_BetaSymmu} (resp. Figure~\ref{fig:tube_s_BetaSymmu}) we display in full line the mean position of each component of the unstable (resp. stable) local manifold with respect to $\theta$: $\E\left(Q^{(i)}\right)(\theta)$, and in dotted line the associated standard deviation: $\E\left(Q^{(i)}\right)(\theta)\pm \sqrt{\V\left(Q^{(i)}\right)(\theta)}$. In Figure~\ref{fig:2D_BetaSymmu} we display manifolds for several values of $p$, sampled according to the beta distribution of parameter $(20,20)$. The same manifolds are represented in Figure~\ref{fig:3D_BetaSymmu} with an additional dimension describing the value of $b=\bb+\sigma p$ corresponding to each sample.}
\label{fig:Manifold_LV_BetaSymmu}
\end{figure}

We are now ready to compute parameterizations of the stable and unstable manifolds of the equilibrium, as described in Section~\ref{sec:inv_manifold}. That is, we look for a parameterization of the form
\begin{equation*}
Q^{K,N}(\theta,p)=\sum_{k=0}^{K-1}\sum_{n=0}^{N-1}Q_{k,n} \phi_n(p)\theta^k.
\end{equation*}
Notice that for our example~\eqref{eq:LV_vector_field} the phase space is two-dimensional and both the stable and the unstable manifolds that we are interested in are one-dimensional. Therefore, $\theta$ will be one-dimensional and the coefficients $Q_{k,n}=\left(Q^{(1)}_{k,n},Q^{(2)}_{k,n}\right)$ belong to $\R^2$. It will be convenient to use the following notation
\begin{equation}
\label{eq:Q_k_hidden_n}
Q^{(1)}_k=\begin{pmatrix}
Q^{(1)}_{k,0} \\ \vdots \\ Q^{(1)}_{k,N-1} 
\end{pmatrix} \quad \text{and}\quad  Q^{(2)}_k=\begin{pmatrix}
Q^{(2)}_{k,0} \\ \vdots \\ Q^{(2)}_{k,N-1} 
\end{pmatrix}   \quad \forall~0\leq k<K.
\end{equation}
According to~\eqref{eq:CI_para} we define
\begin{equation*}
Q^{(1)}_0=\begin{pmatrix}
x_{0} \\ \vdots \\ x_{N-1}\end{pmatrix},\quad Q^{(2)}_0= \begin{pmatrix} y_{0} \\ \vdots \\ y_{N-1}
\end{pmatrix},\quad  Q_1^{(1)}=\gamma_i\begin{pmatrix}
u^{(i)}_{0} \\ \vdots \\ u^{(i)}_{N-1} \end{pmatrix}\quad \text{and}\quad Q_1^{(2)}=\gamma_i\begin{pmatrix}
v^{(i)}_{0} \\ \vdots \\ v^{(i)}_{N-1}
\end{pmatrix},
\end{equation*}
where $\left(x_n\right)_n$ and $\left(y_n\right)_n$ are the PC coefficients of the equilibrium computed in Section~\ref{sec:LV_eq}, and $(u^{(i)}_n)_n$ and $(v^{(i)}_n)_n$ are the PC coefficients of the eigenvectors computed in Section~\ref{sec:LV_eig}. Observe that $i=1$ if we want to compute a parameterization of the stable manifold (resp.~$i=2$ for the unstable manifold). To compute the higher order coefficients $Q_k$, $2\leq k<K$ we use~\eqref{eq:Q_k_rec_PC}, which we now specialize to our example~\eqref{eq:LV_vector_field}. We need to introduce some more notations. Given an expansion basis $\phi_n$ and $z=\left(z_n\right)_{0\leq n<N} \in \R^N$, we denote by $M_z$ the $N\times N$ \emph{ product matrix} such that
\begin{equation}
\label{eq:prod_mat}
(M_z z')_n = (z\ast z')_n \qquad \forall~z'\in\R^n,\ \forall~0\leq n<N.
\end{equation}
This matrix can be computed explicitly from $z$ and the linearization coefficients $\alpha$ (see the Appendix). Given an expansion basis $\phi_n$ and double expansions of the form
\begin{equation*}
q^{K,N}(\theta,p)= \sum_{k=0}^{K-1}\sum_{n=0}^{N-1} q_{k,n} \phi_n(p)\theta^k \qquad \text{and} \qquad q'^{K,N}(\theta,p)= \sum_{k=0}^{K-1}\sum_{n=0}^{N-1} q'_{k,n} \phi_n(p)\theta^k
\end{equation*}
we denote by $\circledast$ the product associated to the basis $\phi_n(p)\theta^k$, that is $\left(q\circledast q'\right)_{k,n}$ is the unique sequence such that
\begin{equation*}
q^{K,N}(\theta,p)q'^{K,N}(\theta,p) = \sum_{k=0}^{2K-2}\sum_{n=0}^{2N-2} \left(q\circledast q'\right)_{k,n} \phi_n(p)\theta^k.
\end{equation*} 
We again refer to the Appendix for more details about this product structure, and discussions on how it can be computed in practice. Similarly to~\eqref{eq:Q_k_hidden_n} we write
\begin{equation*}
 \left(q\circledast q'\right)_k=\begin{pmatrix}
\left(q\circledast q'\right)_{k,0} \\ \vdots \\ \left(q\circledast q'\right)_{k,N-1} 
\end{pmatrix} \qquad \forall~0\leq k<K.
\end{equation*} 
We are now ready to go back to~\eqref{eq:Q_k_rec_PC} and give explicit and commutable formulas for the coefficients $Q_k$ $k\geq 2$ of the parameterization in the case of our example~\eqref{eq:LV_vector_field}. Assuming the coefficients $Q^{(1)}_l$ and $Q^{(2)}_l$ have been computed for $0\leq l\leq k-1$, the next coefficients are given by
\renewcommand*{\arraystretch}{1.5}
\begin{align*}
\begin{pmatrix}
Q^{(1)}_k \\
Q^{(2)}_k
\end{pmatrix} &= 
\begin{pmatrix}
kM_{\lambda^{(i)}} -\left(I_N-2M_x-aM_y\right) & aM_x \\
-\left(-\bb M_y -\sigma M_{py}\right) & kM_{\lambda^{(i)}} - \left(I_N-2M_y-\bb M_x-\sigma M_{px}\right)
\end{pmatrix}^{-1} \\
& \qquad\qquad\qquad\qquad\qquad\qquad \begin{pmatrix}
\left(Q^{(1)}\circledast Q^{(1)}\right)_k - a \left(Q^{(1)}\circledast Q^{(2)}\right)_k \\
\left(Q^{(2)}\circledast Q^{(2)}\right)_k -\bb \left(Q^{(1)}\circledast Q^{(2)}\right)_k -\sigma \left(p\circledast Q^{(1)}\circledast Q^{(2)}\right)_k
\end{pmatrix},
\end{align*}
\renewcommand*{\arraystretch}{1}
where $I_N$ is the $N\times N$ identity matrix and again $i=1$ if we are computing a parameterization of the stable manifold (resp. $i=2$ for the unstable manifold). Here $p$ is identified with the sequence $p_{k,n}=\delta_{k,0}\delta_{n,1}$ (again assuming the $\phi_n$ are normalized such that $\phi_1(p)=p$). Once the coefficients $Q^{(1)}_k$ and $Q^{(2)}_k$ have been computed up to the desired order, we have access to the statistical properties of the manifold. Indeed we can for instance compute
\begin{equation*}
\E\left(Q^{(1)}\right)(\theta) = \sum_{k=0}^{K-1} Q^{(1)}_{k,0} \theta^k 
\end{equation*}
and
\begin{equation*}
\V\left(Q^{(1)}\right)(\theta) = \sum_{k=0}^{2K-2} \left(\sum_{n=1}^{N-1} \sum_{l=0}^k X_{l,n}X_{k-l,n} \right) \theta^k 
\end{equation*}
which describe the mean and variance of the first component of the manifold with respect to $\theta$ (see Figures~\ref{fig:Manifold_LV_BetaCheb1} and~\ref{fig:Manifold_LV_BetaSymmu}). We can also carry out cheap sampling to get a sense of how the manifolds are distributed in phase space (see again Figures~\ref{fig:Manifold_LV_BetaCheb1} and~\ref{fig:Manifold_LV_BetaSymmu}).

\section{Second example: the Lorenz system}
\label{sec:Lorenz}

In this section, we compute periodic orbits, invariant manifolds and heteroclinic orbits for the Lorenz system
\begin{equation}
\label{eq:Lorenz}
\left\{\begin{aligned}
\dot{x}&=\varsigma(y-x), \\
\dot{y}&=\varrho x-y-xz, \\
\dot{z}&=-\beta x +xy,
\end{aligned}\right.
\end{equation}
where $\varrho=\varrho(\omega)$ is a bounded random variable having a given distribution. We renormalize it by writing
\begin{equation*}
\varrho(\omega)=\brho + \sigma p(\omega),
\end{equation*}
where $\sigma>0$ and $p$ is a random variable taking values in $[-1,1]$.

\subsection{Computation of periodic orbits}

\begin{figure}[h]
\centering
\begin{subfigure}{0.49\linewidth}
\includegraphics[width=\linewidth]{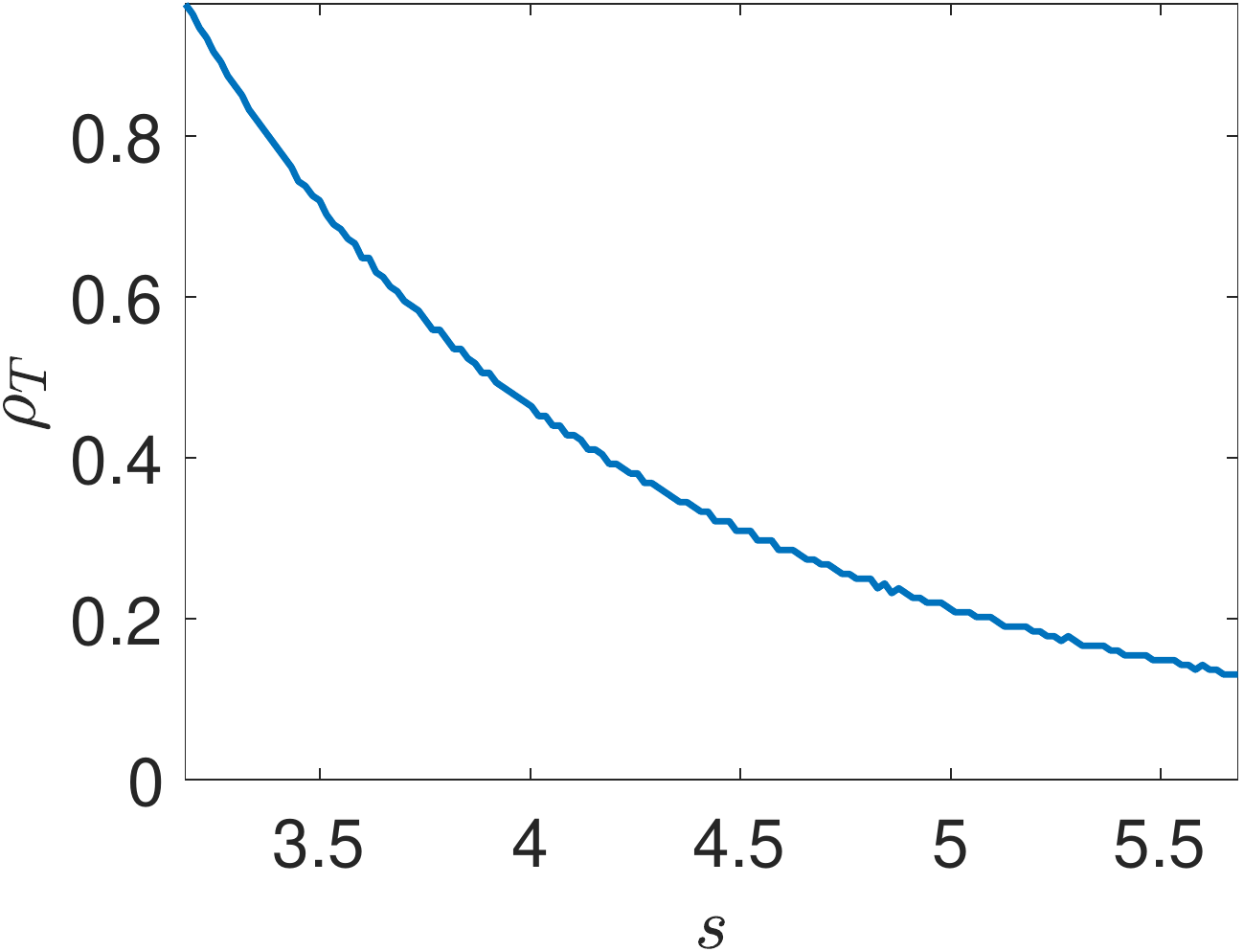}
\caption{PDF of the period when $\rho_p(s)=\frac{1}{2}\1_{s\in(-1,1)}$.}
\label{fig:density_period_Uniform}
\end{subfigure}
\hfill
\begin{subfigure}{0.49\linewidth}
\includegraphics[width=\linewidth]{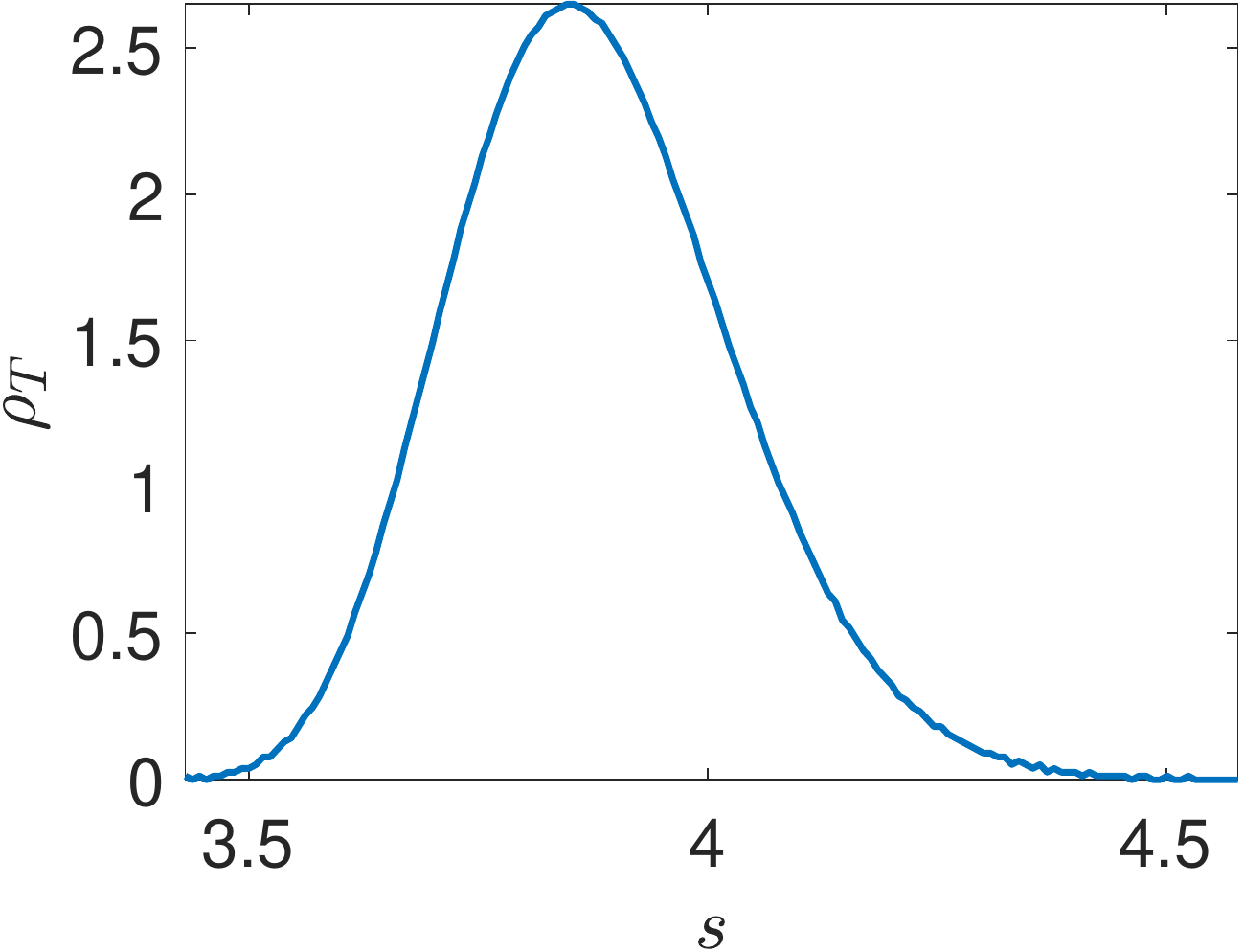}
\caption{PDF of the period when $\rho_p(s)=\1_{s\in(-1,1)}\frac{2^{2\mu-1}\mu B(\mu,\mu)}{\pi}(1-s^2)^{\mu-\frac{1}{2}}$ with $\mu=20$.}
\label{fig:density_period_BetaSymmu}
\end{subfigure}
\begin{subfigure}{0.70\linewidth}
\includegraphics[width=\linewidth]{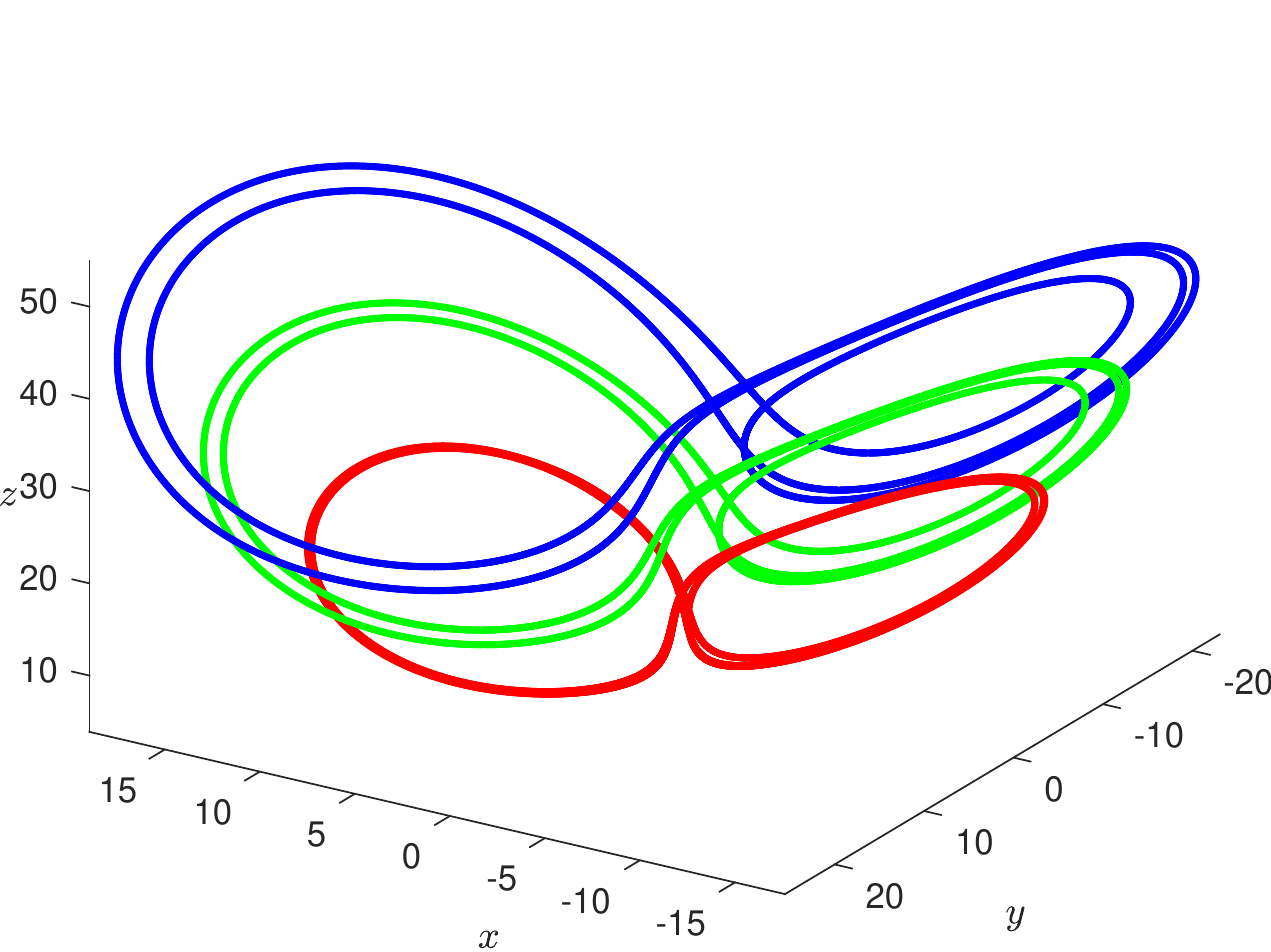}
\caption{Some of the periodic orbits contained in the expansion, in red for $\varrho=18$, in green for $\varrho=28$ and in blue for $\varrho=38$.}
\label{fig:several_orbits}
\end{subfigure}
\caption{We illustrate here some of the information that can be recovered from the Fourier$\times$PC expansion of the periodic orbit. The computations where done for $\varsigma=10$, $\beta=8/3$, $\brho=28$, $\sigma=10$, $K=80$, $N=10$ and several PDF $\rho_p$ of $p$. In Figure~\ref{fig:density_period_Uniform} and Figure~\ref{fig:density_period_BetaSymmu} we display the PDF of the period $T=2\pi/\Omega$ of the orbit. In Figure~\ref{fig:several_orbits}, we display some of the orbits described by the expansion for given values of $\varrho$ (notice that, from the expansion we can again sample cheaply to obtain such orbits for many values of $\varrho$, but the picture becomes cluttered very quickly).}
\label{fig:per_Lorenz}
\end{figure}

Using the framework presented in Section~\ref{sec:inv_periodic}, we compute periodic orbits of~\eqref{eq:Lorenz} via a Fourier$\times$PC expansion. That is, we write
\begin{equation*}
x^{K,N}(t,p)=\sum_{k=-K+1}^{K-1}\sum_{n=0}^{N-1}x_{k,n} \phi_n(p) \txte^{\txti k\Omega(p) t},
\end{equation*}
and similarly for $y^{K,N}$ and $z^{K,N}$, where
\begin{equation*}
\Omega(p)=\sum_{n=0}^{N-1} \Omega_n \phi_n(p).
\end{equation*}
In this subsection, we use $\circledast$ to denote the product associated to the basis $\phi_n(p) \txte^{\txti k\Omega(p) t}$, i.e. $(x\circledast y)_{k,n}$ is the unique sequence such that
\begin{equation}
\label{eq:exp_four}
x^{K,N}(t,p) y^{K,N}(t,p) = \sum_{k=-2K+2}^{2K-2}\sum_{n=0}^{2N-2}(x\circledast y)_{k,n}\phi_n(p) \txte^{\txti k\Omega(p) t}.
\end{equation}
It will again be convenient to write
\begin{equation}
\label{eq:not_x_k_four}
x_k=\begin{pmatrix}
x_{k,0} \\ \vdots \\ x_{k,N-1} 
\end{pmatrix}, \quad 
y_k=\begin{pmatrix}
y_{k,0} \\ \vdots \\ y_{k,N-1} 
\end{pmatrix}, \quad 
z_k=\begin{pmatrix}
z_{k,0} \\ \vdots \\ z_{k,N-1} 
\end{pmatrix}, \qquad \forall~0\leq k<K,
\end{equation}
and similarly
\begin{equation}
\label{eq:not_x_y_k_four}
 \left(x\circledast y\right)_k=\begin{pmatrix}
\left(x\circledast y\right)_{k,0} \\ \vdots \\ \left(x\circledast y\right)_{k,N-1} 
\end{pmatrix} \qquad \forall~0\leq k<K,
\end{equation} 
and so on. Plugging the expansions~\eqref{eq:exp_four} for $x$ , $y$ and $z$ in~\eqref{eq:Lorenz} and using the notations~\eqref{eq:not_x_k_four} and~\eqref{eq:not_x_y_k_four}, we obtain the following system of $(2K-1)\times N$ equations:
\begin{equation}
\label{eq:syst_coef_four}
\left\{\begin{aligned}
&\txti k(\Omega\circledast x)_k -\varsigma(y_k-x_k) = 0 \\
&\txti k(\Omega\circledast y)_k -\left(\brho x_k +\sigma (x\circledast p)_k -y_k -(x\circledast z)_k\right) = 0 \\
&\txti k(\Omega\circledast z)_k -\left(-\beta z_k +(x\circledast y)_k\right) = 0 \\
\end{aligned}\right. 
\qquad \forall~\vert k\vert <K,
\end{equation}
where again $p$ is identified with the sequence $p_{k,n}=\delta_{k,0}\delta_{n,1}$ and $\Omega$ with the sequence $\Omega_{k,n}=\delta_{k,0}\Omega_n$. As mentioned in Section~\ref{sec:inv_periodic}, we complement this system with a phase condition, that we also expand using PC. That is, given
\begin{equation*}
u^{(i)}=\begin{pmatrix}
u^{(i)}_0 \\ \vdots \\ u^{(i)}_{N-1}
\end{pmatrix} \qquad \text{and} \qquad v^{(i)}=\begin{pmatrix}
v^{(i)}_0 \\ \vdots \\ v^{(i)}_{N-1}
\end{pmatrix} \qquad i=1,2,3,
\end{equation*}
we append to~\eqref{eq:syst_coef_four} the $N$ equations
\begin{equation}
\label{eq:phase_cond_PC}
\left(\sum_{k=-K+1}^{K-1} x_k - u^{(1)}\right)\ast v^{(1)} + \left(\sum_{k=-K+1}^{K-1} y_k - u^{(2)}\right)\ast v^{(2)} + 
\left(\sum_{k=-K+1}^{K-1} z_k - u^{(3)}\right)\ast v^{(3)} =0.
\end{equation}
In practice, we solve~\eqref{eq:syst_coef_four}-\eqref{eq:phase_cond_PC} using Newton's method. More precisely, we assume that we are given a deterministic periodic orbit (i.e. a solution for $\sigma=0$) and do a predictor-corrector continuation in $\sigma$ (using Newton's method as the corrector step) until we reach the desired value. In Figure~\ref{fig:per_Lorenz} we illustrate the output of this procedure. 

\subsection{Computation of heteroclinic orbits}

In this subsection, we detail how our approach can be used to compute heteroclinic orbits for the Lorenz system, going between \begin{equation}
\label{eq:eyep}
\left(\sqrt{\beta(\varrho-1)},\sqrt{\beta(\varrho-1)},\varrho-1\right)
\end{equation}
and the origin. We focus on the case where
\begin{equation*}
\varrho >\frac{\varsigma\left(\varsigma+\beta+3\right)}{\varsigma+\beta-1},
\end{equation*}
in which the origin has a two-dimensional stable manifold and~\eqref{eq:eyep} has a two-dimensional unstable manifold. 

\begin{figure}[h]
\centering
\includegraphics[width=\linewidth]{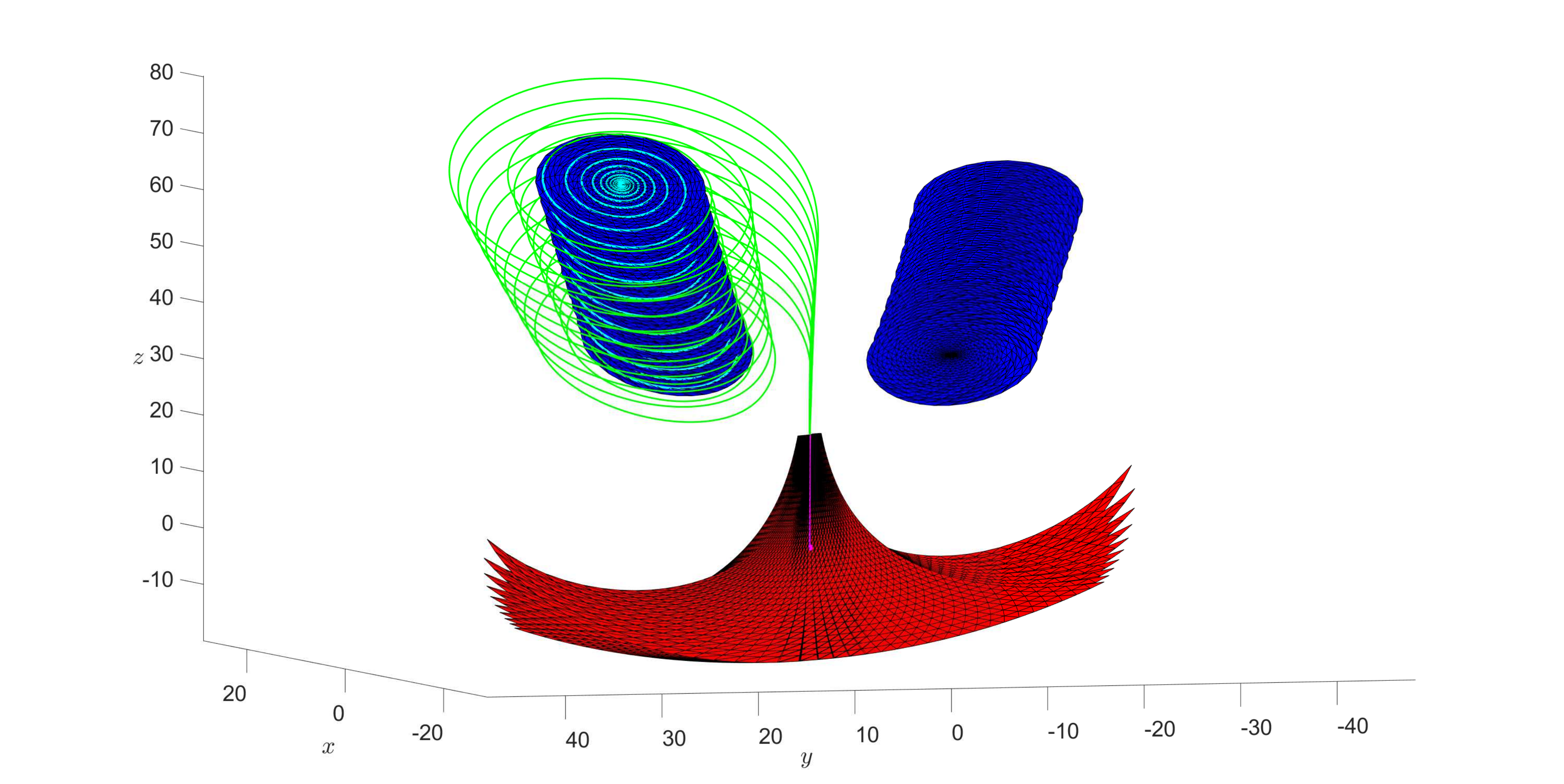}
\caption{Some of the heteroclinic orbits encoded in the piece-wise Chebyshev$\times$PC expansions and in the parameterizations of the local manifolds. The computations where done for $\varsigma=10$, $\beta=8/3$, $\brho=50$, $\sigma=15$, $M=25$, $K=8$, $N=15$ and $\hat K=30$, $\check K=30$. We display the manifolds and orbits for 10 values of $p$, uniformly sampled. The local unstable manifolds of~\eqref{eq:eyep} are in blue and the local stable manifolds of the origin are in red. The part of the heteroclinic orbits that solves the boundary value problem between the two manifolds and is computed via piece-wise Chebyshev$\times$PC expansions is displayed in green. The remaining parts of the heteroclinic orbits (in cyan and magenta), are obtained \emph{for free} via the conjugation properties of the parameterizations (as explained on Figure~\ref{fig:PM}).}
\label{fig:connexion_Lorenz_rho_50}
\end{figure}

As explained in Section~\ref{sec:inv_heteroclinic}, we use these local manifolds to set up our boundary value problem for the heteroclinic orbits. A parameterization of these local manifolds is computed as described in Section~\ref{sec:inv_manifold}, using a Taylor$\times$PC expansion. This procedure was already described in details in Section~\ref{sec:LV_manifolds}, so we omit these details here. We denote by $\hat Q^{\hat K,N}$ (resp.~$\check Q^{\check K,N}$) a local parameterization of the unstable manifold of~\eqref{eq:eyep} (resp. of the stable manifold of the origin) of the form
\begin{equation*}
\hat Q^{\hat K,N} (\theta,p) = \sum_{0\leq \vert k\vert < \hat K} \sum_{0\leq n<N} \hat Q_{k,n} \phi_n(p) \theta^k.
\end{equation*}
Notice that since both manifolds are two-dimensional, we have $\theta=(\theta_1,\theta_2)$. We recall that our goal is to find $\tau(p)$, $\hat\theta(p)$, $\check\theta(p)$ and an orbit $X(t,p)$ that solves the Lorenz system on $[0,\tau(p)]$, and satisfies the boundary conditions $X(0,p)=\hat Q^{\hat K,N} (\hat\theta(p),p)$ and $X(\tau(p),p)=\check Q^{\hat K,N} (\check\theta(p),p)$. As mentioned in Remark~\ref{rem:para_BVP}, this system is underdetermined and we can in fact fix $\hat\theta_1$ and $\check\theta_1$, and only solve for $\tau(p)$, $\hat\theta_2(p)$, $\check\theta_2(p)$ and $X(t,p)$ to recover a unique solution. In practice, we solve for PC expansions of $\tau(p)$, $\hat\theta_2(p)$ and $\check\theta_2(p)$:
\begin{equation*}
\tau^N(p)=\sum_{n=0}^{N-1} \tau_n \phi_n(p),\quad \hat\theta^N_2(p)=\sum_{n=0}^{N-1} \left(\hat\theta_2\right)_n \phi_n(p)\quad \text{and}\quad \check\theta^N_2(p)=\sum_{n=0}^{N-1} \left(\check\theta_2\right)_n \phi_n(p).
\end{equation*}
To compute the heteroclinic orbit (or to be more precise, the part of that orbit that connects the two local manifolds), we use piece-wise Chebyshev$\times$PC expansions, as exposed in Section~\ref{sec:inv_heteroclinic}. That is, we now write
\begin{equation*}
x^{K,N}(t,p)=\sum_{k=-K+1}^{K-1}\sum_{n=0}^{N-1}x^{(j)}_{k,n} \phi_n(p) T_{\vert k\vert}^{(j)}(t,p), \quad \forall~t\in(t^{(j-1)}(p),t^{(j)}(p)),\ \forall~j\in\{1,\ldots,J\},
\end{equation*}
and similarly for $y^{K,N}$ and $z^{K,N}$, where we use the notations of Section~\ref{sec:inv_heteroclinic} for partition of $[0,\tau]$ and the associated rescaled Chebyshev polynomials. In this subsection, we use $\circledast$ to denote the product associated to the basis $\phi_n(p) T_{\vert k\vert}^{(j)}(t,p)$, i.e. $(x\circledast y)_{k,n}$ is the unique sequence such that
\begin{equation}
\label{eq:exp_cheb}
x^{K,N}(t,p) y^{K,N}(t,p) = \sum_{k=-2K+2}^{2K-2}\sum_{n=0}^{2N-2}(x\circledast y)^{(j)}_{k,n}\phi_n(p) T_{\vert k\vert}^{(j)}(t,p), \quad \forall~t\in(t^{(j-1)}(p),t^{(j)}(p)),\ \forall~j\in\{1,\ldots,J\}.
\end{equation}
It will again be convenient to write
\begin{equation*}
\label{eq:not_x_k_cheb}
x^{(j)}_k=\begin{pmatrix}
x^{(j)}_{k,0} \\ \vdots \\ x^{(j)}_{k,N-1} 
\end{pmatrix}, \quad 
y^{(j)}_k=\begin{pmatrix}
y^{(j)}_{k,0} \\ \vdots \\ y^{(j)}_{k,N-1} 
\end{pmatrix}, \quad 
z^{(j)}_k=\begin{pmatrix}
z^{(j)}_{k,0} \\ \vdots \\ z^{(j)}_{k,N-1} 
\end{pmatrix}, \qquad \forall~0\leq k<K,
\end{equation*}
and similarly
\begin{equation*}
\label{eq:not_x_y_k_cheb}
 \hat Q_k=\begin{pmatrix}
 \hat Q_{k,0} \\ \vdots \\  \hat Q_{k,N-1} 
\end{pmatrix} \qquad \forall~0\leq k<K,
\end{equation*} 
\begin{equation*}
\label{eq:not_Q_k}
 \left(x\circledast y\right)^{(j)}_k=\begin{pmatrix}
\left(x\circledast y\right)^{(j)}_{k,0} \\ \vdots \\ \left(x\circledast y\right)^{(j)}_{k,N-1} 
\end{pmatrix} \qquad \forall~0\leq k<K,~\forall~j\in\{1,\ldots,J\},
\end{equation*} 
and so on. To write all three components at once, we also use
\begin{equation*}
X^{(j)}_k = \begin{pmatrix}
x^{(j)}_k \\ y^{(j)}_k \\ z^{(j)}_k
\end{pmatrix}.
\end{equation*}
With these notations, system~\eqref{eq:BVP_coeffs} for the Lorenz vector field is then given by
\begin{equation*}
\label{eq:BVP_coeffs_Lorenz}
\left\{\begin{aligned}
&kx^{(j)}_k-\frac{\tilde t^{(j)}-\tilde t^{(j-1)}}{4}\left(\tau\ast\left(\varsigma (y^{(j)}_{k-1}-x^{(j)}_{k-1})-\varsigma (y^{(j)}_{k+1}-x^{(j)}_{k+1})\right)\right)=0,\quad \forall~k\geq 1,\ \forall~1\leq j\leq J \\
&ky^{(j)}_k-\frac{\tilde t^{(j)}-\tilde t^{(j-1)}}{4}\left(\tau\ast\left(\left(\brho x^{(j)}_{k-1} +\sigma (x\circledast p)^{(j)}_{k-1} -y^{(j)}_{k-1} -(x\circledast z)^{(j)}_{k-1}\right) \right.\right. \\\
&\qquad\qquad\qquad\qquad\qquad\qquad \left.\left.- \left(\brho x^{(j)}_{k+1} +\sigma (x\circledast p)^{(j)}_{k+1} -y^{(j)}_{k+1} -(x\circledast z)^{(j)}_{k+1}\right)\right)\right)=0,\quad \forall~k\geq 1,\ \forall~1\leq j\leq J \\
&zy^{(j)}_k-\frac{\tilde t^{(j)}-\tilde t^{(j-1)}}{4}\left(\tau\ast\left(\left(-\beta z^{(j)}_{k-1} +(x\circledast y)^{(j)}_{k-1}\right) - \left(-\beta z^{(j)}_{k+1} +(x\circledast y)^{(j)}_{k+1}\right)\right)\right)=0,\quad \forall~k\geq 1,\ \forall~1\leq j\leq J \\
&X^{(j)}_0+2\sum_{k=1}^{K-1} X^{(j)}_k = X^{(j+1)}_{0}+2\sum_{k=1}^{K-1} (-1)^kX^{(j+1)}_{k},\quad \forall~1\leq j\leq M-1 \\
&X^{(1)}_0+2\sum_{k=1}^{K-1} (-1)^kX^{(1)}_k = \sum_{0\leq \vert k\vert < \hat K} \hat Q_{k} \ast \hat\theta^k, \\
&X^{(M)}_0(p)+2\sum_{k=1}^{K-1} X^{(J)}_k = \sum_{0\leq \vert k\vert < \check K} \check Q_{k}\ast \check\theta^k,
\end{aligned}\right.
\end{equation*}
where, for $k=(k_1,k_2)$, $\hat\theta^k$ must be understood as 
\begin{equation*}
\underbrace{\hat\theta_1 \ast \ldots \ast \hat\theta_1}_{k_1 \text{ times}} \ast \underbrace{\hat\theta_2 \ast \ldots \ast \hat\theta_2}_{k_2 \text{ times}}.
\end{equation*}
Again, we solve this system for the desired noise level $\sigma$ by doing predictor corrector steps with Newton iterations. In Figure~\ref{fig:connexion_Lorenz_rho_50} we illustrate the output of this procedure. 

\begin{figure}[h]
\centering
\includegraphics[width=\linewidth]{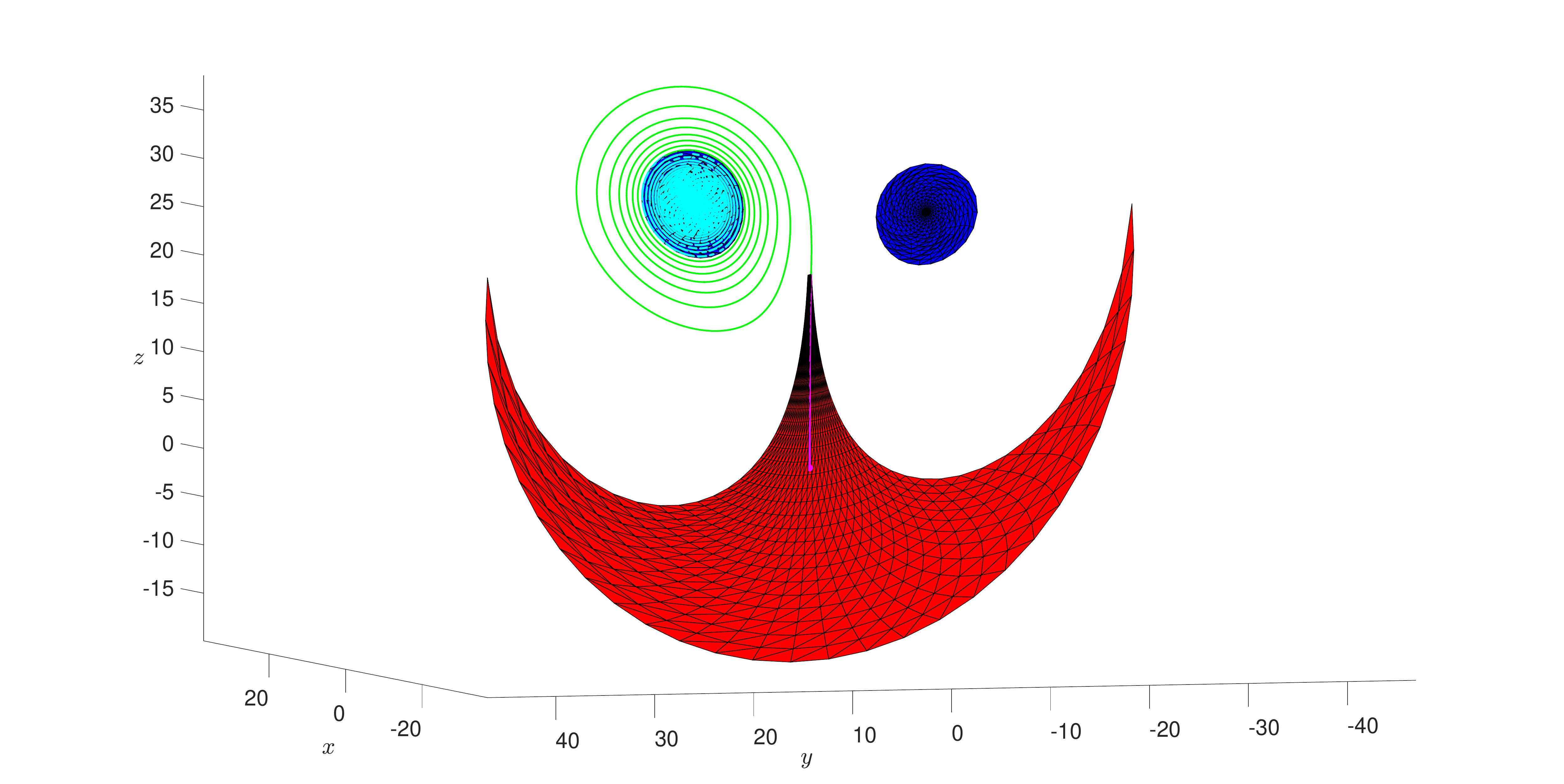}
\caption{An example of heteroclinic orbit encoded in the piece-wise Chebyshev$\times$PC expansions and in the parameterizations of the local manifolds. The computation was done for $\varsigma=10$, $\beta=8/3$, $\brho=28$, $\sigma=1$, $M=25$, $K=8$, $N=15$ and $\hat K=30$, $\check K=90$. Because of the relatively small value of $\sigma$ we only display one orbit. The local unstable manifold of~\eqref{eq:eyep} is in blue and the local stable manifold of the origin is in red. The part of the heteroclinic orbit that solves the boundary value problem between the two manifolds and is computed via piece-wise Chebyshev$\times$PC expansions is displayed in green. The remaining parts of the heteroclinic orbit (in cyan and magenta), are obtained \emph{for free} via the conjugation properties of the parameterizations (as explained on Figure~\ref{fig:PM})}.
\label{fig:connexion_Lorenz_rho_28}
\end{figure}

We point out that, while our method is also successful when $\varrho$ takes values around the classical value $\varrho=28$, the computation is then more challenging due to the proximity of the Hopf bifurcation at $\varrho =\frac{\varsigma\left(\varsigma+\beta+3\right)}{\varsigma+\beta-1}$ ($\approx 24.7$ for the classical values $\varsigma=10$ and $\beta=8/3$), and therefore we are only able to handle smaller noise level (i.e.~$\sigma\approx 1$, see Figure~\ref{fig:connexion_Lorenz_rho_28}). For such parameter values,  the computations (in both the stochastic and the deterministic framework) are intrinsically more difficult because of the strong oscillatory behavior induced by the pair of eigenvalues having a close to zero real part.

\section{Conclusions \& Outlook}

PC expansions have proven very successful in studying evolution equations (be it ODEs or PDEs) with random coefficients, but mostly from the point of view of time integration. In this work we developed a complementary approach, also based on PC but aimed at studying invariant sets of such systems, and applied it to investigate steady states, stable and unstable manifolds, periodic orbits and connecting orbits for ODEs. This approach is driven by the paradigm of nonlinear dynamics to study the structure of invariant sets to understand the relevant effects. Once a PC expansion representation of an invariant set is computed, this expansion contains explicit information about the random invariant sets. Using fast sampling or geometric visualization of the moments, then allows us to understand, how likely different phase space structures are going to be based upon the random input. 

We conclude by providing a brief discussion about some potential generalizations and further directions of research connected to this work.

\paragraph{Extension to PDEs:} A natural extension of this work would be to use similar techniques to study invariant sets of time-dependent PDEs with random coefficients. The steady state case has already been extensively investigated, but we believe that our new framework could allow to complement the already existing studies on periodic orbits, and to explore new problems related to invariant manifolds and connecting orbits for PDEs with random coefficients.

\paragraph{Bifurcations:} If the invariant state we are studying undergoes a bifurcation for some value of the random parameter, then the curve or manifold of invariant state that we are looking for may not be smooth, and the PC expansions will then converge slowly, if at all. Several \emph{multiresolution analysis schemes}, based on subdivising the random space or on different bases such as wavelets, were developed to handle such situations (see e.g.~\cite{WanKar06},\cite[Chapter 8]{LeMKni10} and the references therein), and could be used also in our setting to study bifurcation problems.

\paragraph{Rigorous computation:} A posteriori error analysis for PC expansions is of course critical, as the quantification of uncertainty provided by those expansions is only relevant if the error coming from the discretization/truncation can also be controlled. Techniques of \emph{rigorous computations} have been developed to obtain certified a-posteriori error estimates about numerically computed invariant sets of deterministic systems~\cite{FiGamLesLla17}, and we aim at generalizing them for random systems in a future work.

\section{Appendix}

We discuss here the computations of products of PC expansions, and detail some implementation aspects.

\subsection{Linearization formulas for products}

For any family of orthogonal polynomials $\phi_n$ associated to a weight $\rho$, and any $m$ and $n$ in $\N$, the product $\phi_m\phi_n$ can be written in the original basis:
\begin{equation}
\label{eq:lin1}
\phi_m\phi_n = \sum_{k=0}^{n+m} \alpha^{m,n}_k \phi_k.
\end{equation}
This is often called a \emph{linearization} formula. By orthogonality, the linearization coefficients satisfy
\begin{equation*}
\alpha^{m,n}_k = \frac{1}{h_k}\left\langle \phi_m\phi_n,\phi_k \right\rangle_\rho = \frac{1}{h_k}\left\langle \phi_k\phi_n,\phi_m \right\rangle_\rho,
\end{equation*}
from which we infer that $\alpha^{m,n}_k = \alpha^{n,m}_k$ and that, if $k+n<m$ or $k+m<n$, then $\alpha^{m,n}_k=0$, i.e.
\begin{equation}
\label{eq:lin1bis}
\phi_m\phi_n = \sum_{k=\vert m-n\vert}^{n+m} \alpha^{m,n}_k \phi_k.
\end{equation}
Besides, if the weight $\rho$ is even, then any even function is orthogonal to any odd function and therefore, for all $k$ having a different parity than $m+n$, $\alpha^{m,n}_k=0$. In such case, it can be convenient to eliminate all the coefficients that are a priori equal to zero, and rewrite~\eqref{eq:lin1} as
\begin{equation}
\label{eq:lin2}
\phi_m\phi_n = \sum_{k=0}^{\min(m,n)} \tilde\alpha^{m,n}_k \phi_{m+n-2k}.
\end{equation}

\begin{remark} 
In practice, it is convenient to precompute and store the linearization coefficients $\alpha^{m,n}_k$ (or $\tilde\alpha^{m,n}_k$). These coefficients can be obtained by numerically computing the integrals $\left\langle \phi_m\phi_n,\phi_k \right\rangle_\rho$, for instance using quadrature rules. For the classical orthogonal polynomials, the linearization coefficients can also be computed in closed form (see e.g.~\cite{OlvLozBoiCla10}). For Legendre polynomials $P_n$, we have
\begin{equation*}
\tilde\alpha^{m,n}_k = \frac{\binom{m-k-1/2}{m-k}\binom{n-k-1/2}{n-k}\binom{k-1/2}{k}(n+m-2k+1/2)}{\binom{m+n-k-1/2}{m+n-k}(n+m-k+1/2)}, \qquad \forall~m,n,k\in\N,\ k\leq \min(m,n),
\end{equation*}
where
\begin{equation*}
\binom{z}{k}=\frac{\Gamma(z+1)}{\Gamma(k+1)\Gamma(z-k+1)}.
\end{equation*}
For Chebyshev polynomials of the first kind $T_n$, we have
\begin{equation*}
\tilde\alpha^{m,n}_k = \left\{\begin{aligned}
&1 \qquad &\text{if }k=\min(m,n)=0, \\
&1/2 \qquad &\text{if }\min(m,n)>0,\text{ and } k=0 \text{ or } k=\min(m,n), \\
&0 \qquad &\text{otherwise}.
\end{aligned}\right.
\end{equation*}
For Chebyshev polynomials of the second kind $U_n$, we have
\begin{equation*}
\tilde\alpha^{m,n}_k = 1, \qquad \forall~m,n,k\in\N,\ k\leq \min(m,n).
\end{equation*}
For Gegenbauer polynomials $C^\mu_n$, we have
\begin{equation*}
\tilde\alpha^{m,n}_k = \frac{\binom{m-k+\mu-1}{m-k}\binom{n-k+\mu-1}{n-k}\binom{k+\mu-1}{k}(n+m-2k+\mu)}{\binom{m+n-k+\mu-1}{m+n-k}(n+m-k+\mu)}, \qquad \forall~m,n,k\in\N,\ k\leq \min(m,n),
\end{equation*}
Finally, we point out that linearization coefficients $\alpha^{m,n}_k$ can of course also be defined for arbitrary (i.e. non necessarily orthogonal) bases of polynomials. In particular, for the canonical basis $\phi_n(s)=s^n$ associated to Taylor expansions we have
\begin{equation*}
\tilde\alpha^{m,n}_k = \left\{\begin{aligned}
&1 \qquad &\text{if }k=0, \\
&0 \qquad &\text{otherwise}.
\end{aligned}\right.
\end{equation*}
\end{remark}

\begin{definition}
\label{def:convo}
Given to sequences $u=\left(u_n\right)_{n\in\N}$ and $v=\left(v_n\right)_{n\in\N}$, we define their \emph{convolution product} $u\ast v$ by
\begin{equation*}
\left(u\ast v\right)_k = \sum_{m=0}^\infty\sum_{n=0}^\infty u_m v_n \alpha^{m,n}_k, \qquad \forall~k\in\N,
\end{equation*}
with the convention $\alpha^{m,n}_k=0$ for all $k>m+n$. Notice that this definition depends on the weight $\rho$, or equivalently on the family $\left(\phi_n\right)_{n\in\N}$, via the coefficients $\alpha^{m,n}_k$.
\end{definition}

\begin{lemma}
If the weight $\rho$ is even, the convolution of $u$ and $v$ can also be written as
\begin{equation*}
\left(u\ast v\right)_k = \sum_{p=0}^\infty\sum_{q=0}^k u_{p+q} v_{p+k-q} \tilde\alpha^{p+q,p+k-q}_p.
\end{equation*}
\end{lemma}

This definition is the natural one to describe the product of two functions in the basis given by $\left(\phi_n\right)_{n\in\N}$. Indeed, writing
\begin{equation*}
u(s)=\sum_{n=0}^\infty u_n \phi_n(s) \qquad \text{and}\qquad v(s)=\sum_{n=0}^\infty v_n \phi_n(s),
\end{equation*}
one has, at least formally,
\begin{equation}
\label{eq:prod_star}
u(s)v(s)=\sum_{n=0}^\infty \left(u\ast v\right)_n \phi_n(s).
\end{equation}

With the notations of Definition~\ref{def:convo}, one has
\begin{equation*}
\left\Vert u\ast v\right\Vert_1 \leq \left(\sup_{m,n\in\N} \sum_{k=0}^{m+n} \vert \alpha^{m,n}_k\vert\right) \left\Vert u\right\Vert_1 \left\Vert v\right\Vert_1,
\end{equation*}
where $\left\Vert u\right\Vert_1 = \sum_{n=0}^\infty \vert u_n\vert$. Therefore, as soon as 
\begin{equation}
\label{eq:banach_alg_cond}
\sup_{m,n\in\N} \sum_{k=0}^{m+n} \vert \alpha^{m,n}_k\vert<\infty,
\end{equation}
the space $\ell^1$ of sequences with finite $\left\Vert \cdot\right\Vert_1$ norm is stable under the convolution product, i.e. $\left(\ell^1,\ast\right)$ is a Banach algebra.

\begin{remark}
\label{rem:banach_alg}
If the family $\phi_n$ is such that $\sup_{s\in[-1,1]} \vert \phi_n(s)\vert \leq 1$ for all $n\in\N$, then the $\ell^1$-norm of the coefficients controls the $\cC^0$-norm of the function:
\begin{equation*}
\sup_{s\in[-1,1]} \vert u(s) \vert \leq \sup_{n\in\N} \vert u_n\vert.
\end{equation*}
In particular, if the coefficients associated to $u$ and $v$ have finite $\ell^1$-norm and~\eqref{eq:banach_alg_cond} is satisfied, then the sum in~\eqref{eq:prod_star} is guaranteed to converge, and~\eqref{eq:prod_star} holds not only formally, but also in $\cC^0$.
\end{remark}

For each family of polynomials used in this work, namely the Legendre polynomials $P_n$, the Chebyshev polynomials of the first kind $T_n$, the canonical basis $X^n$, and (suitable renormalization of) the Chebyshev polynomials of the second kind $\tilde U_n=U_n/U_n(1)$ and the Gegenbauer polynomials $\tilde C^\mu_n = C^\mu_n/C^\mu_n(1)$, we have
\begin{equation}
\label{eq:sumto1}
\sum_{k=0}^{m+n} \vert \alpha^{m,n}_k\vert = 1, \qquad \forall~m,n\in\N.
\end{equation}
Indeed, in all these cases $\phi_n(1)=1$ for all $n\in\N$, therefore~\eqref{eq:lin1} yields
\begin{equation*}
\sum_{k=0}^{m+n} \alpha^{m,n}_k = 1, \qquad \forall~m,n\in\N.
\end{equation*}
Besides, in all these cases the linearization coefficients $\alpha^{m,n}_k$ are nonnegative for all $m,n,k\in\N$, and hence we get~\eqref{eq:sumto1}, which implies~\eqref{eq:banach_alg_cond}. Notice that in all those cases $\sup_{s\in[-1,1]} \vert \phi_n(s)\vert = 1$ for all $n\in\N$, and therefore Remark~\ref{rem:banach_alg} applies. For a more in depth discussion about the convolution products and Banach algebra structures that can be associated to orthogonal polynomials, we refer to the lecture notes~\cite{Szw05} and the references therein.

\medskip

In practice we need to implement the linear operator representing the (truncated) convolution with a given sequence (see e.g. Section~\ref{sec:LV_manifolds}). Given $u=\left(u_n\right)_{0\leq n <N}$, and linearization coefficients $\alpha$ associated to a basis $\left(\phi_n\right)$, we therefore define the $N\times N$ matrix $M_u$ by
\begin{equation*}
\left(M_u\right)_{i,j}=\sum_{k=\vert i-j\vert}^{i+j}\alpha^{k,j}_i u_k \qquad \forall~0\leq i,j <N.
\end{equation*}
It follows from~\eqref{eq:lin1bis} that, for any vector $v=\left(v_n\right)_{0\leq n <N}$
\begin{equation*}
\left(M_u v\right)_k = (u\ast v)_k \qquad\forall~0\leq k<N.
\end{equation*}

\medskip

Let us now consider two basis $\left(\phi^{(1)}_{n_1}\right)_{n_1\in\N}$ and $\left(\phi^{(2)}_{n_2}\right)_{n_2\in\N}$ having linearization coefficients $\alpha^{(1)}$ and $\alpha^{(2)}$ and associated convolution products $\ast^{(1)}$ and $\ast^{(2)}$. We encountered this situation in Section~\ref{sec:LV} and Section~\ref{sec:Lorenz}, for instance when dealing with Taylor$\times$PC expansions, but what follows could also be used for multivariate PC expansions. Given bidimensional sequences $u=\left(u_n\right)_{n\in\N^2}$ and $v=\left(v_{n}\right)_{n\in\N^2}$ associated to expansions
\begin{equation*}
u(x,y)=\sum_{n\in\N^2} u_n \phi^{(1)}_{n_1}(x)\phi^{(2)}_{n_2}(y) \qquad\text{and}\qquad v(x,y)=\sum_{n\in\N^2} v_n \phi^{(1)}_{n_1}(x)\phi^{(2)}_{n_2}(y),
\end{equation*}
we introduce the bidimensional product $\circledast$ such that
\begin{equation*}
u(x,y)v(x,y) =\sum_{n\in\N^2} \left(u\circledast v\right)_n \phi^{(1)}_{n_1}(x)\phi^{(2)}_{n_2}(y).
\end{equation*}
The coefficients $\left(u\circledast v\right)_n$ can of course be computed from the univariate convolution products, for instance by looking at $u$ and $v$ as functions of one variable (say $x$) having coefficients depending on the other variable (say $y$):
\begin{equation*}
u(x,y)=\sum_{n_1\in\N} \left(\sum_{n_2\in\N} u_{n_1,n_2} \phi^{(2)}_{n_2}(y)\right)\phi^{(1)}_{n_1}(x) \qquad\text{and}\qquad v(x,y)=\sum_{n_1\in\N} \left(\sum_{n_2\in\N} v_{n_1,n_2} \phi^{(2)}_{n_2}(y)\right)\phi^{(1)}_{n_1}(x).
\end{equation*}
Denoting $u^{(1)}=\left(u^{(1)}_{n_1}\right)_{n_1\in\N}$, where for all $n\in\N$ $u^{(1)}_{n_1}=\left(u_{n_1,n_2}\right)_{n_2\in\N}$, and slightly abusing the notation $\ast^{(1)}$ we get
\begin{align*}
\left(u\circledast v\right)_k &= \left(\left(u^{(1)} \ast^{(1)} v^{(1)} \right)_{k_1}\right)_{k_2} = \sum_{m_1\in\N}\sum_{n_1\in\N} \left(\alpha^{(1)}\right)^{m_1,n_1}_{k_1} \left(u^{(1)}_{m_1} \ast^{(2)} v^{(1)}_{n_1}\right)_{k_2}.
\end{align*}
Of course the role of each variable/basis is interchangeable and we also have
\begin{align*}
\left(u\circledast v\right)_k &= \left(\left(u^{(2)} \ast^{(2)} v^{(2)} \right)_{k_2}\right)_{k_1} = \sum_{m_2\in\N}\sum_{n_2\in\N} \left(\alpha^{(2)}\right)^{m_2,n_2}_{k_2} \left(u^{(2)}_{m_2} \ast^{(1)} v^{(2)}_{n_2}\right)_{k_1}.
\end{align*}

\subsection{Treatment of higher order terms}
\label{sec:higher_order_terms}

Let us now consider situation where product of three or more expansions have to be computed for a given family $\phi_n$. Of course one could consider the coefficients $\beta^{k,m,n}_l$ such that
\begin{equation*}
\phi_k\phi_m\phi_n = \sum_{l=0}^{k+m+n}\beta^{k,m,n}_l \phi_l,
\end{equation*}
and precompute (and store) them using the formula:
\begin{equation*}
\beta^{k,m,n}_l = \frac{1}{h_l} \left\langle\phi_k\phi_m\phi_n,\phi_l\right\rangle_\rho.
\end{equation*}
However, this approach becomes impracticable very fast even for univariate bases when the degree of the nonlinearity increases. Therefore in practice it is more efficient to stick with only the coefficients $\alpha$ for quadratic products, and compute nonlinear term of higher order recursively. That is, given expansions
\begin{equation*}
u(s)=\sum_{n=0}^\infty u_n \phi_n(s), \qquad  v(s)=\sum_{n=0}^\infty v_n \phi_n(s) \qquad \text{and}\qquad w(s)=\sum_{n=0}^\infty w_n \phi_n(s),
\end{equation*}
we get the coefficients of $u(s)v(s)w(s)$ by first computing $u\ast v$ and then $(u\ast v)\ast w$, which by associativity can simply be denoted $(u\ast v\ast w)$.

\begin{remark}
It should be noted that associativity is often only approximately true in practice because of truncation errors, but this is negligible as soon as we use enough coefficients for the truncation error to be small (see e.g.~\cite[Section 4.5.1.2]{LeMKni10}).
\end{remark}

We point out that, even for systems with nonlinear terms of low order, these considerations are also relevant if one wishes to compute higher order moments associated to PC expansions. Indeed the $k$-th moment of $u$ is given by
\begin{equation*}
\E (u^k) = \sum_{l=0}^\infty (u^k)_l \left\langle\phi_l,1\right\rangle_\rho = (u^k)_0,
\end{equation*}
where we have 
\begin{equation*}
(u^k)_l=(\underbrace{x\ast \ldots\ast x}_{k \text{ times}})_l. 
\end{equation*}
Finally, let us point out that handling non polynomial functions of PC expansions is also possible, although less straightforward, see e.g.~\cite{DebNajPebKniGhaLeM04}.

\subsection{Faster computations}
\label{sec:Appendix_fast}

The usual paradigm of PC is that, to minimize computational cost one should use the expansion basis associated to the PDF of the random inputs, as this minimize the number of coefficients needed to reach a given accuracy (see Section~\ref{sec:LV_eq}). However, the cost of computing the nonlinear terms for a given expansion should also be taken into account, and we briefly discuss it here. Let us consider truncated PC expansions of size $N$:
\begin{equation*}
u(s)=\sum_{n=0}^{N-1}u_n\phi_n \qquad \text{and}\qquad v(s)=\sum_{n=0}^{N-1}v_n\phi_n
\end{equation*}
and evaluate the cost of computing $(u\ast v)_n$ for all $0\leq n<N$, assuming the linearization coefficients $\alpha$ have been precomputed and stored. From~\eqref{eq:lin1bis} we have
\begin{equation*}
(u\ast v)_k = \sum_{n=0}^{N-1} \sum_{m=\vert k-n\vert}^{k+n}\alpha^{m,n}_k u_m v_n,
\end{equation*} 
therefore, using directly this formula the total cost of computing $(u\ast v)_k$ for all $0\leq k<N$ is of order $N^3$. However, this cost can easily be reduced by one order for some specific expansions, for which we know a priori that most of the linearization coefficients are zero. In particular, for Taylor expansions the convolution product (usually called Cauchy product in that case) writes
\begin{equation*}
(u\ast v)_k = \sum_{l=0}^k u_l v_{k-l},
\end{equation*}
and thus the total cost of computing $(u\ast v)_k$ for all $0\leq k<N$ is of order $N^2$. Chebyshev expansions of the first kind also enjoy a similar property. Indeed, let us consider $\phi_n$ defined as $\phi_0=T_0$ and $\phi_n=2T_n$ for all $n\geq 1$. Then the associated convolution product is related to the classical \emph{discrete convolution} product and writes
\begin{equation*}
(u\ast v)_k = \sum_{-N+1}^{N-1} u_{\vert l\vert} v_{\vert k-l\vert},
\end{equation*}
and the total cost of computing $(u\ast v)_k$ for all $0\leq k<N$ is again of order $N^2$. Finally, when $N$ becomes large one may want to try and reduce this cost even further. A natural way to do so for Chebyshev or Taylor expansions is to compute the convolution via a Fast Fourier Transform, which then brings down the cost even lower, to an order of $N \log N$  (see e.g.~\cite{Nus12,RaoYip14}). This idea was partially extended to other bases such as the Gegenbauer polynomials, see~\cite{PotSteTas98} and the references therein for more details.

\paragraph*{Acknowledgments:} MB and CK have been supported by a Lichtenberg Professorship of the VolkswagenStiftung.

\bibliographystyle{plain}
\bibliography{bibfile}

\end{document}